\documentclass[11pt]{article}

\usepackage{booktabs}
\usepackage[margin=15pt,font=small,labelfont={bf,sf},justification=justified]{caption}
\usepackage[superscript,biblabel]{cite}
\usepackage{color}
\usepackage{float}
\usepackage{graphicx}
\usepackage{lscape}
\usepackage[bbgreekl]{mathbbol}
\usepackage{mathtools}
\usepackage{multibib}
\usepackage{multirow}
\usepackage{tabularx}
\usepackage{titleref}
\usepackage{url}

\usepackage{bm}
\usepackage{hyperref}
\usepackage{subcaption}

\newif\ifbembo
\newif\ifcharter
\newif\iferewhon
\newif\iflibertine
\newif\iflibertinealt
\newif\ifpalantino
\newif\iftimesnewroman

\bembofalse
\chartertrue
\erewhonfalse
\libertinefalse
\libertinealtfalse
\palantinofalse
\timesnewromanfalse
\ifbembo
  \usepackage[p,osf]{ETbb} %
  \usepackage[scaled=.95,type1]{cabin} %
  \usepackage[varqu,varl]{zi4}%
  \usepackage[T1]{fontenc}
  \usepackage[libertine,vvarbb]{newtxmath}
  \usepackage[cal=boondoxo,bb=boondox]{mathalfa}
\fi

\ifcharter
  \usepackage[scaled=.98,p]{XCharter}
  \usepackage[scaled=1.04,varqu,varl]{inconsolata}
  \usepackage[type1]{cabin}
  \usepackage[uprightscript,xcharter,vvarbb,scaled=1.05]{newtxmath}
  \usepackage[cal=euler,scr=boondoxo,bb=boondox]{mathalpha}  %
\fi

\iferewhon
  \usepackage[p,osf,scaled=.98,space]{erewhon} %
  \usepackage[varqu,varl]{inconsolata} %
  \usepackage[type1,scaled=.95]{cabin} %
  \usepackage[utopia,vvarbb]{newtxmath}
\fi

\iflibertine
  \usepackage[sb]{libertinus}
  \usepackage[T1]{fontenc}
  \usepackage{textcomp}
  \usepackage[varqu,varl]{zi4}%
  \usepackage{libertinust1math} %
  \usepackage[scr=boondoxo,bb=boondox]{mathalpha} %
\fi

\iflibertinealt
  \usepackage[sb]{libertinus}
  \usepackage[T1]{fontenc}
  \usepackage{textcomp}
  \usepackage{libertinust1math} %
  \usepackage[scr=boondoxo,bb=boondox]{mathalpha} %
\fi

\ifpalantino
  \usepackage[largesc]{newpxtext}
  \usepackage[vvarbb]{newpxmath}
\fi

\iftimesnewroman
  \usepackage{newtxtext} %
  \usepackage[scaled=0.95]{inconsolata} %
  \usepackage[vvarbb]{newtxmath} %
\fi

\usepackage[T1]{fontenc}
\usepackage[final,protrusion=true,expansion=true]{microtype}

\usepackage{titling}
\usepackage{authblk} %
\pretitle{\begin{center}\bfseries\sffamily\LARGE}
\posttitle{\end{center}}
\usepackage{abstract}
\usepackage{sectsty} %
\allsectionsfont{\sffamily}

\usepackage[left=1in,right=1in,top=1in,bottom=1in,includefoot,heightrounded]{geometry}

\newcites{supp}{SUPPLEMENTAL REFERENCES}

\makeatletter
\patchcmd{\LS@rot}{90}{-90}{}{}
\patchcmd{\endlandscape}{90}{-90}{}{}
\makeatother

\newcommand*{\tran}{^{\mkern-1.5mu\mathsf{T}}}

\global\long\def\grad{\boldsymbol{\nabla}}

\global\long\def\J{\sM J}
\global\long\def\S{\sM S}

\renewcommand{\vec}[1]{\bm{\mathrm{#1}}}

\def \x{\vec{x}}
\def \xp{\x'}
\def \X{\vec{X}}
\def \Xn{\X_n}
\def \Vn{V_n}
\def \Xl{\X_m}

\def \Xp{\X'}

\def\Fubar{\vec{F}_{\text{PD}}}

\def \Yubar{\underbar{\vec{Y}}}

\def \Tubar{\underbar{\vec{T}}}

\def \horizon{\vec{\mathcal{H}_X}}
\def \dhorizon{\vec{\mathcal{H}}}
\def \PP{\vec{\mathbb{P}}}
\def \FF{\vec{\mathbb{F}}}

\def \BB{\vec{\mathbb{K}}}
\def \Xi{\vec{\xi}}
\def \Eta{\vec{\eta}}

\def \frho{\rho}
\def \fmu{\mu}

\def \u{\vec{u}}
\def \U{\vec{U}}

\def \v{\vec{v}}

\def \f{\vec{f}}
\def \F{\vec{F}}

\def \Omegaf{\Omega^\text{f}}
\def \Omegas{\Omega^\text{s}}
\def \vchi{\vec{\chi}}

\def \N{\boldsymbol{\mathcal{N}}}
\def \S{\boldsymbol{\mathcal{S}}}
\def \J{\boldsymbol{\mathcal{J}}}

\def \horizonsize{\epsilon}

\def \sca{s_{{\text{c}}_1}}
\def \scb{s_{{\text{c}}_2}}

\def \dVp{ \mathrm{d} {\Xp}}

\def \omegah{\hat{\omega}}

\def \a{\vec{A}}

\title{An immersed peridynamics method for fluid-driven damage and failure of anisotropic materials}
\author[1,2]{Keon Ho Kim}
\author[2--7]{Boyce E. Griffith}
\affil[1]{Oden Institute for Computational Engineering and Sciences, University of Texas, Austin, TX, USA}
\affil[2]{Department of Mathematics, University of North Carolina, Chapel Hill, NC, USA}
\affil[3]{Department of Applied Physical Sciences, University of North Carolina, Chapel Hill, NC, USA}
\affil[4]{Department of Biomedical Engineering, University of North Carolina, Chapel Hill, NC, USA}
\affil[5]{Carolina Center for Interdisciplinary Applied Mathematics, University of North Carolina,
Chapel Hill, NC, USA}
\affil[6]{Computational Medicine Program, University of North Carolina, Chapel Hill, NC, USA}
\affil[7]{McAllister Heart Institute, University of North Carolina, Chapel Hill, NC, USA}
\affil[ ]{\href{mailto:kkeonho@utexas.edu}{kkeonho@utexas.edu} and \href{mailto:boyceg@email.unc.edu}{boyceg@email.unc.edu}}

\begin{document}
\maketitle

\begin{abstract}
\noindent The immersed peridynamics (IPD) method is a fluid-structure interaction (FSI) model to simulate fluid-driven material damage and failure of an immersed structure, in which a peridynamic (PD) constitutive correspondence model is employed within a classical immersed boundary (IB)-type framework to describe stresses, forces, and structural deformations of a structural body, instead of classical continuum mechanics.
This paper introduces an extension of the IPD method to simulate fluid-driven structural deformation, damage, and failure of anisotropic materials with complex geometries.
Our prior IPD work for isotropic materials demonstrated grid-convergent fluid-driven failure processes in simple geometries, using uniform lattices to describe the discretized immersed structures. 
However, this uniform discretization approach limits the fidelity of the methodology for complex structural bodies due to resulting stair-stepped artifacts.
To ensure more accurate geometric representations, we use quadrature rules attached to finite element (FE) meshes to generate both the PD points and their associated weights, which are used to approximate the PD integrals.
We demonstrate that non-uniform discretizations improve both accuracy and volume conservation of hyperelastic materials along with accurately represented boundaries.
To capture realistic biomaterial behaviors, we incorporate hyperelastic constitutive models including both isotropy and anisotropy into the proposed IPD method.
In addition, a ductile failure model is adopted to simulate progressive damage and failure in anisotropic materials.
For non-failure cases, the accuracy and convergence of our approach are compared to an FE-based IB method using widely used benchmark problems of nonlinear incompressible elasticity including anisotropy. 
Our numerical simulations demonstrate that the extended IPD method yields comparable accuracy with similar numbers of structural degrees of freedom for different choices of peridynamic horizon sizes. 
For failure tests, we investigate how fiber orientation affects deformation, damage growth, and crack propagation using hyperelastic material models including anisotropy. 
We further demonstrate that the developed method generates grid-converged simulations of damage growth, crack formation and propagation, and rupture under large deformations, including purely fluid-driven failure processes.

\end{abstract}


\section{Introduction}
Fluid-driven material damage and failure have been widely studied in the context of many engineering and industrial fields, including tissue failure~\cite{MacMinn:2015tk, o2018experimental, wang2022fluid}, hydraulic fracturing~\cite{OSIPTSOV2017513, detournay2016mechanics}, blast loading on structures~\cite{Kambouchev_2006, gupta2021response, aune2021influence}, and drainage fracture~\cite{Kobchenko_2013, lee2016pressure}. 
Understanding the mechanisms behind such fluid-driven failure processes is crucial for predicting the safety and efficacy of biomedical devices and industrial systems. However, accurately simulating these phenomena remains a significant computational challenge due to the complex coupling between fluid dynamics and solid mechanics along with accurately modeling structural damage growth and fracture.

Advanced computational fluid-structure interaction (FSI) approaches, such as the immersed boundary (IB) method~\cite{peskin_2002} and arbitrary Lagrangian-Eulerian (ALE) formulations~\cite{TAKASHI1992115}, have achieved significant success in a broad range of applications, including cardiovascular dynamics~\cite{PESKIN1989372, MCQUEEN1989289, McQueen_2000,  McQueen_book, griffith_2009, griffith2012immersed, Lee_2020, LEE202160, Choi_2014, Bailoor_2021, Davey_2024}, gastrointestinal systems~\cite{Kou_2015, Kou_2017, Lee_2022, Kuhar_2022}, and biolocomotion~\cite{Bhalla_2013, BHALLA2013446, HERSCHLAG201184, Koumoutsako_2008, BORAZJANI20087587, Borazjani_2008, fluids3030045, Wang_2005}.
These standard continuum-based frameworks are robust in simulating large structural deformations and moving boundaries in continuous media, but they present substantial mathematical and algorithmic challenges when modeling topological changes, such as fluid-driven rupture, degradation, and crack propagation in structural bodies.
Their governing equations are primarily based on partial differential equations and continuum mechanics (CM), which fundamentally assume spatial continuity.
Consequently, the local definitions of stresses and strains used to describe the dynamics of deformations within these FSI frameworks become ill-conditioned and cause difficulties in simulating structural mechanics near crack fronts due to discontinuities in the displacement field.

Simulating such topological changes (particularly fracture and rupture) within traditional CM-based methods often requires computationally expensive and algorithmically complex treatments, such as remeshing~\cite{LEE199499} or element deletion~\cite{Song2008}, which can compromise numerical stability.
To address such challenges, advanced numerical treatments have been integrated into CM-based frameworks. 
The cohesive zone model (CZM)~\cite{hillerborg1976analysis} explicitly models fracture by inserting cohesive elements along material interfaces to govern surface separation. 
To overcome mesh dependency, the eXtended finite element method (XFEM)~\cite{moes1999finite} was introduced, which enables the representation of arbitrary discontinuities within finite elements through enriched approximation spaces. 
More recently, phase-field (PF) modeling~\cite{Bourdin_2008} has gained traction for its ability to predict complex crack topologies by employing a diffuse interface approach that regularizes sharp discontinuities over a finite length scale. 
Although many numerical methods have improved the capabilities of continuum mechanics for complex fracture simulations, integrating them into FSI simulations remains computationally demanding.
For instance, explicit interface tracking in CZM and XFEM often complicates the coupling with fluid solvers during topological changes, whereas PF methods require prohibitively fine meshes to resolve the damage gradient~\cite{Pezzulli_2025}.

As an alternative to CM-based methods for fracture mechanics, a non-local extension of CM called peridynamics (PD) was introduced by Silling~\cite{silling2000reformulation} to overcome the limitations in the CM-based formulations.
Unlike classical approaches, PD models avoid the use of derivatives in determining strains.
Instead, the governing equations are integral formulations that consider interactions among all material points within a finite range called the PD horizon. 
This mathematical structure remains valid even in the presence of discontinuities in the displacement field. 
As a result, PD naturally simulates spontaneous crack initiation and propagation through the local breakage of interactions (i.e., bonds) between particles, removing the need for complex tracking algorithms such as those required in XFEM or remeshing strategies~\cite{silling2005meshfree}.

Recently, several attempts have been made to develop coupled PD-FSI models for fluid-driven failure processes. 
Fully Lagrangian frameworks, such as coupling PD with smoothed particle hydrodynamics~\cite{FAN2017349, sun2020smoothed}, have been widely adopted for modeling FSI scenarios involving fragmentation, such as blast loading or high-velocity impact. 
Approaches utilizing the lattice Boltzmann method have also been proposed, particularly for simulating flow through fractured porous media~\cite{zhang2021strongly, zhang2026coupled}. 
Furthermore, to leverage the robustness of mature fluid solvers, recent efforts have integrated PD with ALE formulations~\cite{behzadinasab2021iga} or IB methods~\cite{dalla20223d, patel2025fluid}.

Our prior effort on such a hybrid PD-FSI model for fluid-driven material damage and failure of incompressible hyperelastic materials is called the immersed peridynamics (IPD) method~\cite{KIM2023112466, kim2023immersed}, in which we developed an IB-type FSI model integrated with non-ordinary state-based peridynamics (NOSB-PD)~\cite{silling2007peridynamic, warren2009non}.
The IB method uses Lagrangian variables for the deformations, stresses, and resultant forces of the immersed structure and Eulerian variables for the momentum, viscosity, and incompressibility of the coupled fluid-structure system. 
Coupling between Lagrangian and Eulerian variables is mediated by integral transforms with Dirac delta function kernels in the continuous formulation. 
In discretized IB formulations, the singular delta function is replaced by a regularized delta function~\cite{peskin_2002}. 
This coupling approach enables nonconforming discretizations of the fluid and immersed structures~\cite{griffith2020immersed,griffith2017hybrid}. 
Conventional IB methods use the framework of nonlinear continuum mechanics to compute the elastic body forces of the immersed structure~\cite{griffith2017hybrid, BOFFI20082210, ZHANG20042051, WELLS2023111890}, whereas IPD focuses on the constitutive correspondence model of NOSB-PD to compute structural forces instead of continuum mechanics.
The force vectors in NOSB-PD can differ both in magnitude and direction, which is far more general than other PD theories and is helpful for developing PD models for realistic hyperelastic materials. 
In addition, the constitutive correspondence model has the advantage that we can leverage existing well-developed material models in classical continuum theories instead of defining new material characteristics only for PD models.
For non-failure tests, the IPD method yields comparable accuracy under grid refinement to that offered by the immersed finite element-finite difference (IFED) method~\cite{griffith2017hybrid}, and the developed numerical methodology demonstrates convergent
and consistent failure predictions for a nontrivial range of grid spacings.

In this work, we extend the initial version of the IPD method~\cite{KIM2023112466} that has the potential capability of simulating realistic fluid-driven fracture mechanics of biological materials, such as those observed during aortic dissection~\cite{Nienaber_2016, bonfanti2020combined}. 
For simplicity in numerical implementations, our prior IPD method was limited to isotropic materials and uniformly distributed nodal volumes within the immersed body. 
However, such a uniform distribution can cause unrealistic representations of complex material geometries for realistic biomaterials, i.e., stair-step artifacts.
To accurately describe complex structural bodies, we introduce the non-uniform distribution of PD points, and such a non-uniformly distributed point cloud requires an accurate volume computation for PD volumetric forces.
We use a quadrature rule with a finite element (FE) mesh to generate the PD points as well as their associated weights, which are needed to accurately evaluate the PD integrals in the IPD formulation with non-uniform discretizations.
In particular, this approach ensures seamless geometric matching between the meshfree structural model and the FE-based coupling scheme, thereby optimizing computational efficiency.

Another key contribution of this work is the investigation of the large-deformation and failure behaviors of anisotropic hyperelastic materials using NOSB-PD.
Early bond-based PD models represented fiber-reinforced composites through constituent- and orientation-dependent bond properties~\cite{Hu2012}, whereas state-based PD formulations enabled general anisotropic constitutive behavior to be introduced at material points~\cite{Hattori2018}.
Other studies incorporated explicitly anisotropic failure criteria through direction-dependent critical stretches or fracture properties~\cite{Zhang2019,Fang2021,Li2024}.
More recent developments include general ordinary and non-ordinary state-based PD formulations~\cite{Scabbia2024,Wang2024}, as well as finite-deformation models for transversely isotropic and reinforced hyperelastic materials~\cite{Xu2018,deVilliers2025,YIN2025117494}.
Despite these advances, most existing studies focus on composite materials in uncoupled solid-mechanics settings, and applications of PD to the large deformation and failure of anisotropic hyperelastic materials remain limited~\cite{Sienkiewicz2026}.
This study provides insights into fracture mechanics in anisotropic materials by extending the IPD framework to capture the complex nonlinear fracture behavior inherent in biomaterials.

To simulate anisotropic behaviors of immersed structures in the IPD framework, we employ a well-known constitutive model for anisotropic materials for fibrous soft tissues~\cite{Gasser_2005}.
The benchmark problems, including the compression test~\cite{reese1999new}, three-dimensional Cook's membrane~\cite{Wriggers_2016}, and three-dimensional torsion test~\cite{bonet2015computational}, are reinforced with fibers in the immersed structures, and we demonstrate that the proposed method yields comparable accuracy to a FE-based IB method~\cite{griffith2017hybrid} for anisotropic material responses.
In addition, a realistic adventitial tissue model~\cite{Gasser_2005} is adopted to validate the accuracy of the deformation obtained by the proposed IPD method.
We also investigate the effects of fiber orientations on crack formation and propagation by using a two-dimensional fiber-reinforced tissue sheet with a ductile failure model, which allows progressive degradation rather than the abrupt brittle failure considered in our previous isotropic IPD formulation.
For purely fluid-driven fracture mechanics, we adopt an elastic band benchmark~\cite{vadala2020stabilization,KIM2023112466} and demonstrate grid-converged simulations of damage growth, crack formation and propagation, and rupture under large deformations.

\section{Continuous formulations}
This section presents the continuous IPD formulation describing FSI along with material damage and failure.
We first discuss the continuous FSI formulation based on the IB framework and introduce the constitutive correspondence model for NOSB-PD used for immersed structural mechanics along with a specific material model used for realistic biomaterial responses.

\subsection{Immersed peridynamics}\label{s:continuous IPD}
\begin{figure}[t!]
	\centering
	\includegraphics[width = 0.8\textwidth]{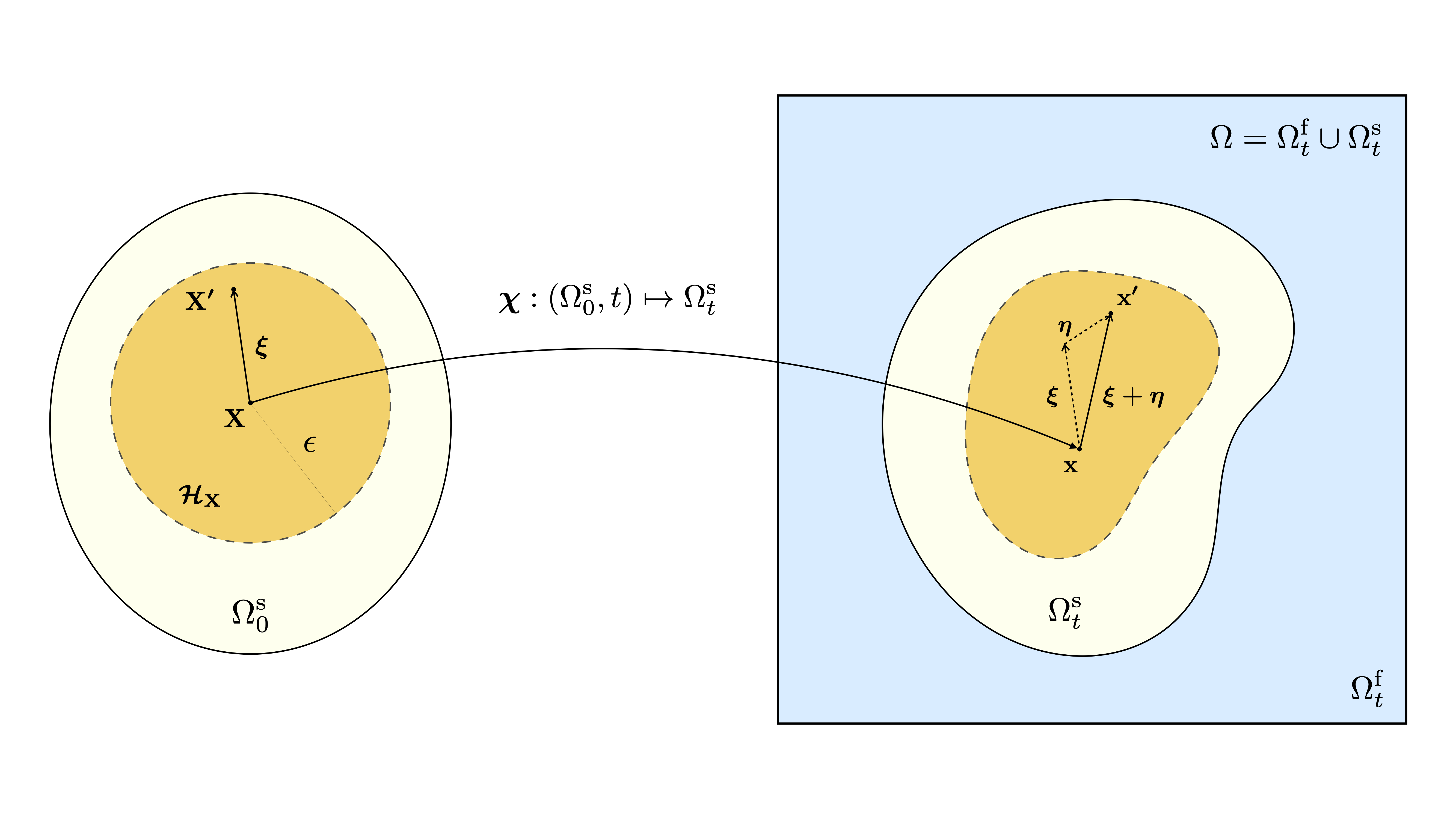}
	\caption{Schematics of the computational domain $\Omega$ and the time-dependent Lagrangian and Eulerian coordinate systems in the IPD formulation along with PD material points, bonds, and horizons. The PD node $\X \in \Omega_0^{\text{s}}$ interacts with its neighbors within a finite range called the horizon denoted by $\horizon \subset \Omega_0^{\text{s}}$.}
	\label{f:IPD_schematics}
\end{figure}

We consider a computational domain $\Omega$ subdivided into a time-dependent fluid subdomain $ \Omega_t^{\text{f}}$ and a solid subdomain $\Omega_t^{\text{s}} $ that are indexed by time $t$. 
We use both fixed spatial coordinates $\x \in \Omega$ and reference coordinates $\X \in\Omega_0^{\text{s}}$, with $\Omega_0^{\text{s}}$ indicating the region occupied by the solid structure at time $t = 0$. 
The deformation of an immersed structural body occupying the region $\Omegas_t$ at time $t$ is described by the reference coordinate $\X \in \Omegas_0$.
We use a deformation mapping $\vchi : \left( \Omegas_0,t \right) \mapsto \Omegas_t$ to relate the reference and deformed coordinate systems, so that $\x = \vchi\left( \X,t \right)$ is the physical position of the material point $\X$ at time $t$.
See Fig.~\ref{f:IPD_schematics}.

The dynamics of the coupled fluid-structure system are described by
\begin{align}
\frho \frac{\mathrm{D} \u}{\mathrm{D} t} \left( \x,t \right)  &= - \grad p\left( \x,t \right) + \fmu \grad^2 \u\left( \x,t \right) + \f \left( \x,t \right), &\x \in \Omega, \label{e:navier_stokes} \\
\grad \cdot \u \left( \x,t \right) & = 0, &\x \in \Omega, \label{e:incompressibility} \\
\f\left( \x,t \right) &= \int_{\Omega_0^{\text{s}}} \F \left( \X,t \right)\, \delta \left(\x - \vchi\left( \X,t \right)\right) \mathrm{d}\X, &\x \in \Omega,  \label{e:force_spreading}\\
\frac{\partial \vchi}{\partial t} \left( \X,t \right) &= \U \left( \X,t \right) =  \int_{\Omega}  \u (\x,t)\, \delta \left(\x - \vchi\left( \X,t \right)\right) \mathrm{d}\x, &\X \in \Omegas_0, \label{e:velocity}
\end{align}
in which $\u \left( \x,t \right)$ is the Eulerian velocity, $\U \left( \X,t \right)$ is the velocity of the structure, $\frho$ is the uniform mass density of both the fluid and the structure, $\fmu$ is the uniform dynamic viscosity, $p \left( \x,t \right)$ is the physical pressure, $\f \left( \x,t \right)$ is the Eulerian structural force density generated by the deformations of the structure, and $\F \left( \X,t \right)$ is the Lagrangian structural force density. 
The operators $\grad$, $\grad^2$, and $\grad \cdot$ are with respect to Eulerian coordinates, and $\frac{\mathrm{D}}{\mathrm{D}t} = \frac{\partial}{\partial t} + \u \cdot \grad$ is the material time derivative in current coordinates.
Eqs.~\eqref{e:incompressibility} and \eqref{e:velocity} imply the exact incompressibility of the immersed structure, and 
Eq.~\eqref{e:velocity} implies that the no-slip boundary condition holds along the fluid-solid interface. 
The coupling between fluid and structure is achieved by interaction equations with the Dirac delta function $\delta$ as in Eqs.~\eqref{e:force_spreading} and \eqref{e:velocity}.

In the IPD formulation, both stationary and flexible bodies are considered.
For stationary structures, we approximate the rigidity constraints $\frac{\partial \vchi}{\partial t} \left( \X,t \right) = \vec{0}$ by using a Lagrange multiplier force:
\begin{align}\label{Lagrangian_multiplier}
\F_{\text{c}} \left( \X,t \right) = \kappa \left( \X - \vchi\left( \X,t \right) \right) - \eta \U \left( \X,t \right),
\end{align}
in which $\kappa$ is a stiffness penalty parameter, $\eta$ is a damping penalty parameter, and $\F_{\text{c}}$ has units of force per unit volume.
Note that as $\kappa \rightarrow \infty$, $\vchi \left( \X,t \right) \rightarrow \X$, and $\frac{\partial \vchi}{\partial t} \left( \X,t \right) \rightarrow \vec{0}$, and we recover an exactly stationary model.
For flexible structures, the Cauchy stress on the computational domain is defined as
\begin{align}\label{e:cauchy stress}
\bbsigma \left( \x, t \right)  = \bbsigma^{\text{f}} \left( \x, t \right)+ 
\begin{cases}
\bbsigma^{\text{s}}\left( \x, t \right), \, &\text{if } \x \in \Omegas_t, \\
\mathbb{0}, \,  &\text{if } \x \in \Omegaf_t,
\end{cases}
\end{align}
in which $\bbsigma^{\text{f}} = -p \, \mathbb{I} + \mu \left(\grad \u + \grad \u^{\tran} \right)$ is the fluid Cauchy stress and $\bbsigma^{\text{s}}$ is the structural Cauchy stress determined by the constitutive correspondence model~\cite{silling2007peridynamic}.
Additionally, the damping term in Eq.~\eqref{Lagrangian_multiplier} is used for a flexible structure to assist the system in reaching steady state.

\subsection{Immersed body force}
\label{s:immersed body forces}
The evaluation of structural forces within the IPD method relies on the PD framework~\cite{silling2000reformulation, silling2005meshfree, silling2007peridynamic}, which provides a non-local extension of classical continuum mechanics.
Among various choices of PD models, we particularly focus on the constitutive correspondence model of NOSB-PD~\cite{silling2007peridynamic}.
This specific constitutive correspondence model is highly advantageous as it allows us to directly use well-established material characteristics in classical continuum mechanics.

Unlike interacting only with adjacent neighbors in the classical continuum theory, a PD material point $\X$ interacts with its neighbors within a finite region around $\X$ called the horizon and denoted as $\horizon \subset \Omegas_0$, whose size provides a length-scale parameter missing in classical continuum mechanics.
Although various geometric definitions of the horizon exist in the PD literature, we particularly focus on a spherical horizon centered at $\X$ with a radius $\horizonsize >0$\footnote{In the peridynamics literature, the horizon size is commonly denoted by $\delta$, but here we use $\horizonsize$ to distinguish the horizon size from the Dirac delta function used in Sec.~\ref{s:continuous IPD}.}.
The PD node $\X$ interacts with its neighbor $\Xp$ in the horizon $\horizon$ through a bond $\Xi = \Xp - \X$ in the reference frame, and the bond deforms to $\Xi + \Eta = \xp - \x$ in the deformed frame, with $\xp = \vchi \left(\Xp,t \right)$ and $\x = \vchi\left( \X,t \right)$.
See Fig.~\ref{f:IPD_schematics}.

To accurately track the kinematic evolution of a material body, state-based peridynamic theories employ the concept of peridynamic states. 
The deformation mapping $\vchi$ of the initial bond  $\Xi$ at point $\X$ to its deformed configuration at time $t$ is captured by the deformation vector state $\Yubar$:
\begin{align}\label{deformation_vec_state}
    \Yubar = \Yubar\left[ \X,t \right] \langle \Xi \rangle = \xp - \x = \vchi \left( \X + \Xi,t \right) - \vchi\left( \X,t \right).
\end{align}
Correspondingly, the internal force generated within this specific bond $\Xi$ associated with material point $\X$ at time $t$ is represented by the force vector state:
\begin{align}
\Tubar = \Tubar\left[ \X,t \right] \langle\Xi\rangle.
\end{align}

The NOSB-PD correspondence model defines a non-local deformation gradient tensor to describe material deformation in a non-local sense around $\X$. 
The non-local deformation gradient tensor serves as the non-local equivalent to the classical local deformation gradient tensor and is defined by
\begin{equation}\label{nonlocal_deformation_gradient}
    \FF  = \left[ \int_{\horizon} \omega \left(| \Xi |\right) \; \left( \Yubar\langle \Xi \rangle \otimes \Xi \right) \; \dVp  \right] \BB^{-1},
\end{equation}	
in which $\omega$ is a non-negative scalar-valued influence function that controls the relative contribution of interacting particles within the horizon, and $\BB$ is the shape tensor defined by
\begin{equation}\label{shape_tensor}
	\BB =  \int_{\horizon} \omega \left(| \Xi |\right) \; \left( \Xi \otimes \Xi \right) \; \dVp.
\end{equation}
Fundamentally, $\FF$ plays a similar role to the classical local deformation gradient, which describes the mapping of the reference PD horizon to its current horizon, thereby providing an averaged deformation at the material point $\X$ through interactions with all neighbors in the PD horizon.
Similarly, the shape tensor represents a weighted volume average of the reference bonds within the horizon, which effectively normalizes this nonlocal integration.
Because these kinematic descriptors rely entirely on integral operators rather than spatial derivatives, they remain mathematically robust and well-posed even when the displacement field exhibits sharp discontinuities in the structural body during the crack propagation.

Under the constitutive correspondence principle, individual non-local bond forces are linked to macroscopic stress measures.
The force vector in the constitutive correspondence model is defined by
\begin{align}\label{general_force_vector}
    \Tubar \langle\Xi\rangle = \omega \left(| \Xi |\right) \PP \BB^{-1} \Xi,
\end{align}
in which $\omega$ is a non-negative scalar-valued function called the influence function, which controls the influence of interacting bonds in the horizon, and $\PP$ is the first Piola-Kirchhoff structural stress tensor of the immersed body.
We evaluate $\PP$ by substituting the non-local deformation gradient tensor into a standard strain energy functional $\Psi$, such that $\PP = \frac{\partial \Psi(\FF)}{\partial \FF}$. 
Note that this specific form of the force vector state ensures the conservation of both linear and angular momentum~\cite{silling2007peridynamic}.
Subsequently, the total interaction between two material points $\X$ and $\Xp$ with a bond $\Xi$ is given by the pairwise bond force function $\Fubar(\Xp,\X,t)$:
\begin{align}
\Fubar \left( \Xp,\X,t \right) &  = \Tubar \left[ \X,t \right] \langle\Xi\rangle - \Tubar\left[ \Xp,t \right] \langle-\Xi\rangle.
 \label{eqn_f_pd}            
\end{align}
Consequently, to determine the macroscopic internal force density at the node $\X$, we integrate the pairwise interactions across the entire horizon $\horizon$:
\begin{align}\label{net_pd_force}
\F \left(\X , t \right) =  \int_{\horizon} \Fubar \left( \Xp,\X,t \right)  \; \dVp.
\end{align}
This net PD body force density serves as the immersed body force density in the IPD formulation.

\subsection{Failure model}
\label{s:failure model}
Peridynamic models represent material damage and failure at the bond level, whereby an interacting bond between two particles can fail under excessive deformation.
Such bond breakage is irreversible, and there are various failure criteria, such as a critical stretch-based criterion, a strain-based criterion, or an energy-based criterion.
In our IPD model, we use a critical stretch-based criterion.
Since we focus exclusively on homogeneous materials, the same spatially uniform critical stretch values are assigned to all bonds for both isotropic and anisotropic materials. 
For anisotropic materials, the directional failure response is implicitly represented through the constitutive response and corresponding bond-stretch distribution rather than from explicitly orientation-dependent critical stretch values.

In general, there are two types of failure behavior: brittle and ductile.
In the brittle failure model, if the bond stretch exceeds its critical value, there is no longer interaction between the two material points connected by the bond; the bond breaks immediately. 
In contrast, in the ductile failure model, the bond gradually transitions from the active state to the inactive state if it exceeds the critical stretch and eventually breaks, i.e., softening.
To represent progressive bond degradation rather than abrupt brittle rupture, we employ a ductile failure model with a finite softening regime~\cite{Li2016, Kamani_2024}.

To account for irreversible damage growth, we define a maximum historical bond stretch as
\begin{align}\label{eq:max_histo_stretch}
s_{\max}\left[\X,t\right]\langle \Xi\rangle = \max_{0\le\tau\le t} \frac{|\Yubar\left[ \X,\tau \right] \langle \Xi \rangle|}{|\Xi|}.
\end{align}
Then a time-dependent bond connectivity is tracked by an indicator function $I \left[\X,t\right]\langle \Xi\rangle$ as follows:
\begin{align}\label{fail_paramter}
I \left[\X,t\right]\langle \Xi\rangle = 
\begin{cases} 
1, &\text{ if } s_{\max} \le \sca, \\
\frac{\scb - s_{\max}}{\scb - \sca}, &\text{ if } \sca < s_{\max} \le \scb, \\
0, &\text{ if } \scb < s_{\max},
\end{cases}
\end{align}
in which $\sca$ and $\scb$ are critical values.
Instead of a binary cutoff in the brittle failure model, we define an intermediate softening regime bounded by a lower threshold (i.e., damage initiation) and an upper limit (i.e., complete rupture).
A smaller softening width produces a more abrupt, brittle-like loss of bond interaction, whereas a larger width produces more gradual bond degradation.

In the IPD formulation, this continuous degradation parameter directly modulates the spatial interaction weights and yields a modified influence function:
\begin{align}\label{modified_inf_func}
\omegah \left(|\Xi|,t\right) =  \omega \left(|\Xi| \right) \, I \left[\X,t\right]\langle \Xi\rangle ,
\end{align} 
which controls both the connectivity and the influence of an interacting particle via the bond $\Xi$ during the deformation.
For notational simplicity, the dependence of $\omegah$ on the material point $\X$ and the bond-specific damage history through $I \left[\X,t\right]\langle \Xi\rangle$ is left implicit.
Consistent with established damage-correspondence formulations in the PD literature~\cite{queiruga2017numerical,behera2020peridynamic}, we seamlessly model localized damage and failure by modifying the non-local deformation gradient, shape tensor, and force vector state in Eqs.~\eqref{nonlocal_deformation_gradient}--\eqref{general_force_vector} through the substitution of the original influence function with its bond-status-dependent counterpart.
Consequently, this failure model changes not only the force carried by the damaged bond but also the non-local kinematic reconstruction of the remaining active neighborhood.

To quantify the local damage at a PD node $\X$ at time $t$, we evaluate a scalar damage field $\varphi(\X,t)$ by tracking the bond connectivity within the peridynamic horizon $\horizon$ as described by Behera et al.~\cite{behera2020peridynamic}:
\begin{align}\label{damage}
\varphi \left( \X, t \right) = 1 - \frac{\int_{\horizon} I \left[\X,t\right]\langle \Xi\rangle \, \dVp}{\int_{\horizon} \dVp},
\end{align}
which measures the volume-weighted loss of bond connectivity at time $t$.
This volume-weighted ratio dynamically maps the physical degradation of the structural body, and $\varphi = 0$ implies a non-failure state, whereas $\varphi = 1$ indicates the total topological separation.

\subsection{Constitutive model}
\label{s:constitutive_laws}
Realistic biomaterials possess complex characteristics such as anisotropy that can be observed at multiple scales, ranging from macro- to microscopic levels. 
For instance, in arterial walls, the collagen fibers are aligned in a highly ordered manner, providing stiffness and strength against external loadings.
Such anisotropy in fiber-reinforced structures can be accounted for by incorporating anisotropic energies into the constitutive models. 
To examine deformations of anisotropic materials in the IPD framework, we consider several closely related anisotropic material models, including the well-known fiber-reinforced material model introduced by Holzapfel, Gasser, and Ogden (HGO)~\cite{Gasser_2005}. 
The HGO model includes the effect of the anisotropy through the anisotropic invariant in an incompressible neo-Hookean material model.

For the ground isotropic material, we adopt a modified neo-Hookean model~\cite{vadala2020stabilization} unless otherwise mentioned, which was previously shown to model incompressible isotropic hyperelastic materials in the IPD framework~\cite{KIM2023112466}.
In the modified neo-Hookean model, we use the Flory decomposition of the deformation gradient tensor~\cite{Flory_1961} and define a modified non-local deformation gradient tensor $\bar{\FF} = J^{-1/d}\, \FF$, in which $J = \det \left( \FF \, \right)$ and $d$ is the spatial dimension. 
The strain energy and elastic stress consist of two parts, isochoric and volumetric,
\begin{align}\label{neo_hookean_strain}
\Psi_{\text{isotropic}} &= \Psi_{\text{isochoric}} + \Psi_{\text{volumetric}},\\
\Psi_{\text{isochoric}} &= \frac{G}{2}\left( \bar{I}_1 - 3 \right), \ d=2,3, \\
\Psi_{\text{volumetric}} &= \frac{K}{2} \left( \ln J \right)^2,
\end{align}
in which $G$ is the shear modulus, $\bar{I}_1 = \text{tr}\left(\mathbb{\bar{C}}\right) $ is the modified first invariant, $\mathbb{\bar{C}} = \bar{\FF}^{\tran} \bar{\FF}$ is the modified right Cauchy-Green tensor, $K$ is the bulk modulus, and $d$ is the spatial dimension.
Note that we replace the bulk modulus with the numerical bulk modulus introduced by Vadala-Roth et al.~\cite{vadala2020stabilization} to achieve better volume conservation in the discretized formulation.

To impose the characteristics of biomaterials, we follow the HGO model and consider two transversely isotropic energy contributions of the families of fibers added to the soft ground material:
\begin{align}\label{HGO_model}
\Psi_{\text{anisotropic},i} &= \frac{k_1}{2k_2} \left[ \exp\left( k_2 \left( \left(\kappa \, \bar{I}_1 + (1-3\kappa) \, I_{4,i}\right) - 1\right)^2 \right) -1 \right], \ i=1,2,
\end{align}
in which $k_1$ is a stress-like parameter, $k_2$ is a dimensionless parameter, $\kappa$ is the fiber dispersion parameter, $I_{4,i} = \mathbb{C} : \left( \a_i \otimes \a_i \right)$ is the fourth invariant that represents the square of each fiber stretch, $\mathbb{C} = \FF^{\tran} \FF$ is the unmodified right Cauchy-Green tensor, and $\a_i$ is the fiber directional unit vector in the initial configuration.
Here we only consider the anisotropic energy if the fiber is active, i.e., $I_{4,i} \ge 1$.
In addition, the analytical study of the volumetric-isochoric split by Sansour~\cite{SANSOUR200828} demonstrates that the fiber-related energy terms should be formulated using the full deformation tensor because incorporating the split in $I_4$ can cause the fiber stresses to be transverse to its direction, which contradicts the expected physical behavior of a fiber.
Thus, we use the unmodified invariant in the fiber direction instead of using modified invariants for both the ground matrix and the fibers.
Our final strain energy functional, incorporating both isotropic and anisotropic contributions for realistic fiber-reinforced materials, is defined as follows:
\begin{align}\label{HGO_strain_energy}
\Psi = \Psi_{\text{isochoric}} + \sum_{i=1,2} \Psi_{\text{anisotropic},i} + \Psi_{\text{volumetric}}.
\end{align}

In the PD constitutive correspondence formulation, material anisotropy does not alter the definition of the non-local deformation gradient tensor $\FF$. 
Instead, the reference fiber directions enter the strain-energy functional through the anisotropic invariants, and the resulting anisotropic stress modifies the PD force state and further affects the deformation field and bond-stretch distribution. 
Although the same critical stretch values are assigned to all bonds as described in Sec.~\ref{s:failure model}, fiber reinforcement generally causes bonds oriented along the reinforced direction to reach the damage threshold later and produces an implicit directional resistance to damage and failure. 
However, the present failure criterion does not explicitly prescribe critical stretches or fracture toughness as functions of the bond-fiber orientation.

\section{Discrete formulations}
This section presents the discrete IPD formulations for FSI with and without structural damage and failure.
The numerical solution to the continuous incompressible Navier-Stokes equations, Eqs.~\eqref{e:navier_stokes}--\eqref{e:incompressibility}, is approximated by discretizing the Eulerian equations on a Cartesian grid and by discretizing the Lagrangian equations as a discrete point cloud in the reference coordinate.
Throughout the section, we describe discrete formulations in three spatial dimensions.
The discrete approximations employed in two spatial dimensions are similar.

\subsection{Eulerian discretization}
The incompressible Navier-Stokes equations, Eqs.~\eqref{e:navier_stokes} and \eqref{e:incompressibility}, are discretized in space using a second-order finite difference scheme on a staggered Cartesian grid~\cite{harlow1965numerical}.
The computational domain $\Omega = \left[0, L\right]^3$ is discretized by an $N \times N \times N$ Cartesian grid with a uniform grid spacing $h = L/N$ in the $x$-, $y$-, and $z$-directions. 
Integer indices $\left(i, j, k\right)$ label the Cartesian grid cells.
The centers of the Cartesian grid cells are located at $\x_{i,j,k} = \left(\left(i+\frac{1}{2}\right)h, \left(j+\frac{1}{2}\right)h, \left(k+\frac{1}{2}\right)h\right)$, where the pressure $p$ is approximated.
The discrete Eulerian velocity $\u = \left(u_1, u_2, u_3\right)$ is defined by vector components that are normal to the faces of the Cartesian grid cells. 
Specifically, $u_1$ is defined at $\x_{i-\frac{1}{2},j,k} = \left(ih, \left(j+\frac{1}{2}\right)h, \left(k+\frac{1}{2}\right)h\right)$, $u_2$ is defined at $\x_{i,j-\frac{1}{2},k} = \left(\left(i+\frac{1}{2}\right)h, jh, \left(k+\frac{1}{2}\right)h\right)$, and $u_3$ is defined at $\x_{i,j,k-\frac{1}{2}} = \left( \left(i+\frac{1}{2}\right)h, \left(j+\frac{1}{2}\right)h, kh\right)$.
The components of the discretized elastic body force density $\f = \left(f_1, f_2, f_3\right)$ are defined at the same locations as the corresponding velocity components.

$\grad_h$, $\grad_h \cdot$, and $\grad_h^2$ are the discrete gradient, divergence, and Laplace operators, respectively.
The discrete gradient of $p$ is approximated at the faces of the grid cells, and the discrete divergence of $\u$ is approximated at the cell centers. 
The discrete Laplacian of $\u$ is approximated component-wise at the faces of the grid cells.
The nonlinear advection term $\u \cdot \grad \u$ is computed using a version of the piecewise parabolic method~\cite{colella1984piecewise}.

\subsection{Lagrangian discretization}
\label{s:Lagrangian_discretization}
\subsubsection{Immersed structures}
\label{s:PD_discretization}
The immersed body $\Omega_0^{\text{s}}$ is discretized as a collection of $M$ PD nodes, i.e., $\X_m \in \Omega_0^{\text{s}}$ and $m = 1,2,\cdots,M$.
Let $\dhorizon_{\Xl} \subset \Omega_0^{\text{s}}$ be the set of interacting neighbors within a distance $\horizonsize$ of the Lagrangian marker $\Xl$.
Then each PD node interacts with a finite number of neighbors in its horizon.
It is natural that the spatial integrals in the continuous NOSB-PD formulation are discretized as volume-weighted sums. Thus, the non-local deformation gradient tensor $\FF$, Eq.~\eqref{nonlocal_deformation_gradient}, and shape tensor $\BB$, Eq.~\eqref{shape_tensor}, at the material point $\Xl$ are
\begin{align}
\FF_m&= \sum_{\Xn \in \dhorizon_{\Xl}}  \ \omegah \left(| \Xn - \Xl |, t\right) \, \Yubar \langle \Xn - \Xl \rangle \otimes \left( \Xn - \Xl \right) \BB_m^{-1} \, \Vn , \label{discrete_nonlocal_deformation_gradient}\\
\BB_m &= \sum_{\Xn \in \dhorizon_{\Xl}} \ \omegah \left(| \Xn - \Xl |, t\right) \left( \Xn - \Xl \right) \otimes \left( \Xn - \Xl\right) \, \Vn , \label{discrete_shape_tensor}
\end{align}
in which $\Xn$ is a neighbor of $\Xl$ in the peridynamic horizon $\dhorizon_{\Xl}$ and $\Vn$ is the volume occupied by material point $\Xn$.
Likewise, the discretized force vector state of material point $\Xl$ is 
\begin{align}\label{discrete_force_vector}
\Tubar\left[\Xl,t\right]\langle\Xn - \Xl\rangle = \omegah \left(| \Xn - \Xl |, t \right) \PP_m \BB_m^{-1} \left( \Xn - \Xl \right),
\end{align}
and the pairwise bond force function is 
\begin{align}\label{discrete_pairwise_bond_f}
\Fubar \left( \Xn, \Xl, t \right)  = \omegah \left(| \Xn - \Xl |, t \right) \left( \PP_m \BB_m^{-1} +  \PP_n \BB_n^{-1} \right) \left( \Xn - \Xl \right),
\end{align}
in which $\PP_n$ and $\BB_n$ are the discretized first Piola-Kirchhoff stress tensor and shape tensor of particle $\Xn$, respectively. 
The discretized first Piola-Kirchhoff stress tensor is computed by the classical constitutive relations as in Sec.~\ref{s:constitutive_laws}, but using the discretized non-local deformation gradient tensor instead of the classical deformation gradient tensor in the continuum theory. 
Consequently, the net internal body force density at material point $\Xl$ is 
\begin{align}\label{discrete_pd_force}
\F \left( \Xl,t \right) =  \sum_{\Xn \in \dhorizon_{\Xl}} \Fubar \left( \Xn,\Xl,t \right) \, \Vn.
\end{align}
This peridynamic net bond force density is used as an elasticity model for the immersed structure at material point $\Xl$ at time $t$. 

Similarly, the local damage at point $\Xl$ is expressed in a discretized form as follows:
\begin{align}\label{discrete_damage}
\varphi \left( \Xl,t \right) = 1 - \frac{\sum_{\Xn \in \dhorizon_{\Xl}} I[\Xl,t]\langle\Xn-\Xl\rangle \, \Vn}{\sum_{\Xn \in \dhorizon_{\Xl}}  \Vn}.
\end{align}
Note that the local damage is equal to 0 if its initial bonds are all active and the value is equal to 1 if all bonds are broken.

We use a cubic B-spline influence function for better smoothness and numerical stability in our numerical tests:
\begin{align}
\omega(r) = 
\begin{cases}
C \left( \frac{2}{3} - r^2 + \frac{r^3}{2} \right) \ &\text{if} \ 0 \le \ r < 1, \\
C\frac{(2-r)^3}{6} \ &\text{if} \ 1 \le \ r < 2, \\
0  \ &\text{otherwise},
\end{cases}
\end{align}
in which $r = \frac{2|\Xi|}{\horizonsize}$ and $C = \frac{15}{7\pi}$ in two spatial dimensions or $C = \frac{3}{2\pi}$ in three spatial dimensions. 
Unlike piecewise constant or linear kernels, the cubic B-spline provides $C^2$-continuity and a smooth decay to zero at the horizon boundary, which reduces spurious stress oscillations and artificial wave reflections~\cite{Seleson_2011,Seleson2018}. 
It also improves the conditioning of the shape tensor and enhances the accuracy of deformation gradient recovery, enabling higher-order consistency and a more stable transition to the local continuum limit~\cite{Hillman_2019}.
Note that the modified influence function, Eq.~\eqref{modified_inf_func}, is used in the IPD formulation.

\subsubsection{Lagrangian nodal volumes}
\begin{figure}[t!]
\centering
    \begin{tabular}{cc}
        \begin{subfigure}{.33\textwidth}
          		\includegraphics[width=\textwidth]{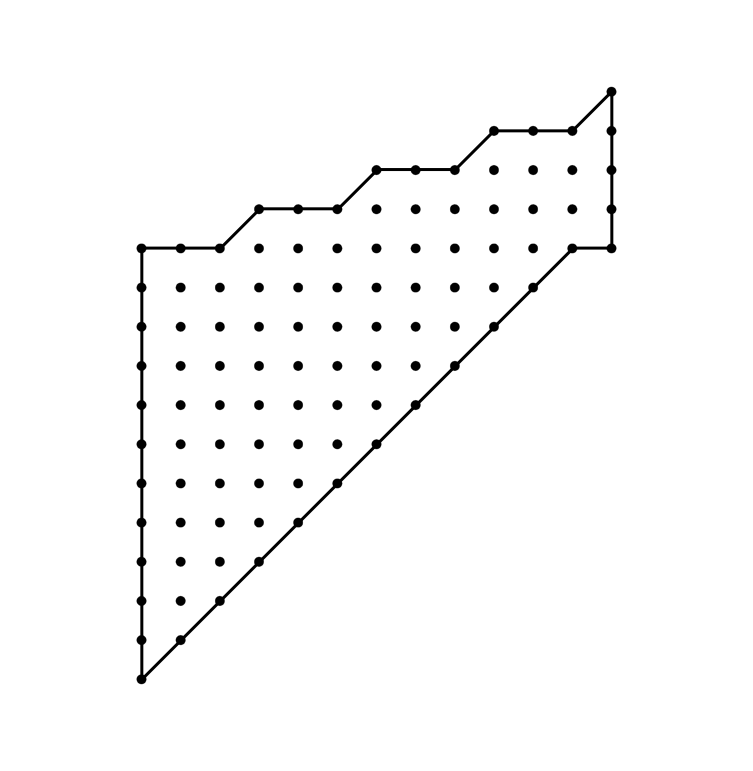}
          		\caption{Uniform lattice}
          		\label{f:uniform_discretization}
        \end{subfigure} 
        \begin{subfigure}{.33\textwidth}
                \includegraphics[width=\textwidth]{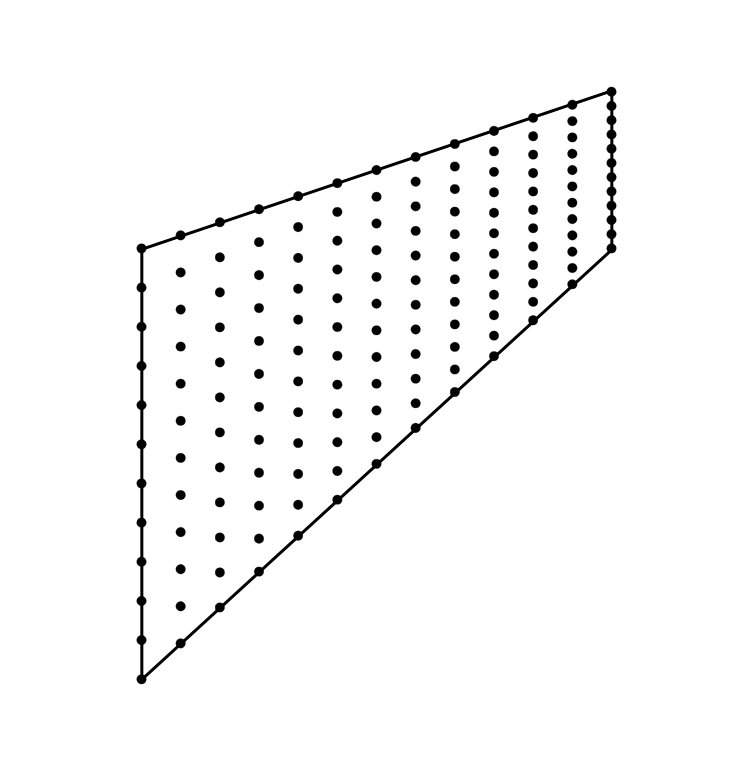}
                \caption{Rescaled lattice}
        \end{subfigure}
        \begin{subfigure}{.33\textwidth}
                \includegraphics[width=\textwidth]{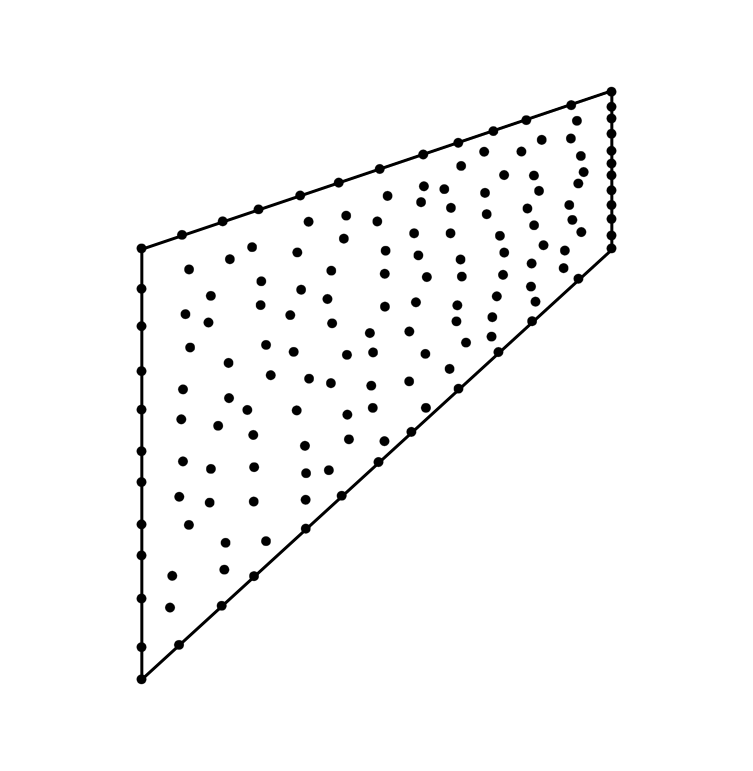}
                \caption{Irregular distribution}
        \end{subfigure}
    \end{tabular}
    \caption{Three different material discretizations of Cook's membrane. The uniform discretization imposes uniformly distributed material volumes along the material body. On the other hand, the material volumes in the rescaled and irregular representations require exact volume computations to avoid nonphysical behaviors. A solid line indicates the boundary represented by each discretization and black circles are PD structural nodes.}
    \label{f:discretizations}
\end{figure}

For simplicity of the geometric representation, the initial version of the IPD method adopts a regular lattice of PD points perfectly aligned with the background Cartesian grid and uniformly distributed volumes along the material body, which is sufficient for simple material geometries such as a rectangle.
However, such a uniform discretization causes the use of a stair-step geometry that represents unrealistic material boundaries, as shown in Fig.~\ref{f:uniform_discretization}.
In addition, the use of regularly spaced structural points can cause a grid bias that complicates the simulation of complex crack patterns.
Specifically, the structured lattice creates preferred directions for bond breakage, often causing fractures to propagate artificially along the grid lines or diagonals rather than following physically driven, arbitrary paths.

To overcome these limitations, herein we employ non-uniform structural discretizations.
We focus on two different non-uniform distributions of PD nodes along the immersed structure to describe realistic material geometries; see Fig.~\ref{f:discretizations}.
PD nodes in the rescaled lattice are structured but have a non-uniformly distributed volume at each node, and the irregular distribution has irregularly displaced nodes and non-uniform volumes along the immersed body.
However, such an unstructured material representation causes a non-uniform volume distribution along the material body and requires an accurate computation of the nodal volume occupied by each material point.

\begin{figure}[t!]
\centering
\includegraphics[width=.5\textwidth]{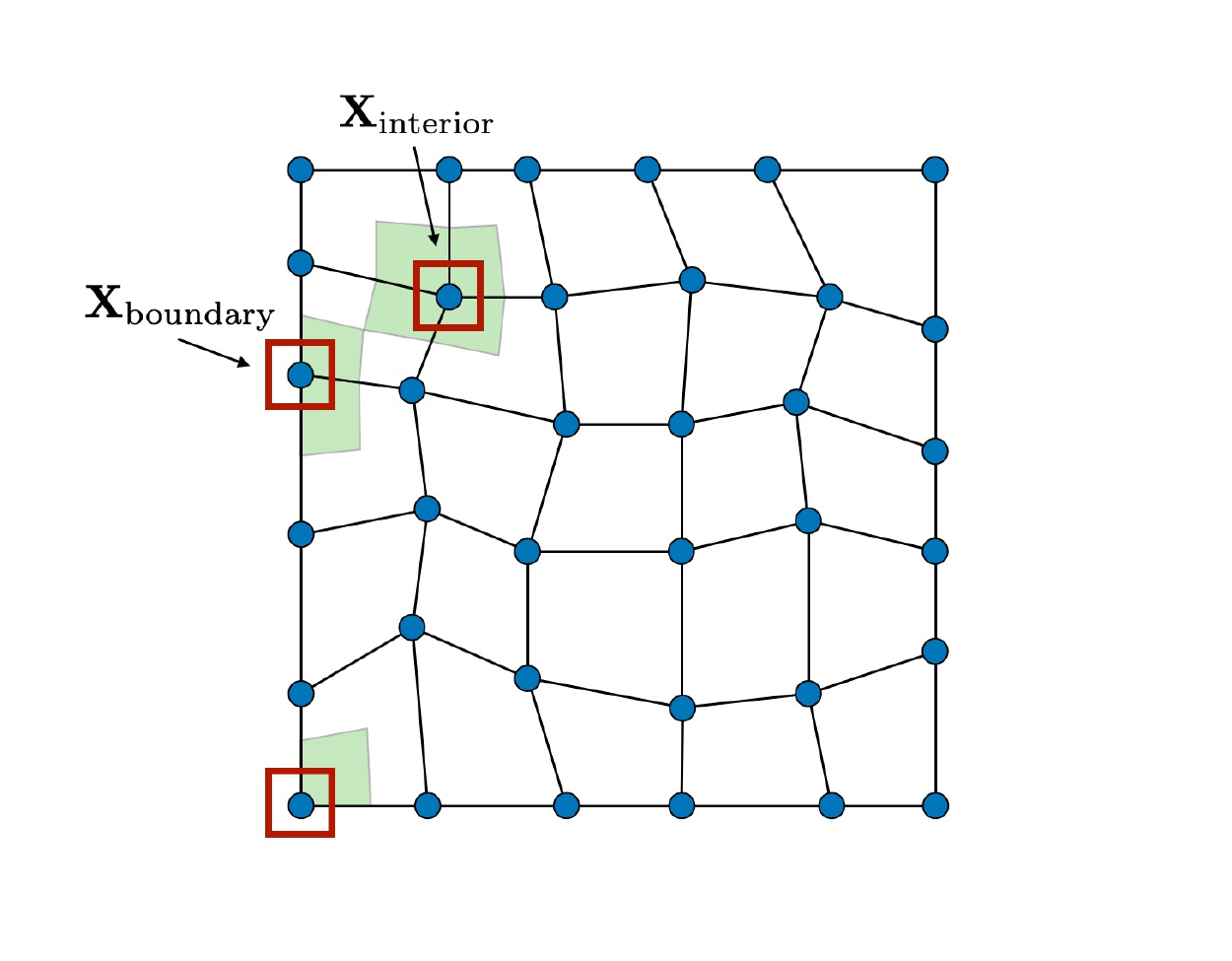} 
\caption{An example of non-uniform discretization with FE quadrilateral meshes and corresponding nodes. The blue points are FE and PD nodes that share the same locations and the green shaded regions are areas (i.e., volumes) computed by FE meshes. The nodal volumes achieved by the FE elements on the boundary need to be corrected for accurate PD volumetric forces to avoid insufficient PD volumetric forces in the IPD formulation.}
\label{f:FE-PD-nodes}
\end{figure}

For the consistent evaluation of the volumes occupied by the Lagrangian points in the non-uniform discretizations, we utilize FE meshes~\cite{hughes2012finite}.
A primary advantage of this approach is that the FE basis functions, which form a partition of unity, can be used to derive consistent quadrature weights.
We assume that the structural body $\Omega_0^{\text{s}}$ is discretized into a collection of FE elements $\{\Omega_{\text{e}}^{\text{s}}\}$.
The volume $V_n$ assigned to a material point $\Xn$ is computed by integrating the basis functions over the elements adjacent to that point:
\begin{align}\label{e:nodal_volume}
V_n = \sum_{e \in \text{adj}\left(\Xn\right)} \int_{\Omega_{\text{e}}^{\text{s}}} \vec{\vec{\phi}}_n^{\text{e}} \left(\X\right) \mathrm{d}\X,
\end{align}
in which $\text{adj}\left(\Xn\right)$ is a set of FE elements sharing $\Xn$ and $\vec{\vec{\phi}}_n^{\text{e}}$ is a basis function associated with the material point $\Xn$ supported on the element $e$.
For a non-affine element mapping, the determinant of the geometric Jacobian varies within the element and is evaluated at each quadrature point using standard FE quadrature.
Consequently, curved boundaries and other geometric variations are incorporated according to their discrete representation in the undeformed reference FE mesh.
Note that the nodal volumes are computed only once during preprocessing and remain fixed throughout the simulation.

The nodal volume computed by Eq.~\eqref{e:nodal_volume} in the two-dimensional formulation is an area occupied by $\Xn$, which must be multiplied by the out-of-plane thickness of the structural body to obtain the volume used in the IPD force calculation.
In addition, Lagrangian points on the boundary require a volume treatment that can recover their full volumes on the boundary layer after computing their partial volumes from FE meshes; see Fig.~\ref{f:FE-PD-nodes}.
Otherwise, because of the lack of PD volumetric forces near the boundary, the structure undergoes unstable deformations.
For the locally regular quadrilateral meshes, nodal weights on boundary edges are multiplied by two, whereas those at corners are multiplied by four.
Similarly, for the hexahedral meshes, the nodal weights associated with boundary faces, edges, and corners are multiplied by two, four, and eight, respectively.
For a sufficiently resolved smooth curved boundary, these correction factors provide a locally planar approximation of the missing nodal support.
More general geometries with boundaries of high curvature or non-orthogonal edges and corners may require geometry-dependent correction factors.

Such an FE-based geometric representation offers the potential for matching material geometries in a coupling scheme between the IPD method and FE-based immersed boundary methods.
In general, a meshfree method requires high computational costs compared to FE-based methods~\cite{MACEK20071169,Kilic_2010}.
Consequently, we expect that a hybrid coupling strategy based on this shared structural discretization could improve computational efficiency by using the IPD formulation in regions where crack initiation and propagation are expected and an FE-based formulation in the remainder of the structural domain.

In the current IPD formulation, the FE-based nodal-volume computation introduces a one-time preprocessing cost that is absent from the uniform discretization.
For $M$ PD material points, the uniform discretization stores a single volume value, whereas the non-uniform discretization computes and stores $M$ distinct nodal-volume values based on the reference configuration.
Since these values remain fixed throughout the simulation, no additional FE integration or volume update is required during the simulation.
For a fixed set of interacting bonds, both formulations perform the same volume-weighted PD summations as introduced in Sec.~\ref{s:PD_discretization}, with the non-uniform formulation requiring only an additional array access to retrieve the corresponding nodal volume.
Therefore, the computational overhead attributable specifically to the FE-based nodal-volume representation is negligible.
However, a non-uniform point distribution may alter, and locally increase, the number of interacting neighbors within each PD horizon.
Consequently, the total cost of the PD interaction calculation may increase from that of a uniform discretization because of differences in bond connectivity rather than the FE-based computation of the nodal volumes.

In addition, the volume associated with a material point $\Xn$ near the boundary of the horizon $\dhorizon_{\Xl}$ is often only partially located inside the $\epsilon$-ball, and summing the full contributions of these boundary nodes leads to an overestimation of the net PD body force density.
To improve the accuracy of numerical solutions, we use a volume correction method~\cite{hu2010numerical, Seleson2018}:
\begin{align}\label{e:horizon_vol_correction}
\Vn^{(m)} = \begin{cases}
\Vn, \ &\text{if} \ |\Xi| \le \horizonsize - \frac{\Delta X}{2}, \\
\frac{1}{\Delta X} \left[ \horizonsize - \left(|\Xi| - \frac{\Delta X}{2}\right) \right] \Vn, \ &\text{if} \ \horizonsize - \frac{\Delta X}{2} < |\Xi| \le \horizonsize, \\
0, \ &\text{otherwise},
\end{cases}
\end{align}
in which $\Xi$ is a bond connecting two PD nodes $\Xl$ and $\Xn$ in the reference configuration.
Instead of using a full nodal volume $\Vn$ near the boundary of the PD horizon, this partial volume correction $\Vn^{(m)}$ is used in Eqs.~\eqref{discrete_nonlocal_deformation_gradient}--\eqref{discrete_damage}.

\subsubsection{Volumetric stabilization}
Conventional IB-type methods frequently struggle with poor volumetric conservation due to FSI coupling operators, regularized delta functions, the discrete divergence, and time-stepping errors~\cite{vadala2020stabilization, lee2022lagrangian, GRUNINGER2026114472}.
To prevent non-physical volume leakage, we utilize a stabilization technique previously developed for the IFED method by Vadala-Roth et al.~\cite{vadala2020stabilization} by incorporating a volumetric stabilization term into the strain energy functional.
This introduces an additional restoring force in the structural domain and resists spurious compressible motions under large deformations.
The penalty bulk modulus $\kappa_{\mathrm{stab}}$ is parameterized using a numerical Poisson's ratio $\nu_{\mathrm{stab}}$:
\begin{align}\label{numerical_bulk_modulus}
\kappa_{\mathrm{stab}} = \frac{2G\left(1+\nu_{\mathrm{stab}}\right)}{3\left(1 - 2 \nu_{\mathrm{stab}}\right)},
\end{align}
in which $G$ is a shear modulus of the immersed structure.
This volumetric parameter controls the discrete incompressibility of the immersed structure.
Note that the limit as $\nu_{\mathrm{stab}}$ approaches $0.5$ yields exact incompressibility.

\subsection{Lagrangian--Eulerian coupling}
In the continuous equations, coupling between Eulerian and Lagrangian variables is achieved by integral transforms with Dirac delta function kernels as in Eqs.~\eqref{e:force_spreading}--\eqref{e:velocity}.
In the discrete formulation, the singular delta function is replaced by a regularized delta function $\delta_h$, which is formed as a tensor product of one-dimensional kernel functions and defined by
\begin{align}\label{regularized_delta_func}
\delta_h \left(\x\right) = \prod_{k=1}^3 \delta_h \left(x_k\right) = \frac{1}{h^3} \phi \left( \frac{x_1}{h}\right) \phi \left( \frac{x_2}{h}\right) \phi \left( \frac{x_3}{h}\right),
\end{align}
in which $\phi\left(r\right)$ is a basic one-dimensional immersed boundary kernel function~\cite{peskin_2002}. 

Similar to PD volume integrals in Sec.~\ref{s:PD_discretization}, the volume integral Eq.~\eqref{e:force_spreading} can be approximated by 
\begin{align}\label{IB_force_spreading}
\left( f_1 \right)_{i - \frac{1}{2},j,k} &= \sum\limits_{m}  F_{m,1} \, \delta_h \left( \x_{i - \frac{1}{2},j,k} - \vchi(\Xl,t) \right) V_m,\\
\left( f_2 \right )_{i,j - \frac{1}{2},k} &= \sum\limits_{m} F_{m,2} \, \delta_h \left( \x_{i,j - \frac{1}{2},k} -  \vchi(\Xl,t) \right) V_m, \\
\left( f_3 \right )_{i,j,k - \frac{1}{2}} &= \sum\limits_{m} F_{m,3} \, \delta_h \left( \x_{i,j,k - \frac{1}{2}} -  \vchi(\Xl,t) \right) V_m,
\end{align}
in which $\F_m = \left( F_{m,1}, F_{m,2}, F_{m,3} \right)$ is the Lagrangian force density at a Lagrangian marker of index $m$ and $V_m$ is a nodal volume occupied by $\Xl$.
Each nodal volume is accurately computed using the boundary-corrected FE weight.
Note that this structural quadrature volume is independent of the Eulerian grid spacing $h$ and is distinct from the bond-dependent partial-volume correction used to approximate PD integrals near the horizon boundary as detailed in Eq.~\eqref{e:horizon_vol_correction}.
We use the notation
\begin{align}\label{force_spreading_operator}
\f = \S\left[ \vchi \left(\cdot,t\right)\right] \F,
\end{align}
in which $\S\left[\vchi \left(\cdot,t\right)\right]$ is the discrete force-spreading operator. The structural body interacts with the surrounding fluid by spreading the force to the Eulerian grid and moves with the local fluid velocity. Similarly, Lagrangian and Eulerian velocities can be related by
\begin{align}
U_{m,1}  \left(\X,t\right) &= \sum\limits_{i,j,k} (u_1)_{i-\frac{1}{2},j,k} \, \delta_h \left( \x_{i - \frac{1}{2},j,k} - \vchi \left( \X,t \right) \right) h^3, \label{velocity_interpolation_1}\\
U_{m,2}  \left(\X,t\right) &= \sum\limits_{i,j,k} (u_2)_{i,j-\frac{1}{2},k} \, \delta_h \left( \x_{i,j-\frac{1}{2},k} - \vchi \left( \X,t \right) \right) h^3, \label{velocity_interpolation_2} \\
U_{m,3}  \left(\X,t\right) &= \sum\limits_{i,j,k} (u_3)_{i,j,k-\frac{1}{2}} \, \delta_h \left( \x_{i,j,k-\frac{1}{2}} - \vchi \left( \X,t \right) \right) h^3, \label{velocity_interpolation_3}
\end{align}
in which $\U_m = \left( U_{m,1}, U_{m,2}, U_{m,3} \right)$ is the Lagrangian velocity at a Lagrangian marker $\Xl$. We use the notation
\begin{align}
\U = \J \left[ \vchi \left(\cdot,t\right) \right] \u,
\end{align}
in which $\J$ is the discrete velocity restriction or interpolation operator. 

To establish weighted adjointness in three spatial dimensions, we define the Eulerian inner product on the staggered Cartesian grid by
\begin{align}
\left\langle \u,\f \right\rangle_{\mathrm{E}}
=\sum_{i,j,k}
\left(
(u_1)_{i-\frac{1}{2},j,k}
(f_1)_{i-\frac{1}{2},j,k}
+
(u_2)_{i,j-\frac{1}{2},k}
(f_2)_{i,j-\frac{1}{2},k}
+
(u_3)_{i,j,k-\frac{1}{2}}
(f_3)_{i,j,k-\frac{1}{2}}
\right) \, h^3 ,
\label{e:eulerian_inner_product}
\end{align}
and the corresponding volume-weighted Lagrangian inner product by
\begin{align}
\left\langle \U,\F \right\rangle_{\mathrm{L}}
=
\sum_m
\left(
U_{m,1}F_{m,1}
+
U_{m,2}F_{m,2}
+
U_{m,3}F_{m,3}
\right)V_m.
\label{e:lagrangian_inner_product}
\end{align}
Note that the Eulerian inner product is weighted by the Cartesian cell volume $h^3$, whereas the Lagrangian inner product is weighted by the boundary-corrected structural nodal volume $V_m$.
Substituting the discrete force-spreading equations into Eq.~\eqref{e:eulerian_inner_product} yields
\begin{align}
\left\langle \u,\S\F \right\rangle_{\mathrm{E}}
=\sum_{i,j,k}
\biggl(
&(u_1)_{i-\frac{1}{2},j,k}
\sum_m F_{m,1}
\delta_h
\left(
\x_{i-\frac{1}{2},j,k}
-
\vchi(\Xl,t)
\right)V_m
\nonumber\\
\quad+
&(u_2)_{i,j-\frac{1}{2},k}
\sum_m F_{m,2}
\delta_h
\left(
\x_{i,j-\frac{1}{2},k}
-
\vchi(\Xl,t)
\right)V_m
\nonumber\\
\quad+
&(u_3)_{i,j,k-\frac{1}{2}}
\sum_m F_{m,3}
\delta_h
\left(
\x_{i,j,k-\frac{1}{2}}
-
\vchi(\Xl,t)
\right)V_m
\biggr) \, h^3.
\end{align}
Interchanging the Eulerian and Lagrangian summations gives
\begin{align}
\left\langle \u,\S\F \right\rangle_{\mathrm{E}}
={}
\sum_m V_m
\biggl(
&F_{m,1}
\sum_{i,j,k}
(u_1)_{i-\frac{1}{2},j,k}
\delta_h
\left(
\x_{i-\frac{1}{2},j,k}
-
\vchi(\Xl,t)
\right)h^3
\nonumber\\
\quad+
&F_{m,2}
\sum_{i,j,k}
(u_2)_{i,j-\frac{1}{2},k}
\delta_h
\left(
\x_{i,j-\frac{1}{2},k}
-
\vchi(\Xl,t)
\right)h^3
\nonumber\\
\quad+
&F_{m,3}
\sum_{i,j,k}
(u_3)_{i,j,k-\frac{1}{2}}
\delta_h
\left(
\x_{i,j,k-\frac{1}{2}}
-
\vchi(\Xl,t)
\right)h^3
\biggr).
\end{align}
The three Eulerian summations are precisely the interpolated velocity components $U_{m,1}$, $U_{m,2}$, and $U_{m,3}$, respectively.
Consequently,
\begin{align}
\left\langle \u,\S\F \right\rangle_{\mathrm{E}}
=
\sum_m
\left(
F_{m,1}U_{m,1}
+
F_{m,2}U_{m,2}
+
F_{m,3}U_{m,3}
\right)V_m
=
\left\langle \U,\F \right\rangle_{\mathrm{L}}
=
\left\langle \J\u,\F \right\rangle_{\mathrm{L}},
\label{e:weighted_adjointness}
\end{align}
which demonstrates that $\S$ and $\J$ satisfy $\J=\S^\ast$ with respect to the Eulerian and volume-weighted Lagrangian inner products defined above.
This remains valid for non-uniform nodal volumes since the same structural quadrature weights $V_m$ are used in the force-spreading operator and the Lagrangian inner product~\cite{griffith2017hybrid}.
Note that both coupling operators must be evaluated using the same structural configuration and regularized delta function.
In addition, the weighted adjoint relation directly implies the discrete power balance for the Lagrangian--Eulerian coupling in the discrete IPD formulation.

\subsection{Time-stepping algorithm}
We now outline the numerical implementation of the IPD method used for the simulations presented in Sec.~\ref{s:ansio_benchmark}. 
Let $\u^n$ and $\vchi^n$ denote the discrete Eulerian velocity field and the structural deformation at time $t^n = n \Delta t$, respectively. 
The coupled system is advanced in time using a formally second-order accurate scheme that employs a midpoint rule for temporal discretization~\cite{griffith2017hybrid}.
We first approximate the deformed structural configuration at time $t^{n+\frac{1}{2}}$ using the current velocity field:
\begin{align}
\frac{\vchi^{n+\frac{1}{2}} -\vchi^n}{\Delta t / 2} &= \J \left[ \vchi^n \right] \u^n. \label{intermediate_lag}
\end{align}
Next, the Eulerian momentum equation is discretized using a Crank-Nicolson scheme for the viscous term and an explicit Adams-Bashforth approximation for the nonlinear advection term as follows:
\begin{align}
\rho \left( \frac{\u^{n+1} -\u^n}{\Delta t} + \N^{\left(n + \frac{1}{2}\right)}\right) &= - \grad_h p^{n + \frac{1}{2}} + \mu \grad_h^2  \left( \frac{\u^{n+1} + \u^n}{2}\right) + \f^{n + \frac{1}{2}}, \label{discretized_navier_stokes}\\
\grad_h \cdot \u^{n+1} &= 0, \label{discretized_incompressibility}
\end{align}
in which $\N^{\left(n + \frac{1}{2}\right)} = \frac{3}{2}\u^{n} \cdot \grad_h \u^{n} - \frac{1}{2} \u^{n -1} \cdot \grad_h \u^{n -1}$ is an explicit approximation to the nonlinear advection term in Eq.~\eqref{e:navier_stokes}
 and  $\f^{n + \frac{1}{2}} = \S \left[ \vchi^{n + \frac{1}{2}} \right] \F^{n + \frac{1}{2}}$ is the spreading of the Lagrangian force density. %
To achieve global second-order accuracy in time, a two-step predictor-corrector scheme is used to initialize the simulation for the first step~\cite{griffith2017hybrid}.
Eqs.~\eqref{discretized_navier_stokes}--\eqref{discretized_incompressibility} are approximated by a projection method, which decouples the velocity and pressure updates.
Finally, the structural configuration is updated to the next time step using the averaged velocity field:
\begin{align}
\frac{\vchi^{n+1} -\vchi^n}{\Delta t} &= \J \left[ \vchi^{n + \frac{1}{2}} \right] \left( \frac{\u^{n+1} + \u^n}{2} \right). \label{update_lag}
\end{align}

During the simulation, the non-local deformation gradient tensor, shape tensor, and PD pairwise bond forces are updated at every time step to describe large deformations and topological changes (i.e., bond breakage).
The net Lagrangian force density $\F^{n + \frac{1}{2}}$ is derived using the constitutive correspondence framework as illustrated in Sec.~\ref{s:immersed body forces}.

\section{Benchmarks}
\label{s:ansio_benchmark}
Prior to investigating the performance of the IPD method with anisotropic materials, we first evaluate the influence of different Lagrangian point distributions. 
In particular, we use the Cook's membrane benchmark~\cite{cook1974improved}.
We then consider anisotropic versions of several standard benchmark cases in solid mechanics, including the compressed block~\cite{reese1999new}, a three-dimensional version of the Cook's membrane test with anisotropy~\cite{cook1974improved}, anisotropic torsion test~\cite{bonet1997nonlinear}, and the adventitial strip model introduced by Gasser, Ogden, and Holzapfel~\cite{Gasser_2005}.
We compare the accuracy and convergence of the proposed IPD method against an FE-based IB method called the IFED method introduced by Griffith and Luo~\cite{griffith2017hybrid}.
To impose homogeneous material characteristics throughout the immersed structure, we use a uniform ball-shaped horizon like that used in the isotropic tests in our previous work.
To further demonstrate grid-converged damage growth and failure processes, we use a two-dimensional tissue strip with different choices of fiber orientations and show the effects of fibers on damage growth and failure. 
We also investigate a purely fluid-driven failure process by modifying the elastic band test~\cite{vadala2020stabilization,WELLS2023111890}.
All IPD simulations utilize the IBAMR software~\cite{griffith2007adaptive, ibamr}, which is a distributed-memory parallel implementation of the IB method with support for Cartesian grid adaptive mesh refinement (AMR). 

Unless otherwise noted, we set the density and viscosity of the fluid to $\rho = 1.0 \, \frac{\text{g}}{\text{cm}^3}$ and  $\mu = 0.01 \,  \frac{\text{dyn$\cdot$s}}{\text{cm}^2}$, respectively. 
The computational domain is $\Omega = [0 , L]^d$, in which $d$ is the spatial dimension and $L$ is the length of the domain. 
The relative grid spacing between the Cartesian grid and structural grid is defined by the mesh factor ratio $M_{\mathrm{FAC}} = \frac{\Delta X}{\Delta x}$, in which $\Delta X$  and $\Delta x$ are the Lagrangian and Eulerian grid spacings, respectively. 
The Eulerian grid size is $\Delta x = \frac{L}{N}$, in which $N$ is the number of Cartesian grid cells in one spatial direction. 
We use $M_{\mathrm{FAC}} = 0.5$ in our IPD simulations, unless otherwise noted, so that the structural discretization is twice as fine as the background Cartesian grid.
For the regularized delta kernel function in our IPD simulations, we use the four-point IB kernel function introduced by Peskin~\cite{peskin_2002} unless otherwise mentioned. 

We also investigate the effect of peridynamic horizon size in our benchmarks, and our grid refinement studies consider the asymptotic convergence of the IPD formulation.
For simplicity, a uniform $\horizonsize$-ball is used for the peridynamic horizon. 
The peridynamic horizon size $\horizonsize$ is always taken to be a constant multiple of the Lagrangian mesh spacing $\Delta X$ in the reference configuration, which is commonly used to define the $\horizonsize$-ball in the PD literature~\cite{madenci2014peridynamic,behera2020peridynamic,Wang_2020}.
A larger horizon implies more non-local effects in the structural body, which requires more computations compared to a smaller horizon size. 
Our simulations examine different peridynamic horizon sizes for the constitutive correspondence: $\horizonsize = 1.015  \Delta X, \ 2.015  \Delta X, \ 3.015 \Delta X$.
Such non-integer coefficients are chosen to have all interacting PD particles within the PD horizon.
We employ a sufficiently large horizon size $\horizonsize = \ 3.015 \Delta X$ for failure tests, as suggested in the PD literature~\cite{madenci2014peridynamic}.

For the failure tests, the critical stretch values are treated as phenomenological numerical parameters rather than experimentally calibrated material properties.
Based on preliminary tests, $\sca$ is selected separately for each benchmark so that bond damage occurs within the prescribed loading range.
Additionally, we set $\scb = 1.2 \, \sca$, which provides a numerically well-resolved softening transition and consistent damage and failure patterns under grid refinement.

As specified in Sec.~\ref{s:PD_discretization}, the cubic B-spline influence function is held fixed throughout all numerical tests. 
Because the influence function determines the relative weighting of bonds within the PD horizon, alternative choices may affect the discretization error and observed grid-refinement behavior as well as damage evolution. 
Therefore, the reported grid consistency and selected critical stretch values are specific to the chosen cubic B-spline influence function, and the critical stretch values and normalized softening width may require recalibration if a different influence function is adopted.

\subsection{Two-dimensional isotropic Cook's membrane}
\label{s:2d_cooks_discretizations}
\begin{figure}[t!]
\centering
    \includegraphics[width=.45\textwidth]{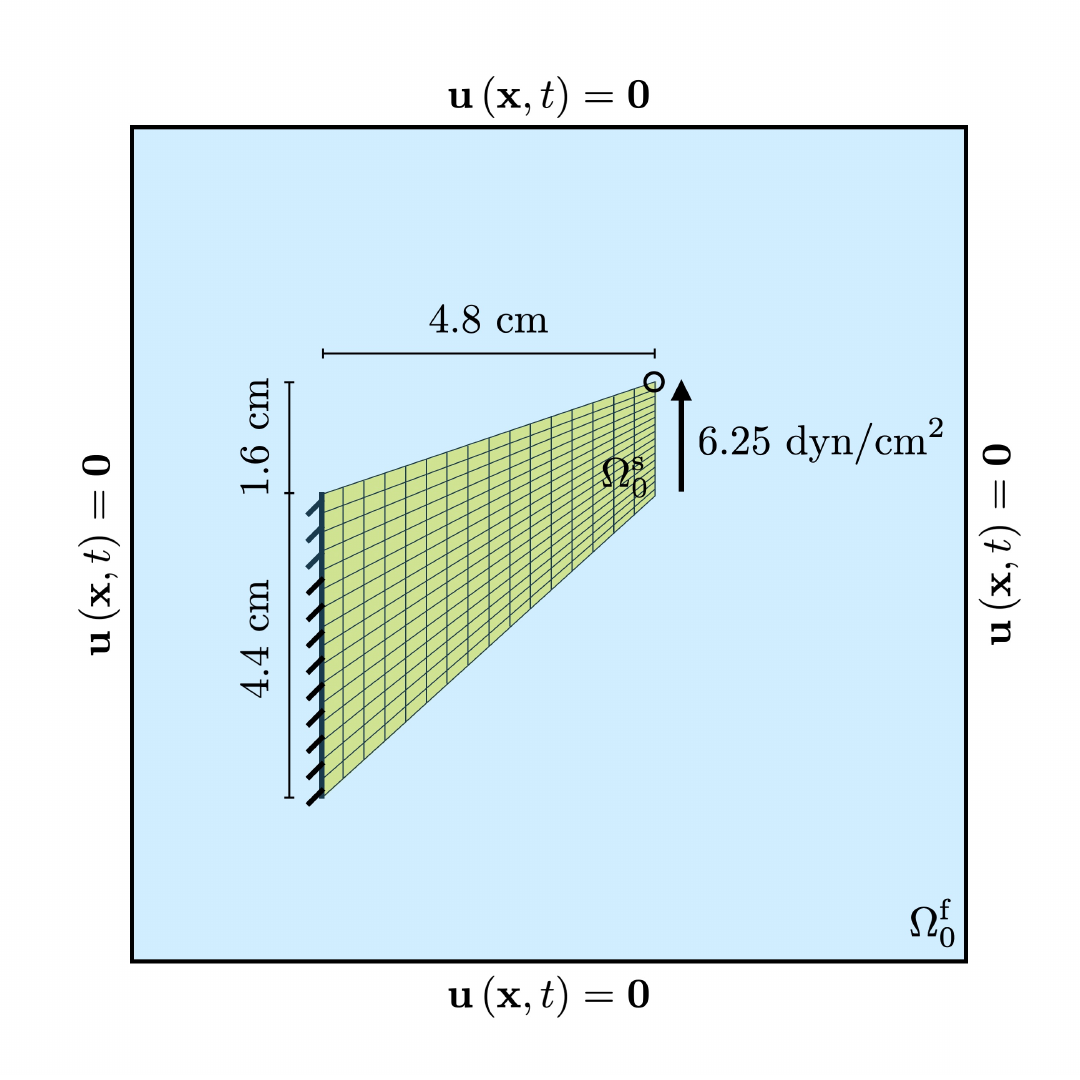}
    \caption{Schematic diagram for the Cook's membrane benchmark. The initial configurations of the immersed structure and a fluid are denoted by $\Omega_0^{\text{s}}$ and $\Omega_0^{\text{f}}$, respectively. The entire computational domain is $\Omega = \Omega_0^{\text{s}} \cup \Omega_0^{\text{f}}$. Zero fluid velocity is enforced on the outer boundaries of the computational domain.}
    \label{f:2d_Cooks_schematics}
\end{figure}

The computational domain is $\Omega = \left[0, L\right]^2$, with $L = 40 \,  \text{cm}$. 
Zero displacement conditions are imposed on the left boundary of the structure, and an upward traction of $6.25 \, \frac{\text{dyn}}{\text{cm}^2}$ is applied along the right boundary. 
Otherwise, stress-free boundary conditions are assumed.  
Fig.~\ref{f:2d_Cooks_schematics} provides a schematic of this test case. 
We use a modified neo-Hookean material model:
\begin{align}
\Psi &=\frac{G}{2}\left( \bar{I}_1 - 2 \right) + \frac{\kappa_{\text{stab}}}{2} \left( \ln J \right)^2.
\end{align}
The modified neo-Hookean material model~\cite{vadala2020stabilization} utilized in the stabilized IFED method uses $\bar{\FF} = J^{-1/3} \, \FF$, but we use $\bar{\FF} = J^{-1/2} \, \FF$ for two-dimensional tests.
A shear modulus of $G = 83.3333 \, \frac{\text{dyn}}{\text{cm}^2}$ is used for the incompressible hyperelastic material.
The isotropic benchmark tests in our previous work show that the numerical Poisson's ratio of $\nu_{\text{stab}} = 0.4$ yields comparable results to the conventional method, so we only consider $\nu_{\text{stab}} = 0.4$ to evaluate the numerical bulk modulus $\kappa_{\text{stab}}$ in this test.
The load time is $T_{\text{l}} = 20 \, \text{s}$, the final time is $T_{\text{f}} = 50 \, \text{s}$, and the damping parameter is set to $\eta = 4.16667 \,  \frac{\text{g}}{\text{cm}^3 \cdot \text{s}}$. 
Three different grid representations  (uniform, rescaled, and irregular lattices) are examined as illustrated in Fig.~\ref{f:discretizations}.
We focus on the vertical displacements of the top-right corner of the membrane to assess convergence by using similar numbers of solid degrees of freedom (DoF). 
Here, the solid DoF are defined as the number of Lagrangian material points.
In addition, the horizon size must be large enough ($\horizonsize \ge \sqrt{2} \Delta X$) to ensure adequate bond connectivity throughout the material. 
Without sufficient bond connectivity, particularly on the boundary, we observe structural instabilities because of the lack of PD volumetric forces.
Bond breakage is not considered in this problem.

 \begin{figure}[t!]
\centering
	\includegraphics[width=\textwidth]{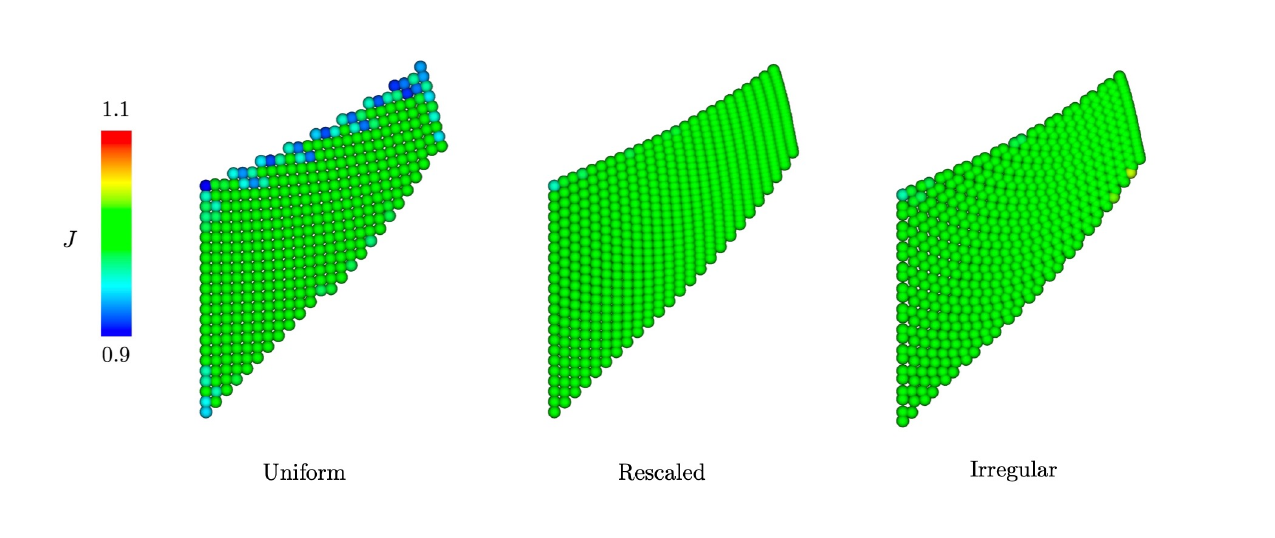} 
    \caption{Deformations of two-dimensional Cook's membrane with the values of $J$ at material points using the neo-Hookean material model with $G = 83.3333 \, \frac{\text{dyn}}{\text{cm}^2}$. The left panel shows uniformly distributed volumes, the center panel shows non-uniformly distributed volumes, and the right panel shows irregularly distributed volumes. The deformations are represented using $381$ solid DoF for the uniform discretization and $575$ solid degrees of freedom (DoF) for the rescaled and irregular discretizations. We set $\horizonsize = 2.015 \Delta X$.  The numerical Poisson's ratio is fixed at $\nu_{\mathrm{stab}}  = 0.4$.}
    \label{f:2d_Cooks_deformation}
\end{figure}

\begin{figure}[t!]
\centering
    \begin{tabular}{cc}
        \begin{subfigure}{.4\textwidth}
          		\includegraphics[width=\textwidth]{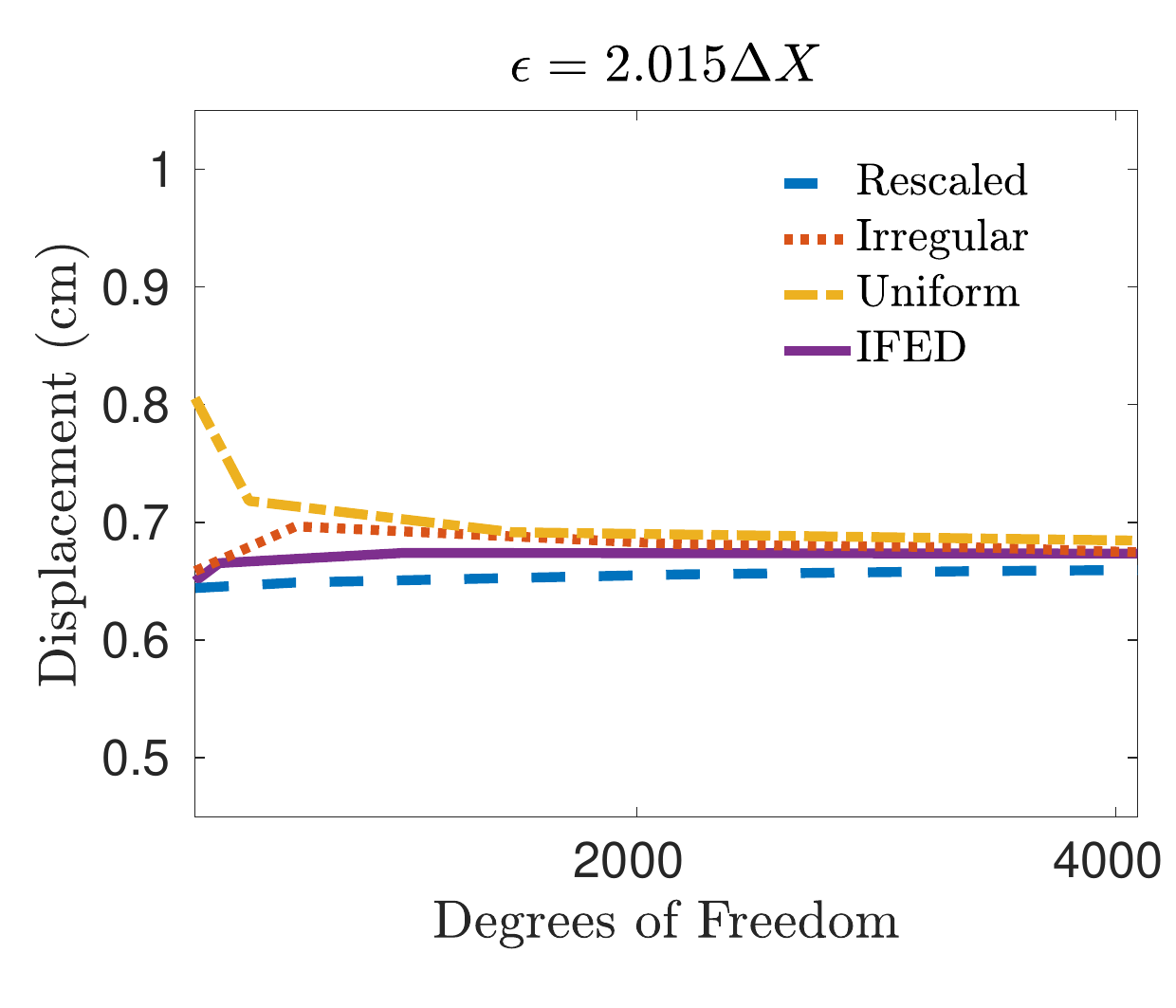}
        \end{subfigure}
        \hspace{.03\textwidth}
        \begin{subfigure}{.4\textwidth}
                \includegraphics[width=\textwidth]{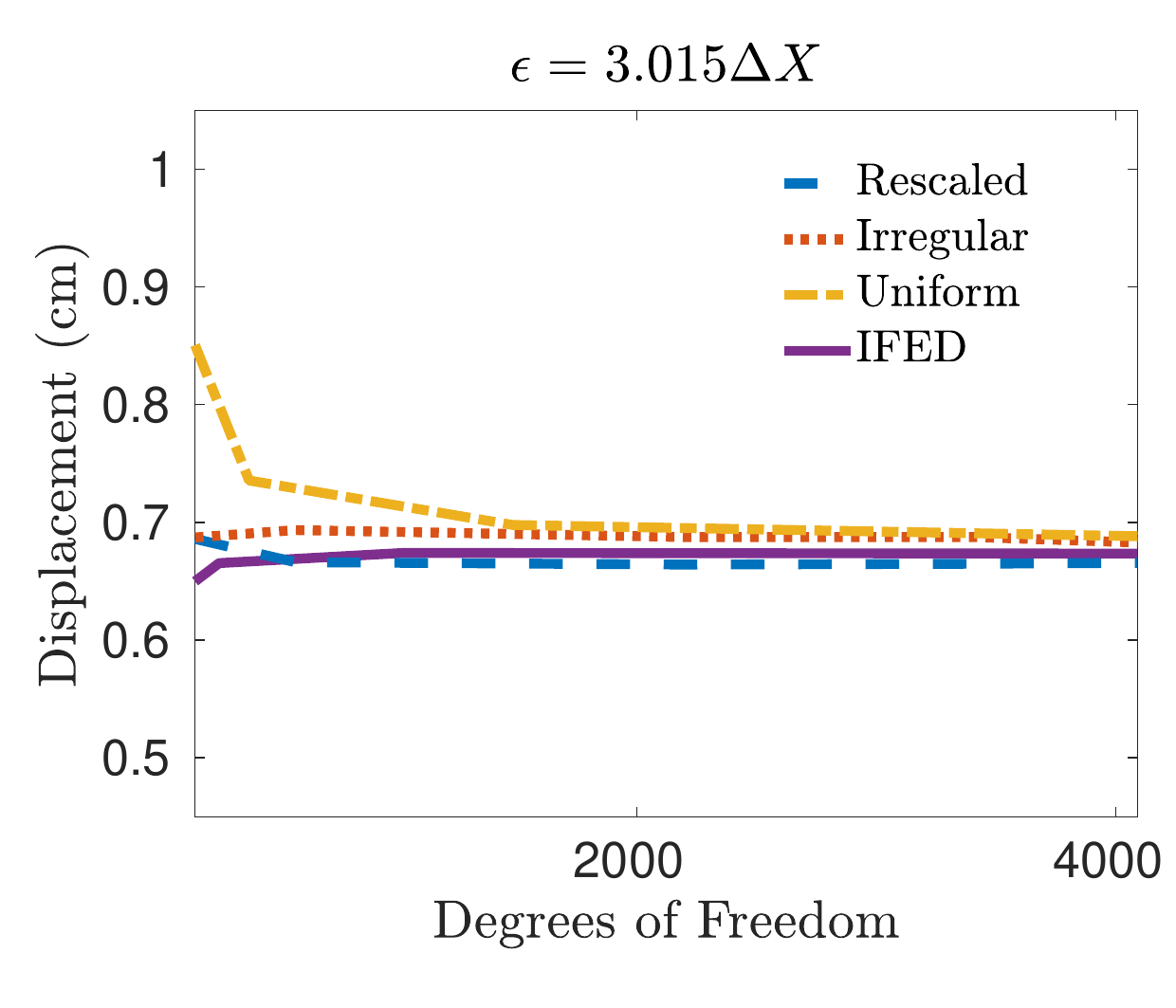}
        \end{subfigure}
    \end{tabular}
    \caption{Vertical displacements of the top corner point of the Cook's membrane benchmark, highlighted in Fig.~\ref{f:2d_Cooks_schematics}, for different choices of discretizations and peridynamic horizon size $\horizonsize$. The numerical Poisson's ratio is set to $\nu_{\mathrm{stab}}  = 0.4$. The solid DoF range from $156$ to $4096$. }
    \label{f:2d_Cooks_disp}
\end{figure}

Fig.~\ref{f:2d_Cooks_deformation} shows the membranes with different material discretizations at the steady state along with pointwise values of the Jacobian determinant of the non-local deformation gradient tensor at each material point. 
Fig.~\ref{f:2d_Cooks_disp} shows the $y$-displacement of the top-right corner (highlighted in Fig.~\ref{f:2d_Cooks_schematics}) after deformations for various PD horizon sizes under grid refinement. 
All three discretizations give results comparable to the result obtained using the IFED method, and they converge under grid refinement to approximately $0.67 \, \text{cm}$.

\begin{figure}[t!]
\centering
    \begin{tabular}{cc}
         \begin{subfigure}{.4\textwidth}
          		\includegraphics[width=\textwidth]{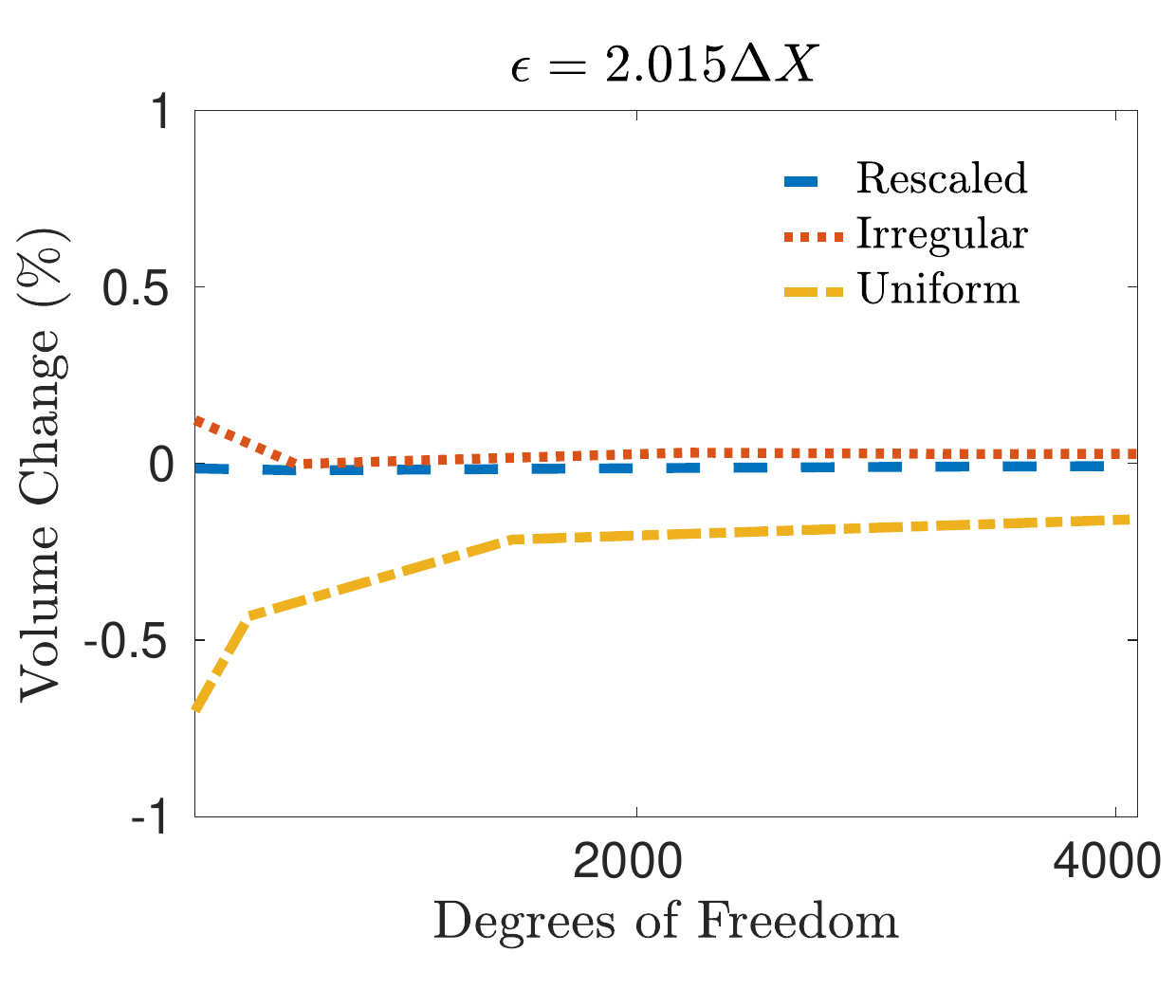}
        \end{subfigure}
        \hspace{.03\textwidth}
        \begin{subfigure}{.4\textwidth}
                \includegraphics[width=\textwidth]{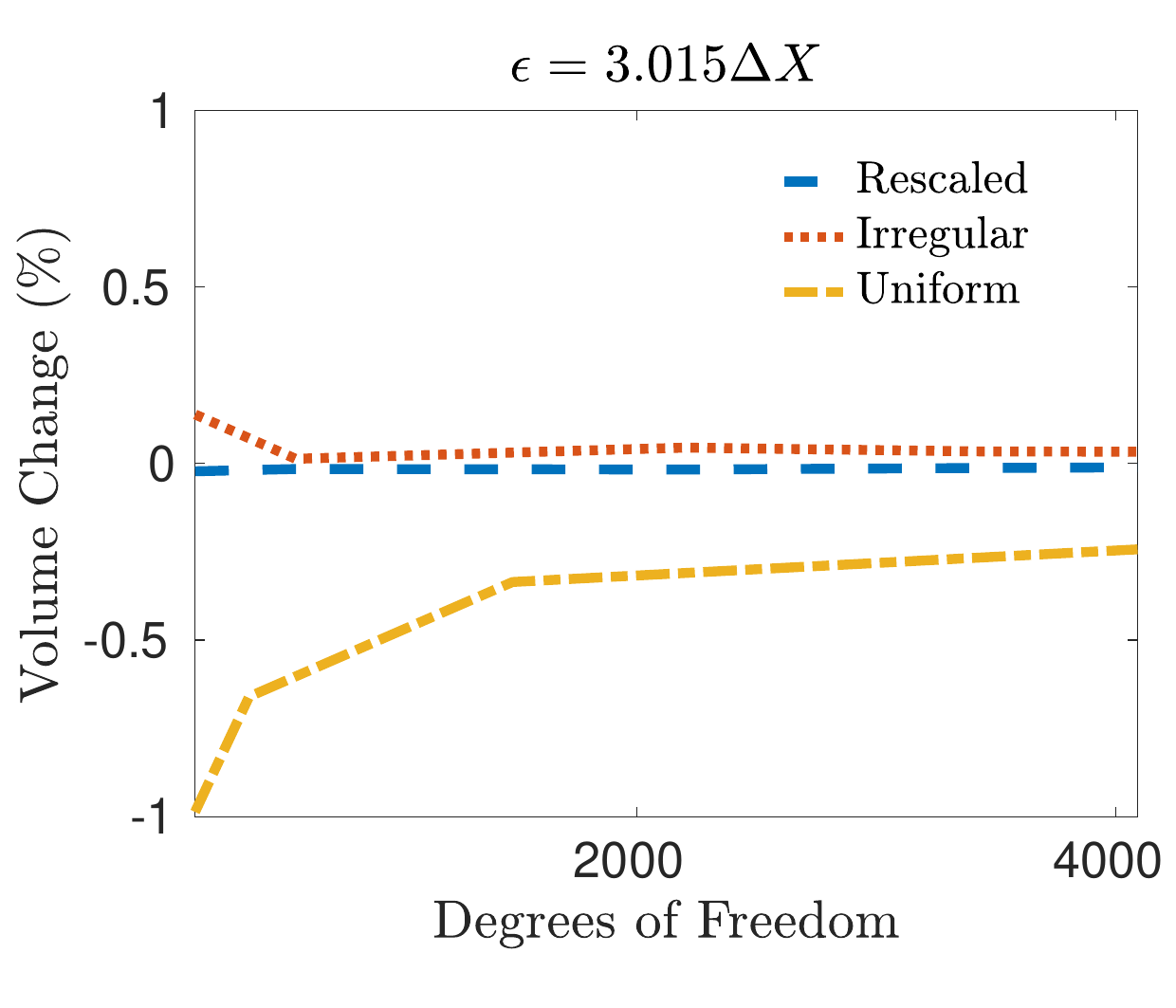}
        \end{subfigure}
    \end{tabular}
     \caption{Volume change of the Cook's membrane benchmark for different choices of discretization and horizon size $\horizonsize$. The numerical Poisson's ratio is set to $\nu_{\mathrm{stab}}  = 0.4$. The solid DoF range from $156$ to $4096$. The largest change is approximately $0.1 \%$ in the rescaled and irregular discretizations.}
    \label{f:2d_Cooks_vol}
\end{figure}

Fig.~\ref{f:2d_Cooks_vol} shows volume changes observed under deformation for different discretizations. 
With a smaller number of degrees of freedom, slight volume changes are observed under loading. 
The volume change becomes negligible (approximately $0.001\%$) under grid refinement in this benchmark test. 
The rescaled and irregular discretizations show better volume conservation as compared to the uniform discretization. 
This is also clear in Fig.~\ref{f:2d_Cooks_deformation}. 
However, the volume change using the uniformly distributed volumes is also comparable to the results obtained using IFED (between $0.000021\%$ and $0.1\%$). 
In addition, relatively consistent results are obtained for all considered PD horizon sizes.

\subsection{Anisotropic compressed block}
\label{s:aniso_Benchmark_Compression}
\begin{figure}[t!]
\centering
    \includegraphics[width=.46\textwidth]{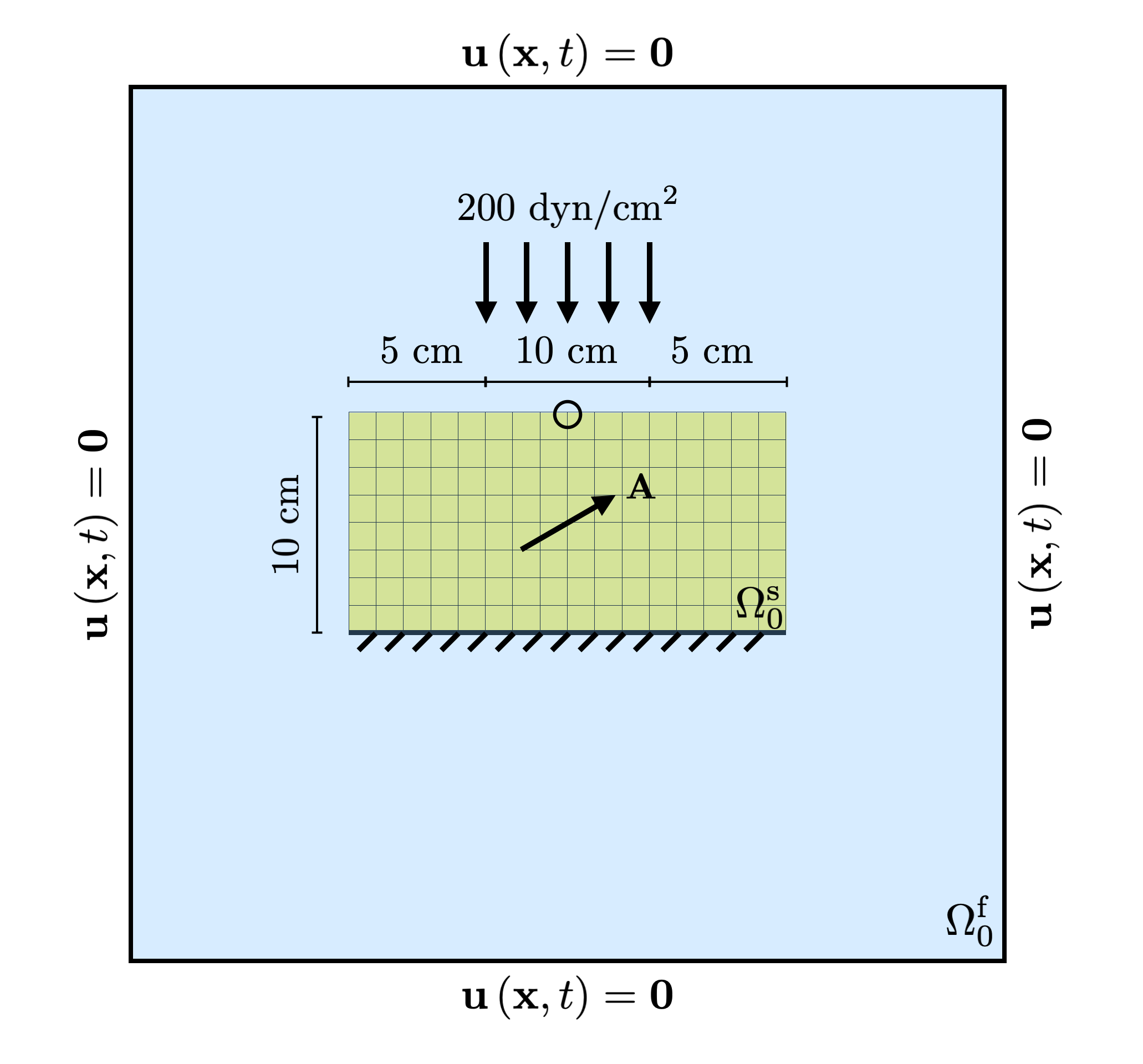}
    \caption{Schematic diagram for the anisotropic compression test. The initial configurations of the immersed structure and a fluid are denoted by $\Omega_0^{\text{s}}$ and $\Omega_0^{\text{f}}$, respectively. The entire computational domain is $\Omega = \Omega_0^{\text{s}} \cup \Omega_0^{\text{f}}$. Zero fluid velocity is enforced on the outer boundaries of the computational domain. The fiber direction is set to $\a = \left( \frac{\sqrt{3}}{2}, \frac{1}{2} \right)$.}
    \label{f:aniso_Compression_schematics}
\end{figure}
 
 We simulate the compression of a rectangular block to investigate the anisotropic hyperelastic material response under plane strain.
The isotropic version of this test was originally introduced by Reese et al.~\cite{reese1999new} to test a stabilization technique for low-order finite elements, and later this compression test was modified to simulate anisotropic hyperelastic material responses by Thekkethil et al.~\cite{THEKKETHIL2023115877}. 
 
The computational domain is defined as $\Omega = \left[0, L\right]^2$, where $L = 40 \, \text{cm}$. 
A downward uniaxial traction of $200 \, \frac{\text{dyn}}{\text{cm}^2}$ is applied to the center of the top surface.
The load time is $T_{\text{l}} = 100 \, \text{s}$ and the final time is $T_{\text{f}} = 500 \, \text{s}$.
The damping coefficient is set to $\eta = 4.0097 \, \frac{\text{g}}{\text{cm}^3 \cdot \text{s}}$.  
We constrain the horizontal displacement along the top boundary and the vertical displacement along the bottom boundary, while maintaining zero traction on all other boundaries. 
A schematic of the problem setup is illustrated in Fig.~\ref{f:aniso_Compression_schematics}.
To demonstrate the correspondence to benchmark IFED results, material failure (i.e., bond breakage) is not allowed.
 
For anisotropy, we consider a modified standard reinforced model
\begin{align}\label{standar_reinforced}
\Psi &= \Psi_{\text{isotropic}} + \Psi_{\text{anisotropic}},\\
\Psi_{\text{isotropic}} &= \frac{G}{2} (\bar{I}_1 - 2) + \frac{\kappa_{\text{stab}}}{2} \left( \ln J \right)^2,\\
\Psi_{\text{anisotropic}} &= \frac{G_{\text{f}}}{2} \left(I_4  -1 \right)^2, 
\end{align} 
with a shear modulus of $G = G_{\text{f}} = 80.194 \, \frac{\text{dyn}}{\text{cm}^2}$, the fourth invariant $I_4 =  \max\left( 1, \mathbb{C} : \left( \a \otimes \a  \right) \right)$, and the initial fiber direction $\a = \left( \frac{\sqrt{3}}{2}, \frac{1}{2} \right)$.
The standard reinforced model for a stabilized finite element method~\cite{THEKKETHIL2023115877} uses $\bar{\FF} = J^{-1/3} \, \FF$, but we use $\bar{\FF} = J^{-1/2} \, \FF$ for two-dimensional computations.

 \begin{figure}[t!]
\centering
	\begin{tabular}{cc}
	 \includegraphics[width=\textwidth]{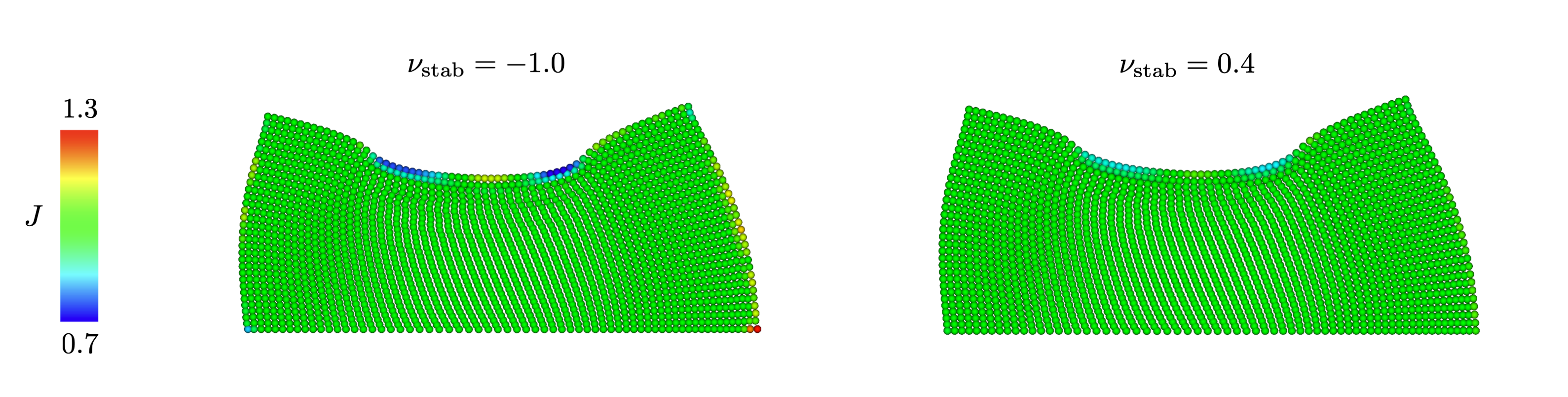}
	\end{tabular}
	    \caption{Deformations of the hyperelastic block along with the values of $J$ at material points using the standard fiber-reinforced model with $G = G_{\text{f}} = 80.194 \, \frac{\text{dyn}}{\text{cm}^2}$. The deformations are computed using $2145$ solid DoF and $\horizonsize = 2.015 \Delta X$. The left panel shows the deformation obtained using $\nu_{\mathrm{stab}}  = -1.0$, and the right panel shows the result for $\nu_{\mathrm{stab}}  = 0.4$.}
    \label{f:aniso_Compression_deformation}
\end{figure} 
 
\begin{figure}[t!]
\centering
    \begin{tabular}{cc}
    \hspace{-.1in}
        \begin{subfigure}{.4\textwidth}
          		\includegraphics[width=\textwidth]{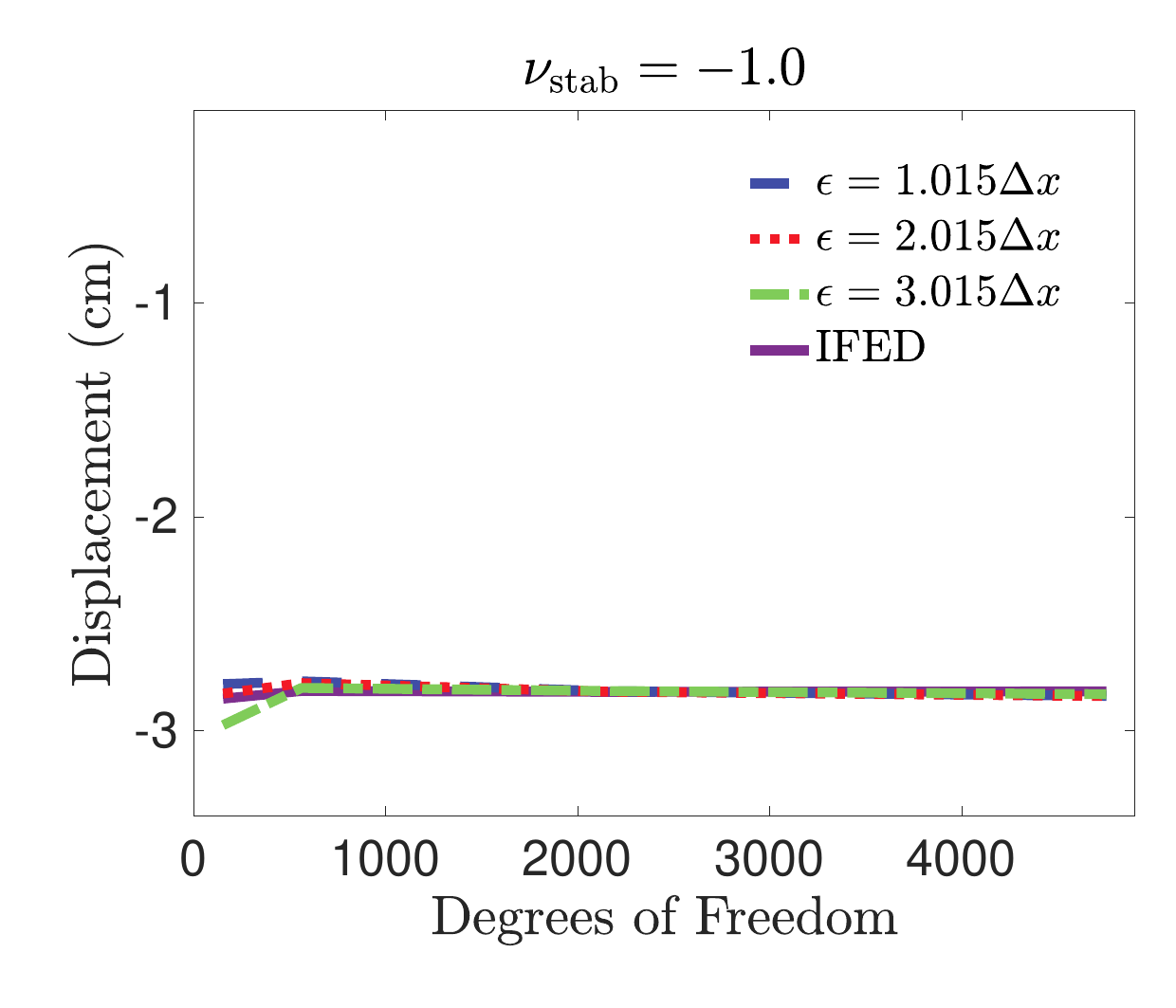}
        \end{subfigure}
        \hspace{.03\textwidth}
        \begin{subfigure}{.4\textwidth}
                \includegraphics[width=\textwidth]{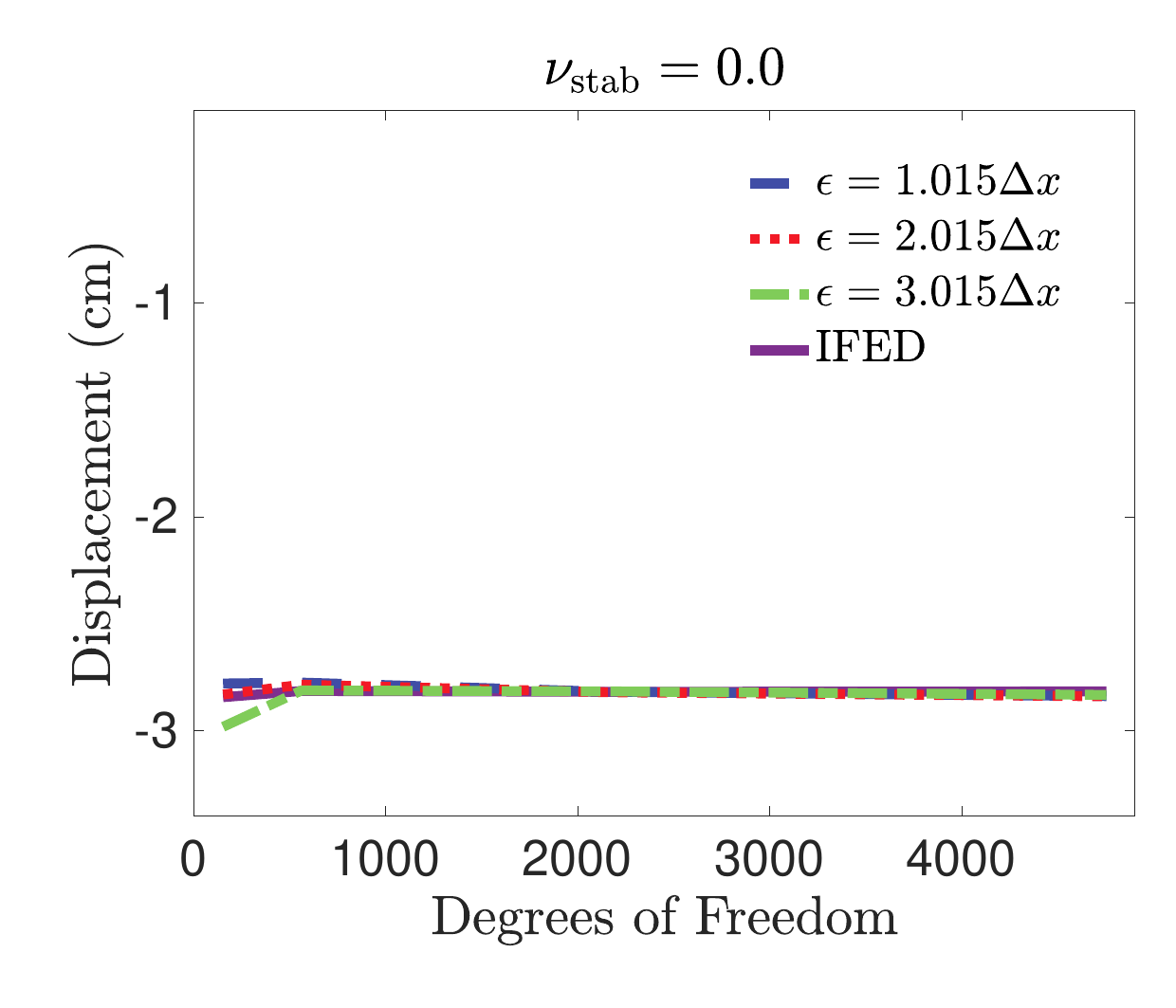}
        \end{subfigure} \\
        \begin{subfigure}{.4\textwidth}
                \includegraphics[width=\textwidth]{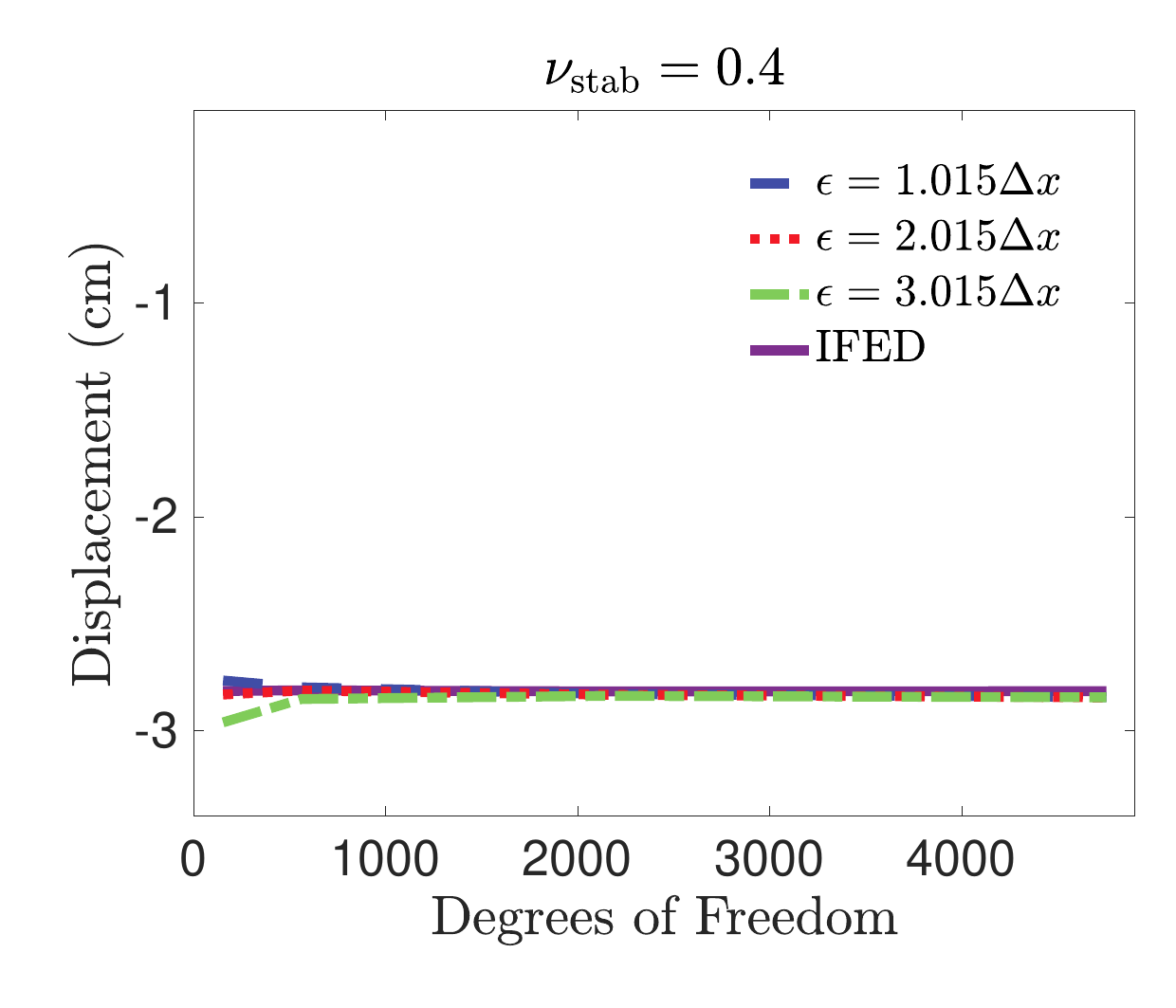}
        \end{subfigure}
        \hspace{.03\textwidth}
        \begin{subfigure}{.4\textwidth}
                \includegraphics[width=\textwidth]{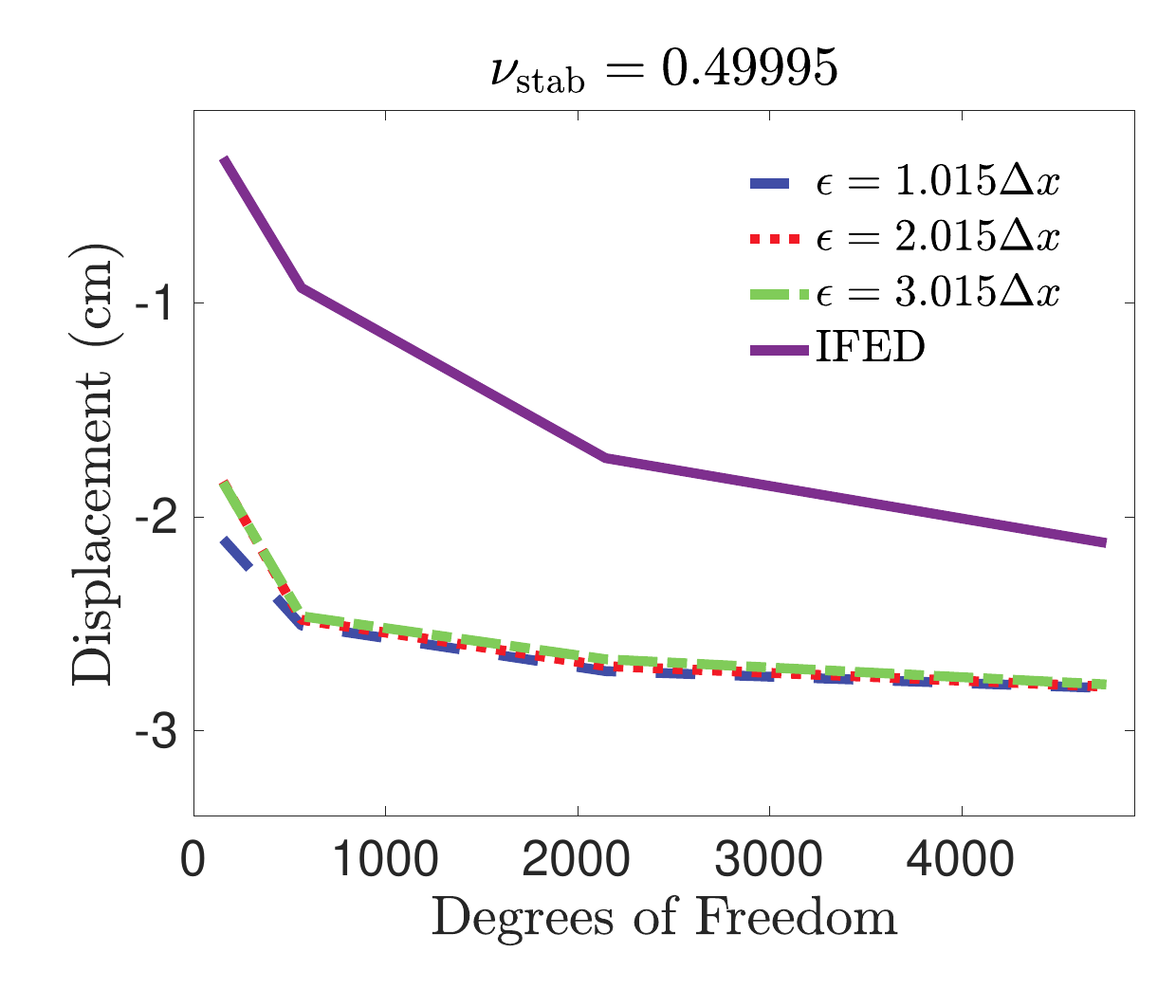}
        \end{subfigure}
    \end{tabular}
    \caption{Vertical displacements of the top center point of the compressed block, highlighted in Fig.~\ref{f:aniso_Compression_schematics}, for different choices of peridynamic horizon size $\horizonsize$ and numerical Poisson's ratio $\nu_{\mathrm{stab}}$. The solid DoF range from $153$ to $4753$.}
    \label{f:aniso_Compression_disp}
\end{figure}

Fig.~\ref{f:aniso_Compression_deformation} illustrates the material body after the deformation along with the pointwise Jacobian determinant $J$. 
Fig.~\ref{f:aniso_Compression_disp} shows the vertical displacements of the top center material point, highlighted in Fig.~\ref{f:aniso_Compression_schematics}, for various numerical Poisson's ratios $\nu_{\mathrm{stab}}$ and peridynamic horizon sizes $\horizonsize$ under grid refinement. 
The maximum displacement of the point obtained using the IPD method is in good agreement with that obtained using the IFED method, which converges under grid refinement to approximately $2.82 \, \text{cm}$.
The maximum displacement of the point of interest obtained is relatively small (between $1.83 \, \text{cm}$ and $2.48 \, \text{cm}$) at low grid resolutions if $\nu_{\mathrm{stab}}  = 0.49995$. 
This issue is commonly known as volumetric locking and occurs with large volumetric penalties in the computational mechanics literature, which is also shown in the results obtained by IFED.
Note that $\kappa_{\mathrm{stab}}  \rightarrow \infty$ as $\nu_{\mathrm{stab}} \rightarrow 0.5$. 
Under grid refinement, we ultimately recover accurate deformations for fixed finite values of $\kappa_{\text{stab}}$, as in standard methods for nearly incompressible elasticity.

\begin{figure}[t!]
\centering
    \begin{tabular}{cc}
    \hspace{-.2in}
        \begin{subfigure}{.4\textwidth}
          		\includegraphics[width=\textwidth]{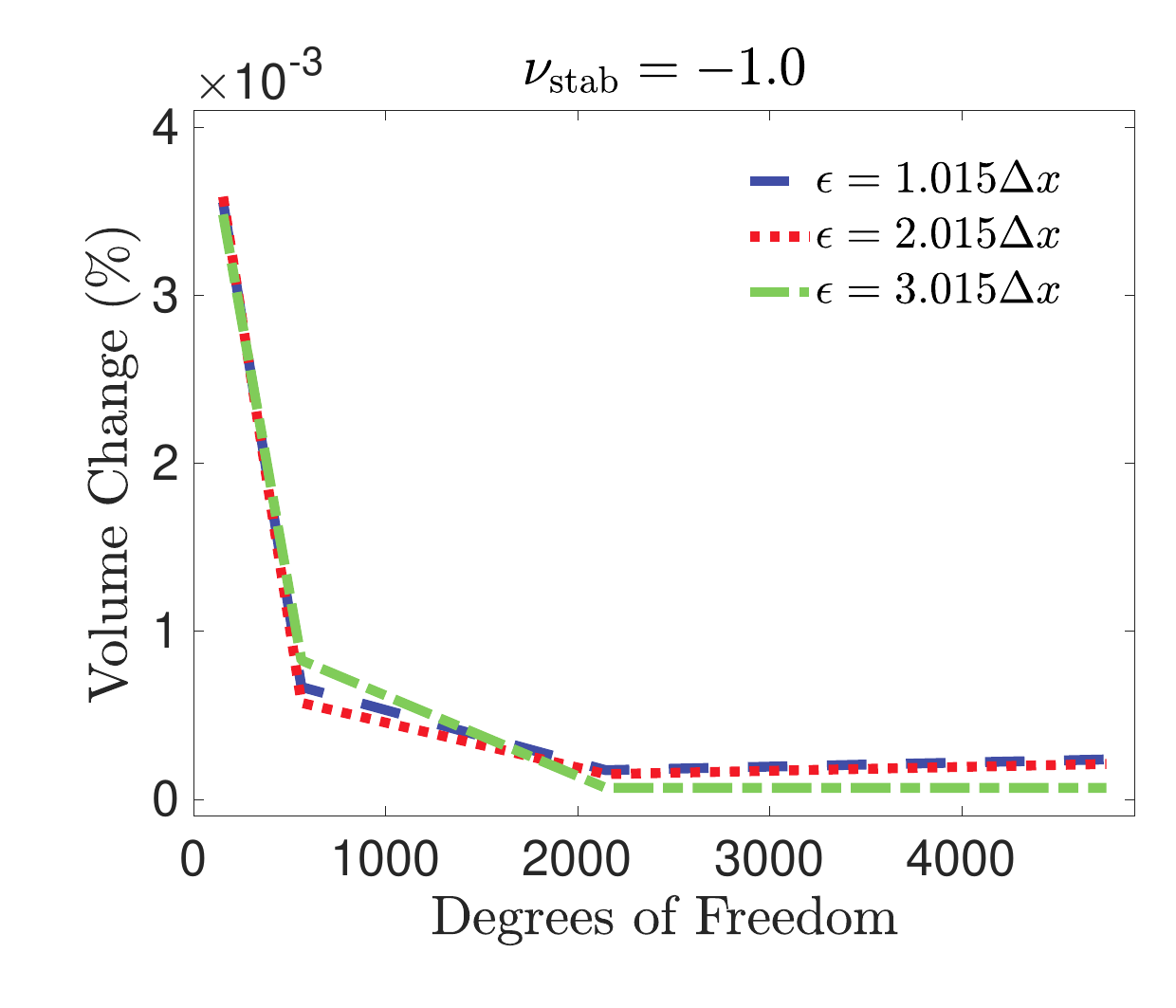}
        \end{subfigure}
        \hspace{.03\textwidth}
        \begin{subfigure}{.4\textwidth} 
                \includegraphics[width=\textwidth]{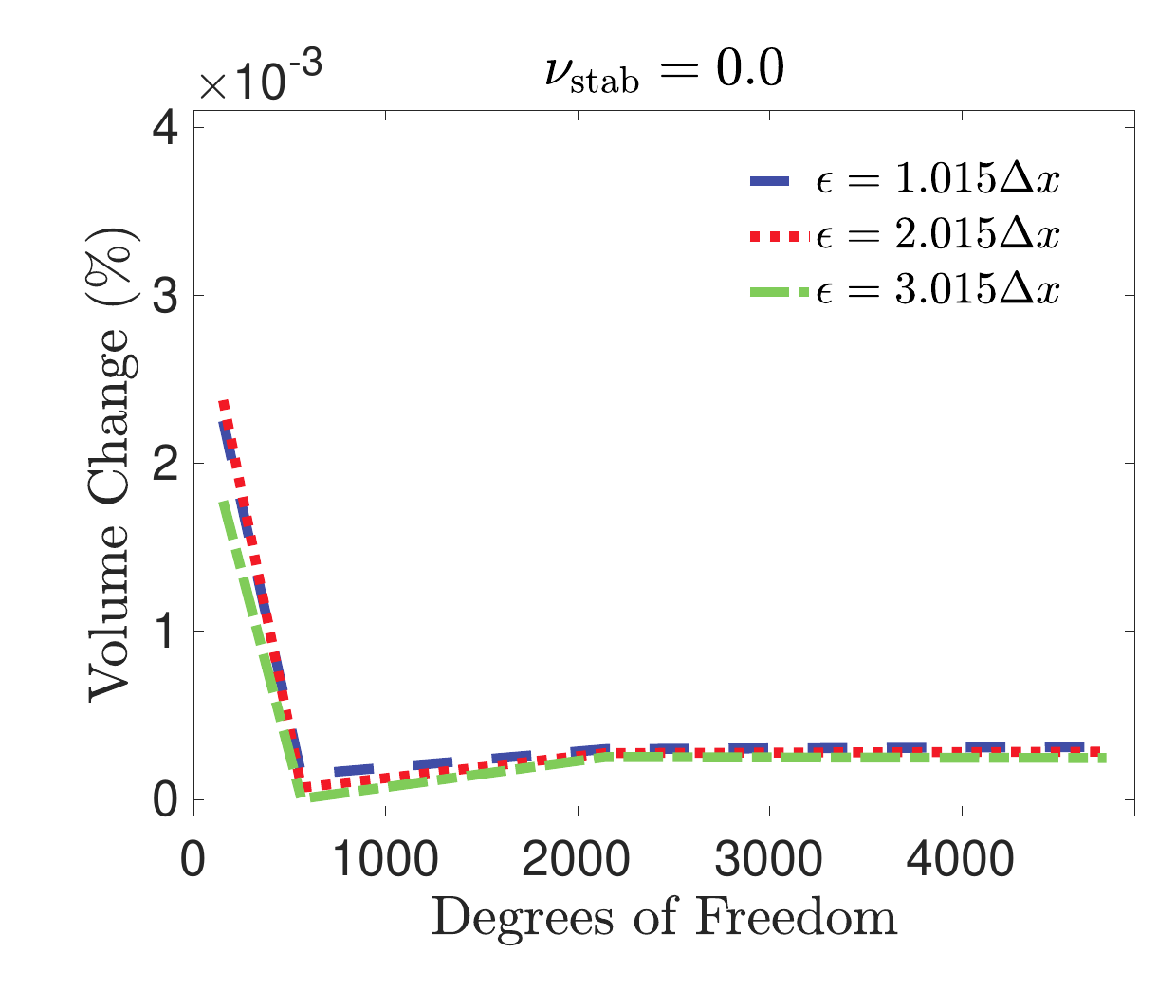}
        \end{subfigure} \\
        \begin{subfigure}{.4\textwidth}
                \includegraphics[width=\textwidth]{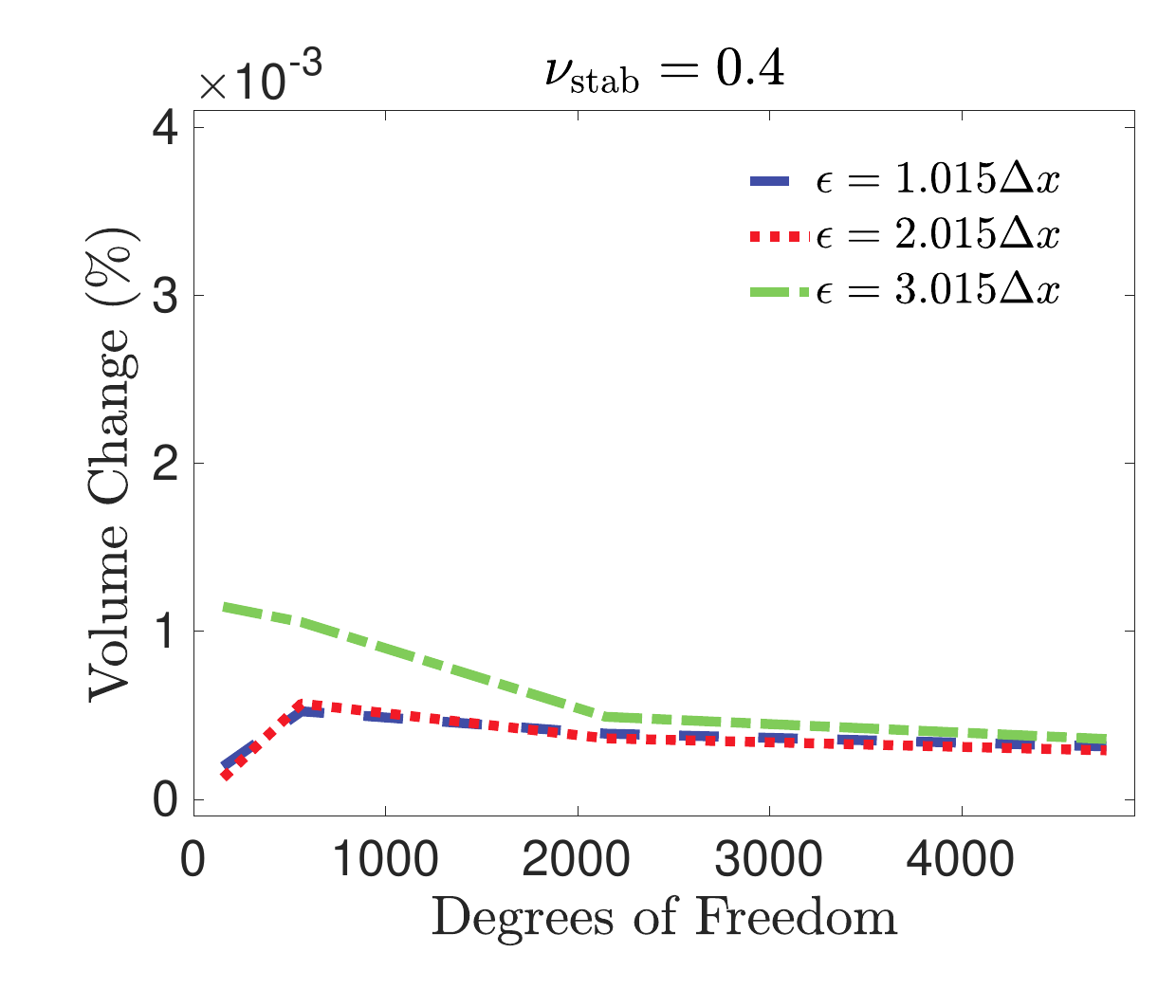}
        \end{subfigure}
        \hspace{.03\textwidth}
        \begin{subfigure}{.4\textwidth}
                \includegraphics[width=\textwidth]{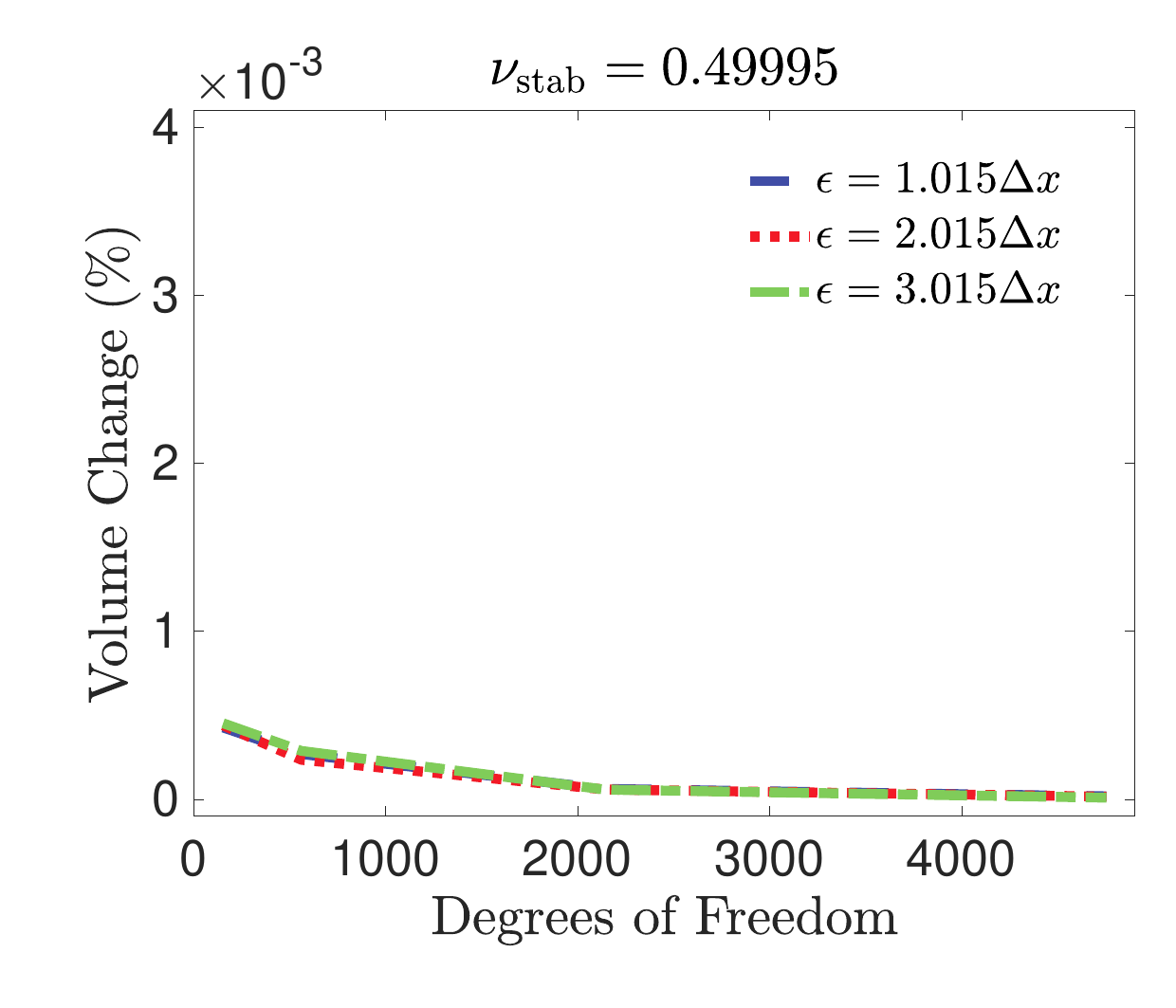}
        \end{subfigure}
    \end{tabular}
    \caption{Volume change of the compressed block for different choices of horizon size $\horizonsize$ and numerical Poisson's ratio $\nu_{\mathrm{stab}}$. The solid DoF range from $153$ to $4753$. The largest volume change is approximately $0.0036\%$.}
    \label{f:aniso_Compression_vol}
\end{figure}

Fig.~\ref{f:aniso_Compression_vol} shows the volume change observed under deformation. 
In contrast to the results from the prior study of isotropic material responses, the volume change with a lower Poisson's ratio is relatively small due to anisotropy along the immersed body and becomes negligible (up to $0.00001 \%$).
This is also clear in Fig.~\ref{f:aniso_Compression_deformation}. 
In addition, relatively consistent results are obtained for all choices of the PD horizon sizes considered in the tests.

\begin{figure}[t!]
\centering
	\begin{tabular}{cc}
	 \includegraphics[width=0.55\textwidth]{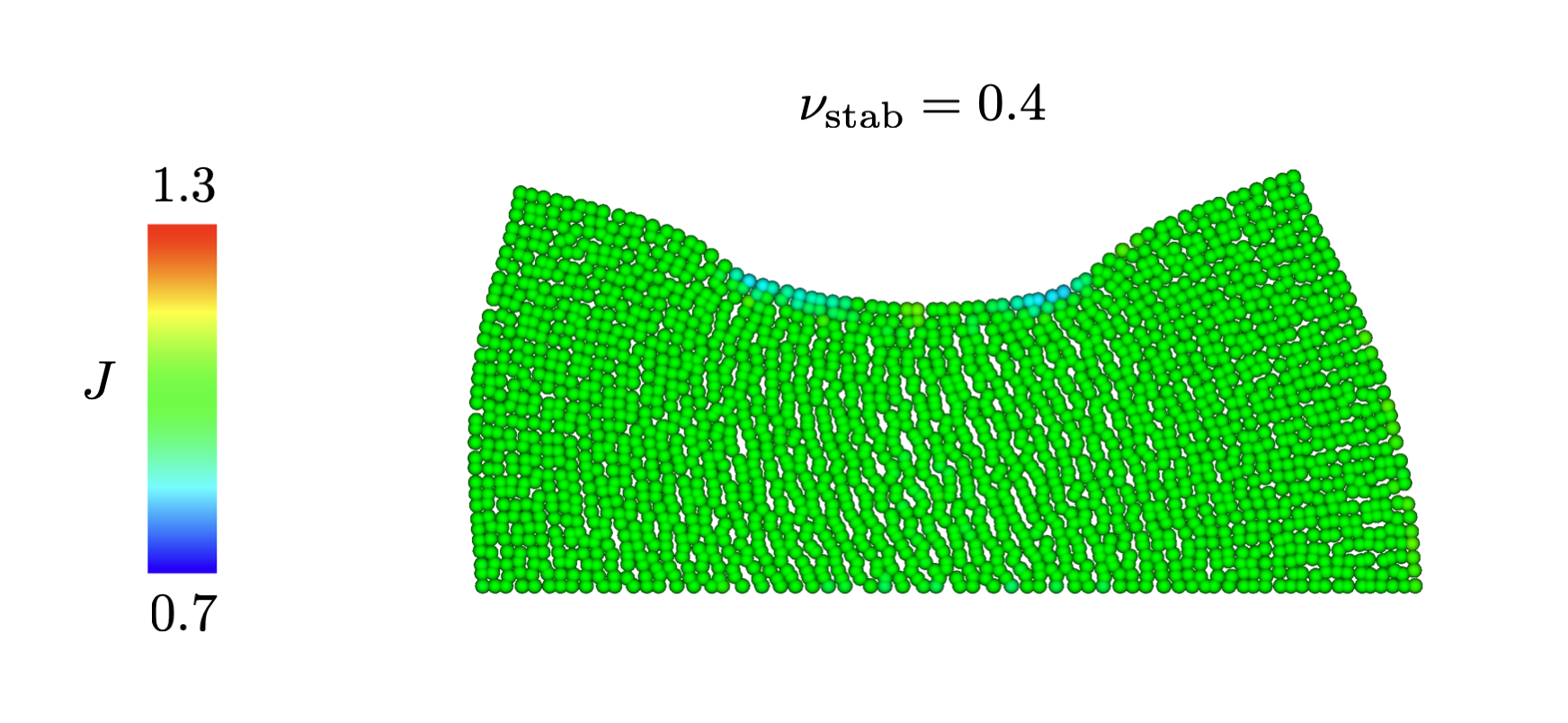}
	\end{tabular}
	    \caption{A deformation of the hyperelastic block along with the values of $J$ at material points using the standard fiber-reinforced model with $G = G_{\text{f}} = 80.194 \, \frac{\text{dyn}}{\text{cm}^2}$. The deformation is computed using irregularly distributed $2145$ solid DoF with $\horizonsize = 2.015 \Delta X$ and $\nu_{\mathrm{stab}}  = 0.4$.}
    \label{f:aniso_Comp_irr}
\end{figure}

\begin{figure}[t!]
\centering
    \begin{tabular}{cc}
    \hspace{-.2in}
        \begin{subfigure}{.4\textwidth}
          		\includegraphics[width=\textwidth]{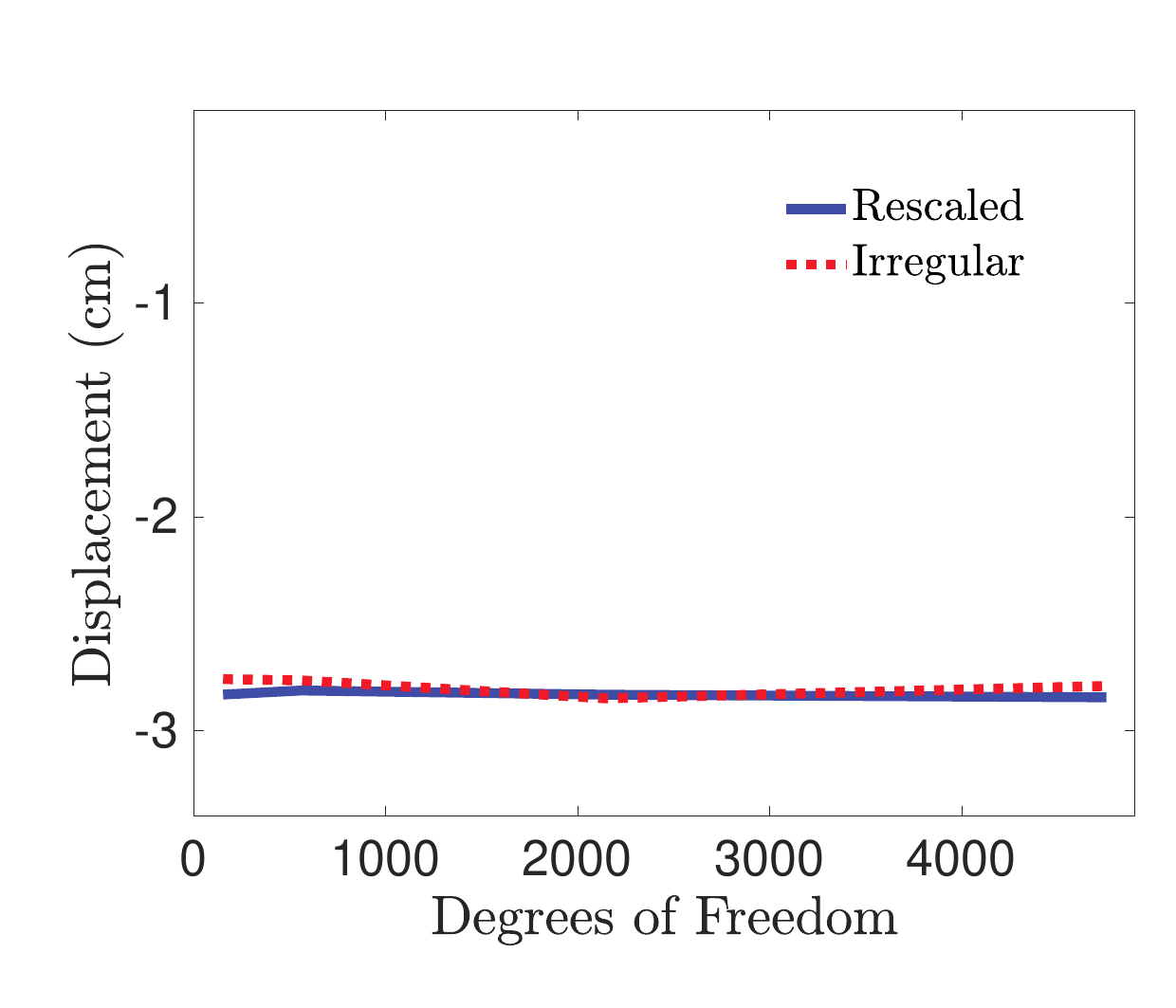}
          		\caption{}
        \end{subfigure}
        \hspace{.03\textwidth}
        \begin{subfigure}{.4\textwidth} 
                \includegraphics[width=\textwidth]{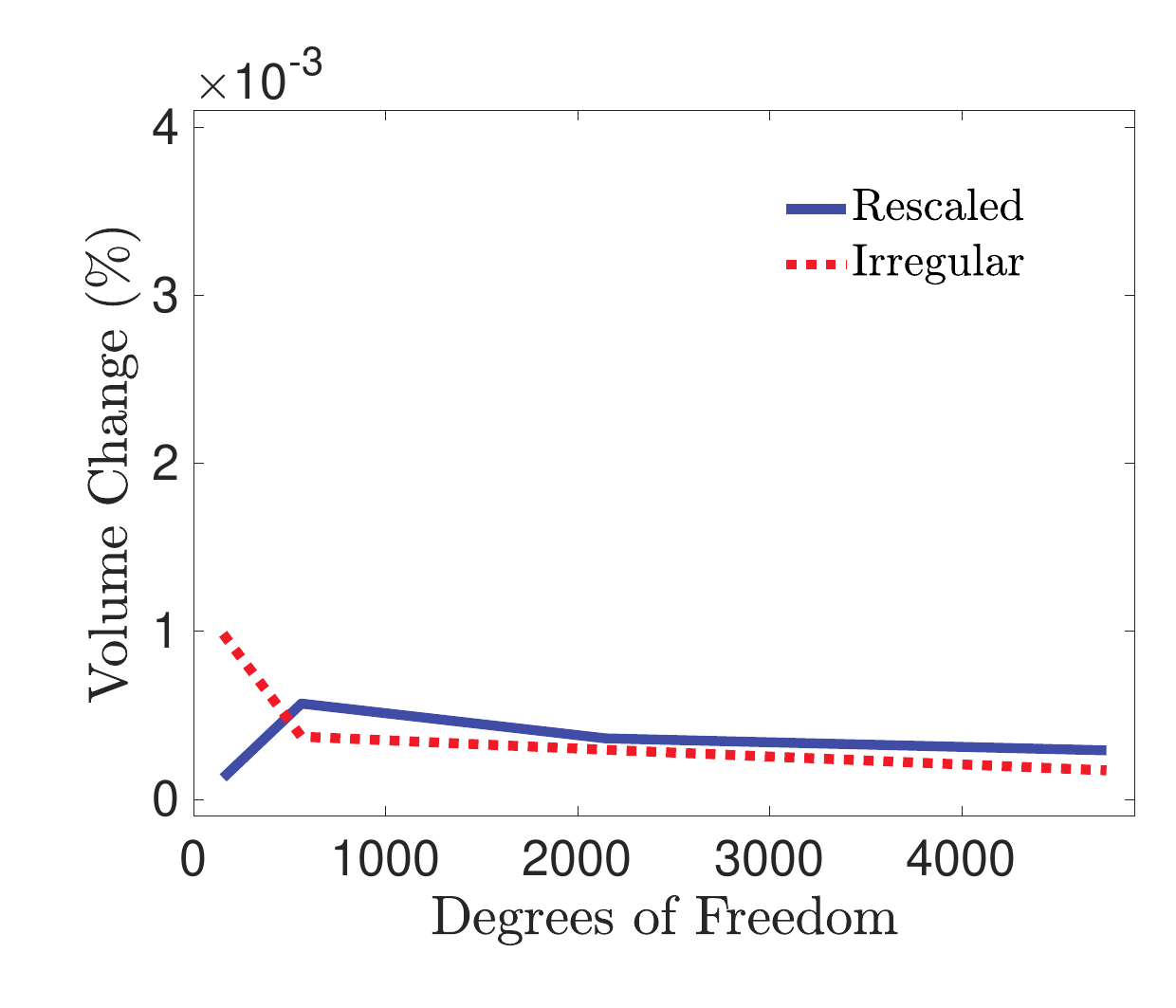}
                \caption{}
        \end{subfigure} 
    \end{tabular}
    \caption{Comparison between two structural point distributions with $\horizonsize = 2.015 \Delta X$ and $\nu_{\mathrm{stab}} = 0.4$. The number of solid DoF ranges from $153$ to $4753$. (a) Vertical displacement of the top center point of the compressed block, highlighted in Fig.~\ref{f:aniso_Compression_schematics}. (b) Volume change of the compressed block. }
    \label{f:aniso_Comp_irr_disp_vol}
\end{figure}

To evaluate the IPD performance with an irregular point cloud, we randomly perturb the initial point cloud and perform the same compression test using $\nu_{\mathrm{stab}} = 0.4$ and $\horizonsize = 2.015 \Delta X$.
Interior points are perturbed in both coordinate directions. 
To preserve an exact representation of the material boundary, each noncorner boundary point is perturbed only in the direction tangent to its corresponding boundary, so that it remains aligned with that boundary. 
The four corner points and the top-center point used to measure displacement are left unperturbed.
Fig.~\ref{f:aniso_Comp_irr} illustrates the material body after the deformation along with the pointwise Jacobian determinant $J$. 
Fig.~\ref{f:aniso_Comp_irr_disp_vol} compares the steady-state displacements and volume changes under grid refinement against the results obtained using the rescaled discretization.
Note that the rescaled and uniform distributions are equivalent in this case.
We demonstrate that the accuracy achieved by the rescaled distribution is comparable to that of the irregular distribution, as also shown for the benchmark in Sec.~\ref{s:2d_cooks_discretizations}.
Therefore, we utilize a rescaled point distribution for the remainder of this study unless an irregular distribution is specifically required.

\subsection{Three-dimensional anisotropic Cook's membrane}
\label{s:Aniso_Cooks}

\begin{figure}[t!]
\centering
    \includegraphics[width=.4\textwidth]{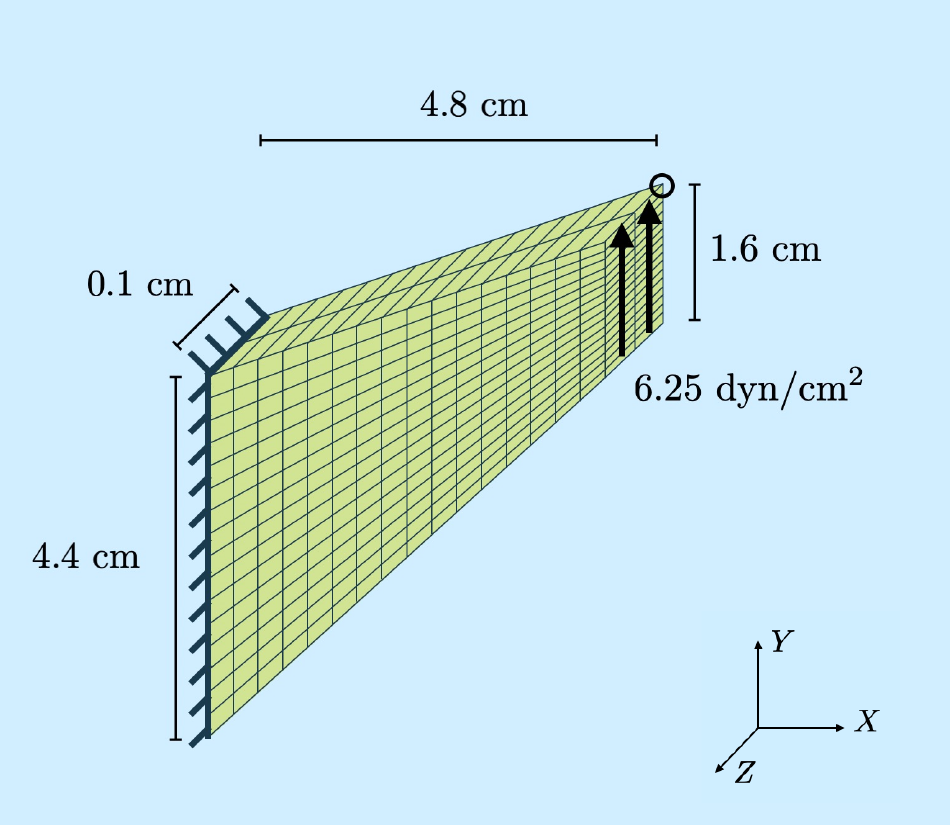}
    \caption{Schematic diagram for the three-dimensional anisotropic Cook's membrane. The computational domain is $\Omega = [0,L]^3$, with $L = 12  \, \text{cm}$. Zero fluid velocity is enforced on the outer boundaries of the computational domain. The fiber direction is set to $\a =\left(\frac{1}{\sqrt{3}},\frac{1}{\sqrt{3}},\frac{1}{\sqrt{3}} \right) $.}
    \label{f:aniso_cooks_schematics}
\end{figure}

We use a three-dimensional version of Cook's membrane~\cite{vadala2020stabilization} to investigate three-dimensional hyperelastic anisotropic material responses under bending and shearing. 
The computational domain is the cube $\Omega = \left[0, L\right]^3$, with $L = 12 \,  \text{cm}$. 
Zero fluid velocity is applied to the boundaries of the computational domain.
The density and viscosity of the surrounding fluid are set to $\rho = 1.0 \,  \frac{\text{g}}{\text{cm}^3}$ and  $\mu = 0.16 \, \frac{\text{dyn$\cdot$s}}{\text{cm}^2}$, respectively. 
A larger viscosity is chosen to accelerate reaching steady state.
An upward traction is set to $6.25 \, \frac{\text{dyn}}{\text{cm}^2}$ on the right face, and a zero displacement condition is applied to the left surface.
Otherwise, stress-free boundary conditions are prescribed.
The load time is $T_{\text{l}} = 14 \, \text{s}$, the final time is $T_{\text{f}} = 35 \, \text{s}$, and the damping parameter is set to zero.
Fig.~\ref{f:aniso_cooks_schematics} provides a schematic of this test case. 
To verify the correspondence to benchmark IFED results, structural damage and failure are not allowed.
Additionally, two-dimensional Cook's membrane tests in Sec.~\ref{s:2d_cooks_discretizations} indicate that rescaled and irregular discretizations give good agreement with classical hyperelasticity compared to the uniform geometry. Thus, we only use the rescaled material discretization for the rest of the numerical tests.

For anisotropy, we use a transversely isotropic material model~\cite{vadala2020stabilization}
\begin{align}
\Psi &= \Psi_{\text{isotropic}} + \Psi_{\text{anisotropic}},\\
\Psi_{\text{isotropic}} &= \frac{G_{\text{T}}}{2} \left( \bar{I}_1 - 3 \right) + \frac{\kappa_{\text{stab}}}{2} \left( \ln J \right)^2,\\
\Psi_{\text{anisotropic}} &= \frac{G_{\text{T}} - G_{\text{L}}}{2} \left( 2I_4 - I_5 - 1 \right) + \frac{E_{\text{L}} + G_{\text{T}} - 4 G_{\text{L}}}{8} \left(I_4 - 1 \right)^2,
\end{align}
in which $G_{\text{T}}$ is the shear modulus of the material in the plane transverse to the fibers, $G_{\text{L}}$ is the shear modulus along the fibers, $E_{\text{L}}$ is a Young's modulus in the fiber direction, $I_4 =  \max\left( 1, \mathbb{C} : \left( \a \otimes \a  \right) \right)$ is the fourth invariant, $\a$ is the initial fiber direction, and $I_5 = \a^{\tran} \mathbb{C}^2 \a$ is the fifth invariant, which encodes the information related to both shear and stretch.
We set $G_{\text{T}} = 8 \, \frac{\text{dyn}}{\text{cm}^2}$, $G_{\text{L}} = 160 \, \frac{\text{dyn}}{\text{cm}^2}$, $E_{\text{L}} = 1200 \, \frac{\text{dyn}}{\text{cm}^2}$, and $\a =\left(\frac{1}{\sqrt{3}},\frac{1}{\sqrt{3}},\frac{1}{\sqrt{3}} \right) $.

Since overall deformations and volume conservation with the numerical Poisson's ratio of $0.4$ in Sec.~\ref{s:aniso_Benchmark_Compression} are comparable to the results achieved by IFED, we only consider two cases: $\nu_{\text{stab}} = 0.4$ for nearly incompressible elasticity and $\nu_{\text{stab}} = -1$ for comparison.

 \begin{figure}[t!]
\centering
	\includegraphics[width=0.8\textwidth]{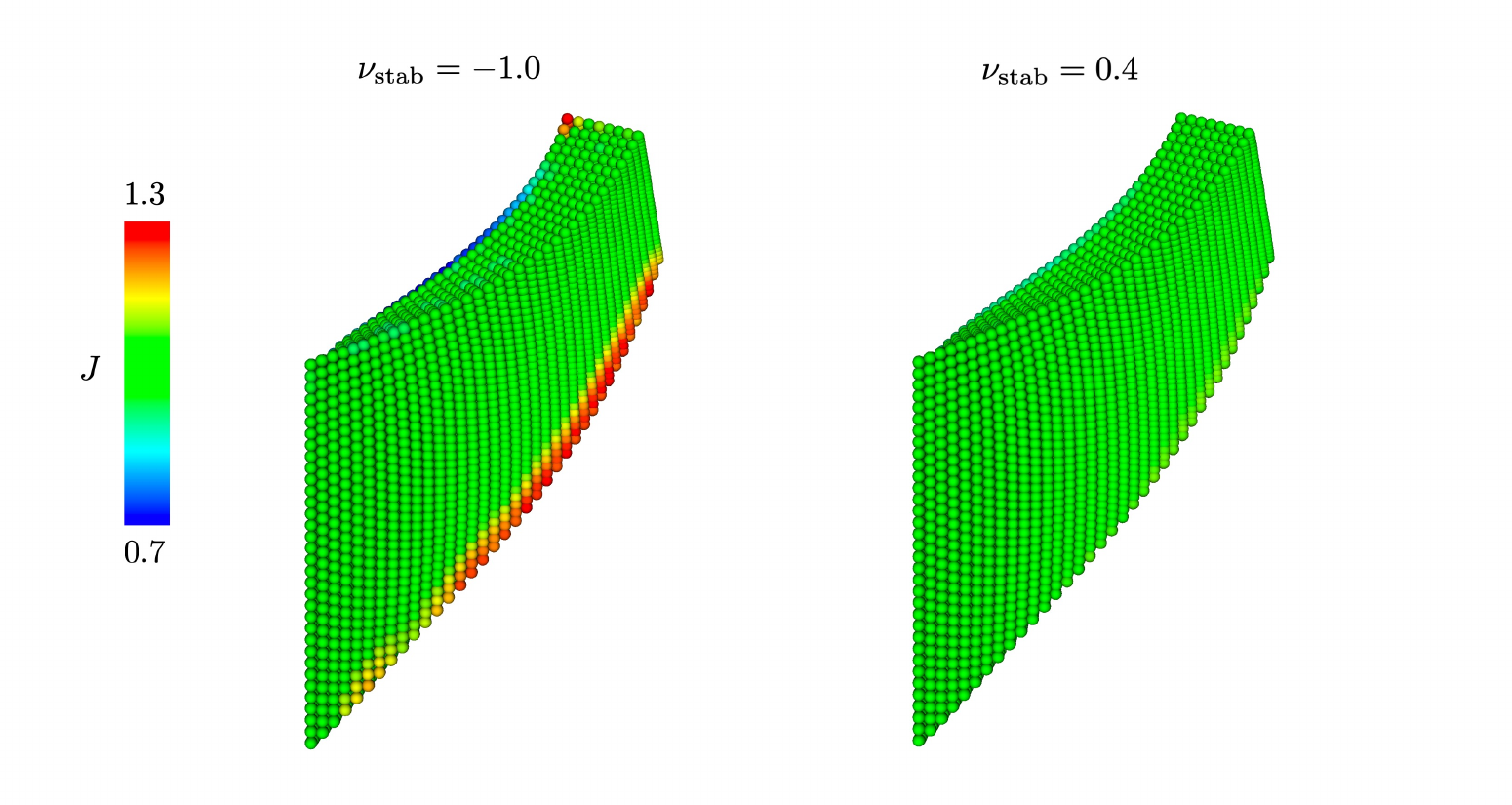} 
    \caption{Deformations of three-dimensional anisotropic Cook's membrane with the values of $J$ at material points using the transversely isotropic material with $G_{\text{T}} = 8 \, \frac{\text{dyn}}{\text{cm}^2}$, $G_{\text{L}} = 160 \, \frac{\text{dyn}}{\text{cm}^2}$, and $E_{\text{L}} = 1200 \, \frac{\text{dyn}}{\text{cm}^2}$. The deformations are represented using $3450$ solid DoF and $\horizonsize = 2.015 \Delta X$ with a non-uniform volume distribution.  }
    \label{f:Aniso_Cooks_deformation}
\end{figure}

\begin{figure}[t!]
\centering
    \begin{tabular}{cc}
        \begin{subfigure}{.4\textwidth}
          		\includegraphics[width=\textwidth]{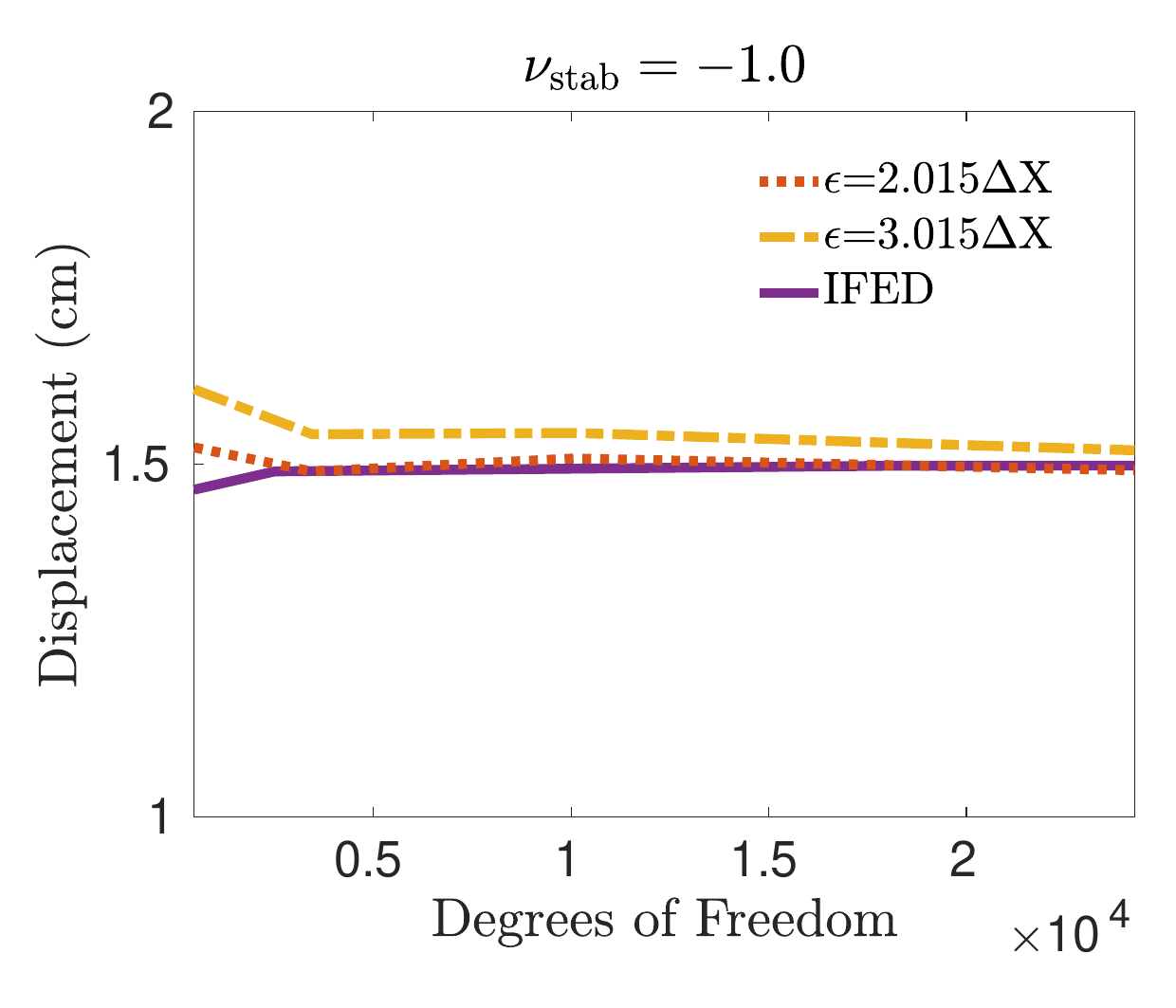}
        \end{subfigure}
        \hspace{.03\textwidth}
        \begin{subfigure}{.4\textwidth}
                \includegraphics[width=\textwidth]{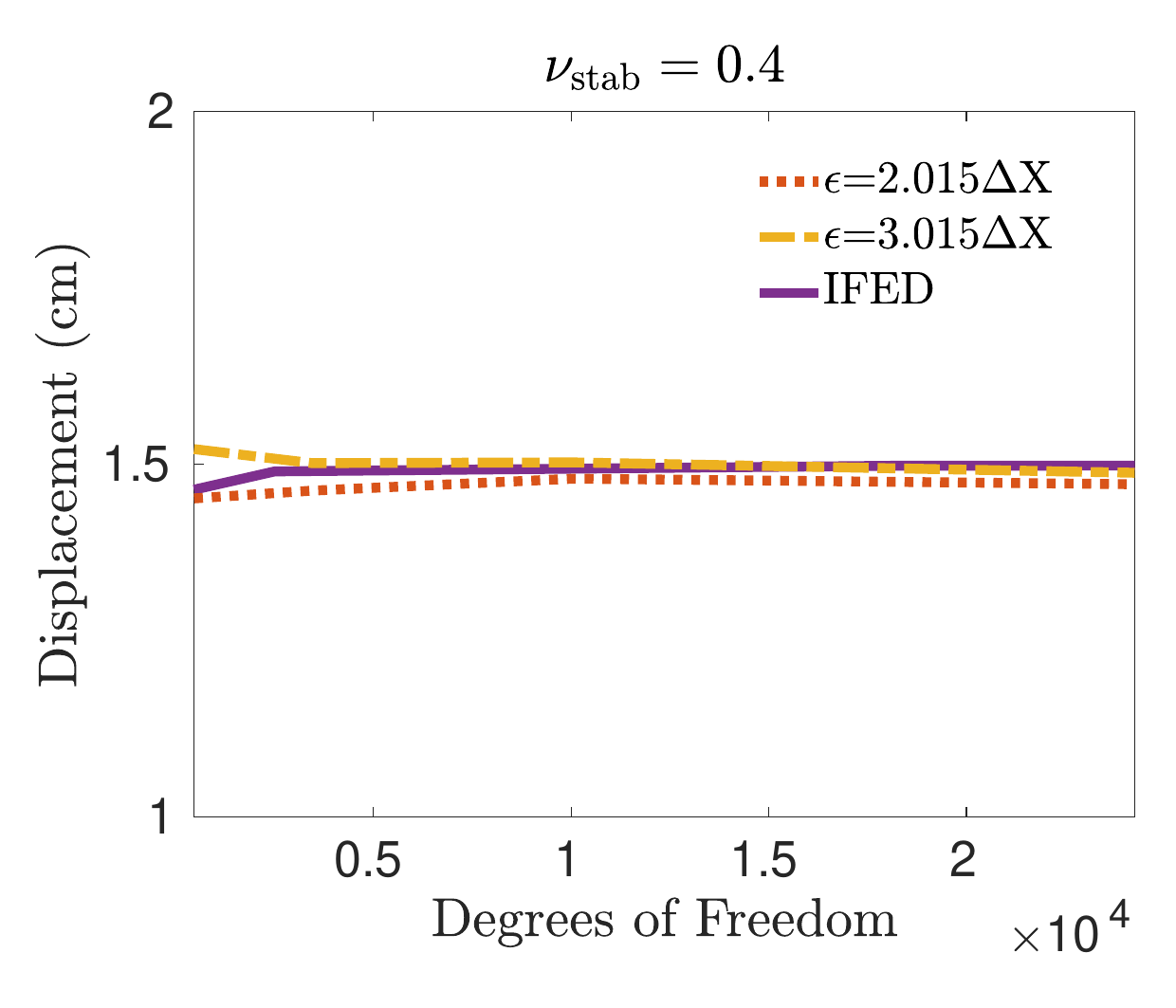}
        \end{subfigure}
    \end{tabular}
    \caption{Vertical displacements of the right top corner point of the Cook's membrane benchmark, highlighted in Fig.~\ref{f:aniso_cooks_schematics}, for different choices of peridynamic horizon size $\horizonsize$ and numerical Poisson's ratios. The solid DoF range from $468$ to $24255$. }
    \label{f:Aniso_Cooks_disp}
\end{figure}

Fig.~\ref{f:Aniso_Cooks_deformation} illustrates the deformations of the membrane under the traction loading along with the Jacobian determinant of the non-local deformation gradient tensor at each material point. 
Fig.~\ref{f:Aniso_Cooks_disp} shows the maximum displacement of the right corner of the right surface in Fig.~\ref{f:aniso_cooks_schematics} at the steady states for various sizes of $\horizonsize$ under grid refinement. 
The displacements using IPD are comparable to the results generated by the IFED method and converge under grid refinement to a value of approximately $1.49 \, \text{cm}$. 

\begin{figure}[t!]
\centering
    \begin{tabular}{cc}
         \begin{subfigure}{.4\textwidth}
          		\includegraphics[width=\textwidth]{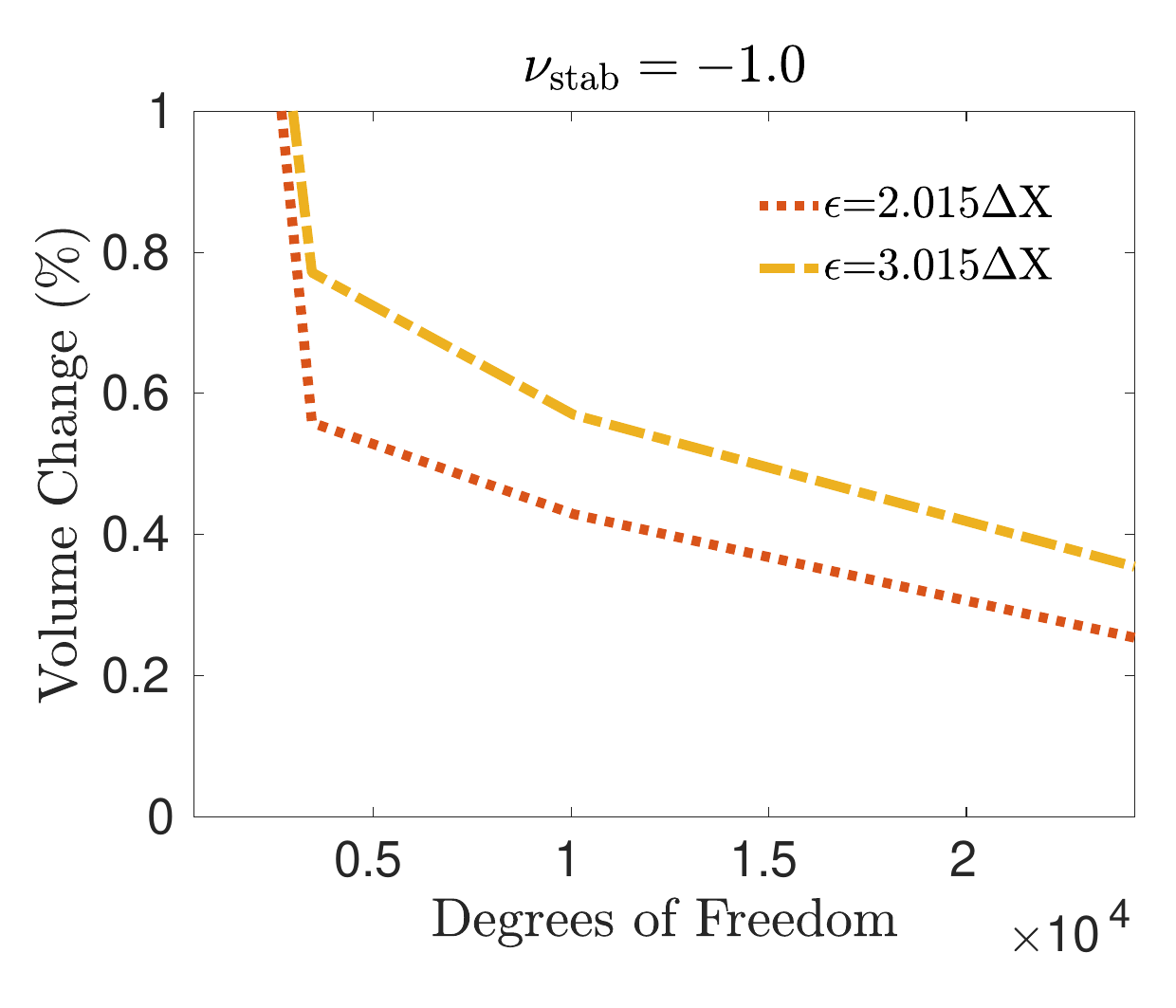}
        \end{subfigure}
        \hspace{.03\textwidth}
        \begin{subfigure}{.4\textwidth}
                \includegraphics[width=\textwidth]{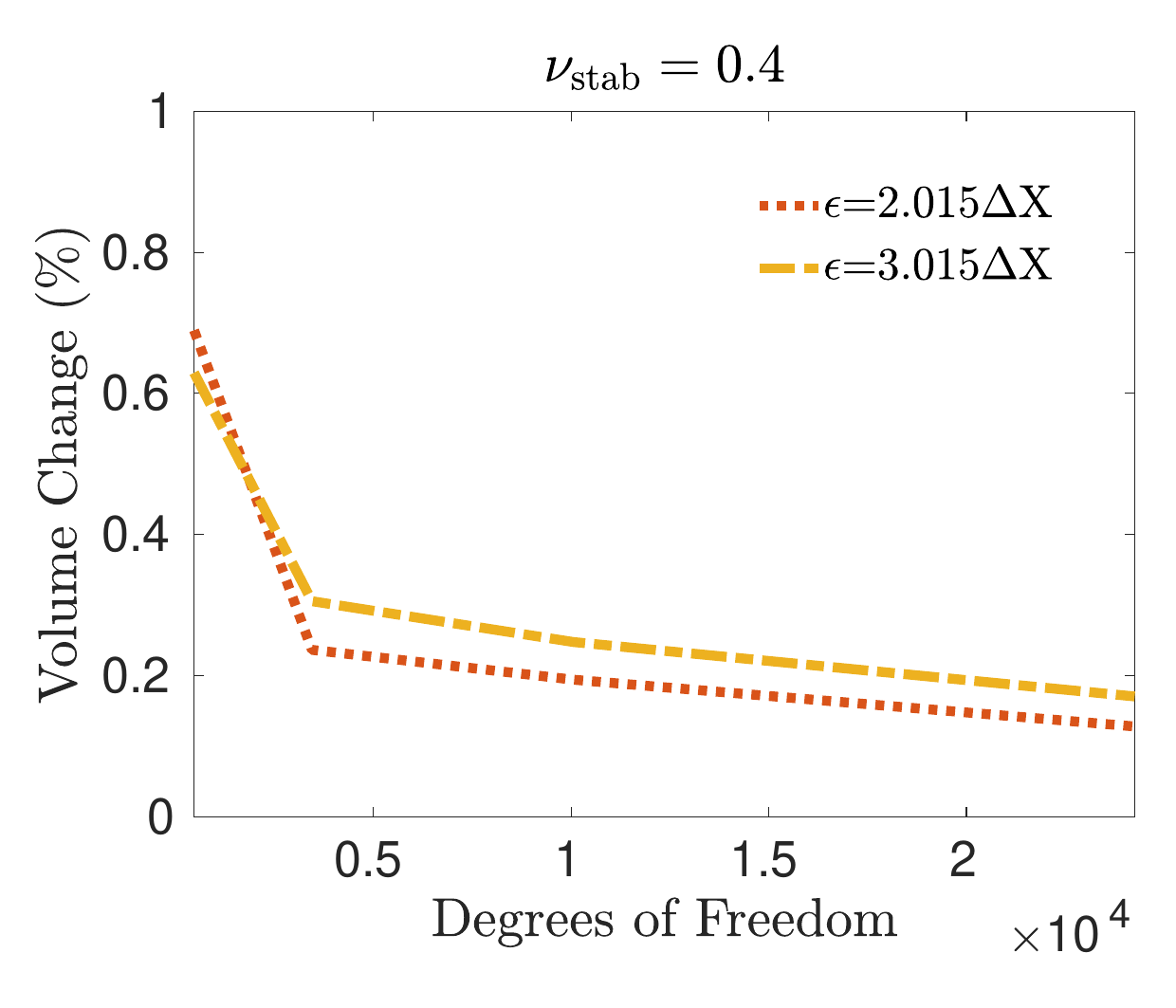}
        \end{subfigure}
    \end{tabular}
     \caption{Volume change of the Cook's membrane benchmark for different choices of horizon size $\horizonsize$ and numerical Poisson's ratios. The solid DoF range from $468$ to $24255$. The largest change is approximately $2.26\%$ in the non-uniform discretizations, but it decreases under grid refinement to as low as $0.12\%$.}
    \label{f:Aniso_Cooks_vol}
\end{figure}

Fig.~\ref{f:Aniso_Cooks_vol} shows the volume change of the membrane for different numbers of the solid DoF, ranging from $0.12 \%$ to $2.26\%$. 
The maximum volume change obtained by the IFED method is $2.1 \%$, which is comparable to the IPD simulation.

\subsection{Anisotropic torsion}
\label{s:aniso_Benchmark_Torsion}
\begin{figure}[t!]
\centering
    \includegraphics[width=.7\textwidth]{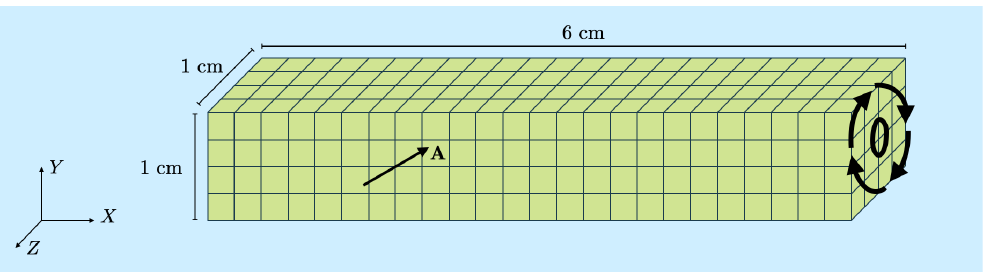}
    \caption{Schematic diagram for the three-dimensional anisotropic torsion benchmark. The computational domain is $\Omega = [0,L]^3$, with $L = 9  \, \text{cm}$, and the three-dimensional beam is placed at the center of the domain. Zero fluid velocity is enforced on the outer boundaries of the computational domain. The three-dimensional beam is reinforced with fibers and the fiber direction is set to $\a = \left(\frac{\sqrt{3}}{2}, \frac{1}{2}, 0 \right)$.}
    \label{f:aniso_torsion_schematics}
\end{figure}

We use a three-dimensional beam to investigate three-dimensional anisotropic material responses under torsion. 
The original isotropic version of this benchmark was suggested by Bonet et al.~\cite{bonet2015computational} and later modified to simulate three-dimensional anisotropic hyperelastic material responses by Thekkethil et al.~\cite{THEKKETHIL2023115877}.
We modify the three-dimensional torsion test in our previous work for isotropic materials by incorporating fiber reinforcement.

The computational domain is the cube $\Omega = \left[0, L\right]^3$, with $L = 9 \,  \text{cm}$. 
The density and viscosity of the surrounding fluid are set to $\rho = 1.0 \,  \frac{\text{g}}{\text{cm}^3}$ and  $\mu = 0.04 \, \frac{\text{dyn$\cdot$s}}{\text{cm}^2}$, respectively. 
One side of the beam is fixed in place, and torsion is applied to the opposite end through displacement boundary conditions. 
Fig.~\ref{f:aniso_torsion_schematics} illustrates the schematic of the test case.
The right surface is rotated by the linear function in time:
\begin{align}
\theta(t) = 
\begin{cases}
 \frac{t}{T_{\text{l}}} \theta_{T_{\text{f}}},   &\text{ if } t \le T_{\text{l}}, \\
 \theta_{T_{\text{f}}},  &\text{ otherwise},
 \end{cases}
\end{align} 
 in which $\theta_{T_{\text{f}}} = 2.5 \pi$.
The maximum angle of rotation is achieved at $T_{\text{l}} = 0.4 T_{\text{f}}$, with $T_{\text{f}} = 5 \,  \text{s}$. 
All other surfaces have zero traction boundary conditions. 
Material damage and failure are not allowed. 

We use a standard reinforced material model for anisotropy and a modified Mooney--Rivlin material model~\cite{vadala2020stabilization} for the hyperelastic isotropic response of a ground material:
\begin{align}
\Psi  &= \Psi_\text{isotropic} + \Psi_\text{anisotropic}, \\
\Psi_\text{isotropic} &= c_1\left( \bar{I}_1 - 3 \right) + c_2 \left( \bar{I}_2 - 3 \right)  + \frac{\kappa_{\text{stab}}}{2} \left( \ln J \right)^2,\\
\Psi_{\text{anisotropic}} &= \frac{G_{\text{f}}}{2} \left(I_4  -1 \right)^2,
\end{align}
in which $c_1$ and $c_2$ are material parameters, $G_\text{f}$ is the fiber modulus, $\bar{I}_2 = \frac{1}{2}\left(\text{tr}\left(\mathbb{\bar{C}}\right)^2 - \text{tr}\left(\mathbb{\bar{C}}^2\right) \right)$, $I_4  =  \max\left( 1, \mathbb{C} : \left( \a \otimes \a  \right) \right)$, and $\a$ is the fiber direction in the reference configuration.  
The material parameters are set to $c_1 = c_2 = 9000 \, \frac{\text{dyn}}{\text{cm}^2}$.
The fiber modulus is set equal to the infinitesimal shear modulus of the Mooney--Rivlin ground material, $G_\text{f}=G=2(c_1+c_2)=36000\,\mathrm{dyn}/\mathrm{cm}^2$, and this value is also used to determine the numerical bulk modulus.

\begin{figure}[t!]
\centering
\hspace{-.1in}
	\begin{tabular}{cc}
	\includegraphics[width=\textwidth]{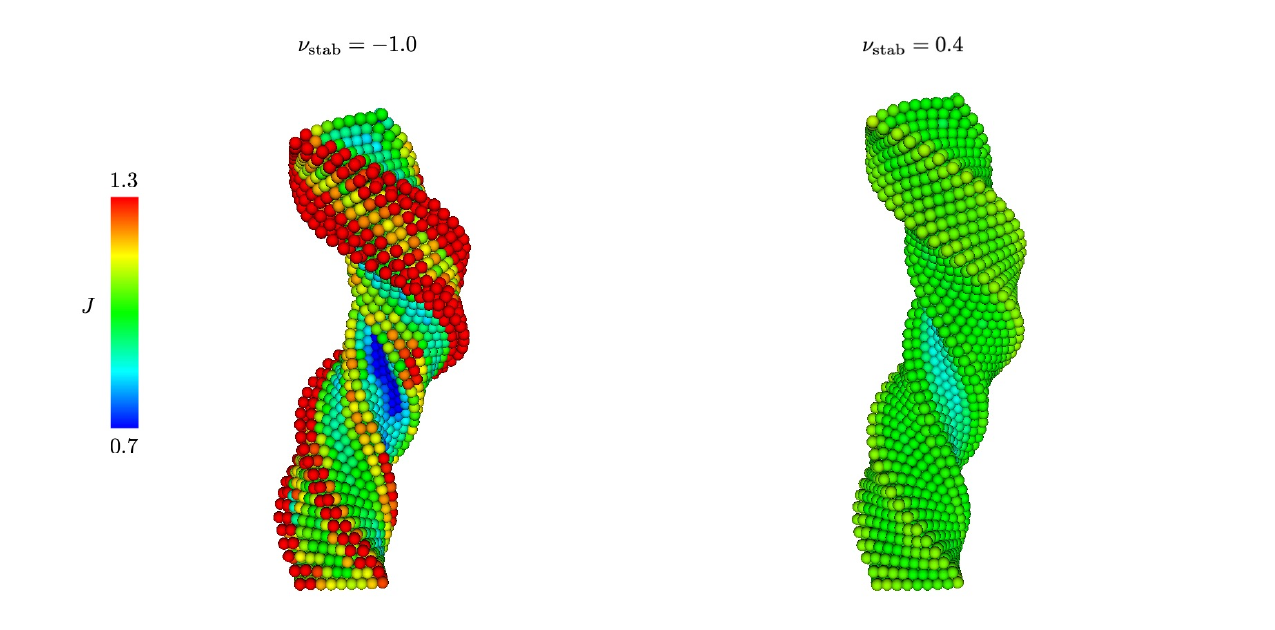}
	\end{tabular}
    \caption{Deformations of the three-dimensional anisotropic beam with the values of $J$ at material points using the Mooney--Rivlin and standard reinforced material models. The deformations are computed using $3969$ solid DoF and $\horizonsize = 2.015 \Delta X$. The left panel shows the deformation obtained using $\nu_{\text{stab}}  = -1.0$, and the right panel shows the result for $\nu_{\text{stab}}  = 0.4$.}
    \label{f:aniso_torsion_deformation}
\end{figure}

\begin{figure}[t!]
\centering
   \begin{tabular}{ccc}
        \begin{subfigure}{.33\textwidth}
          		\includegraphics[width=\textwidth]{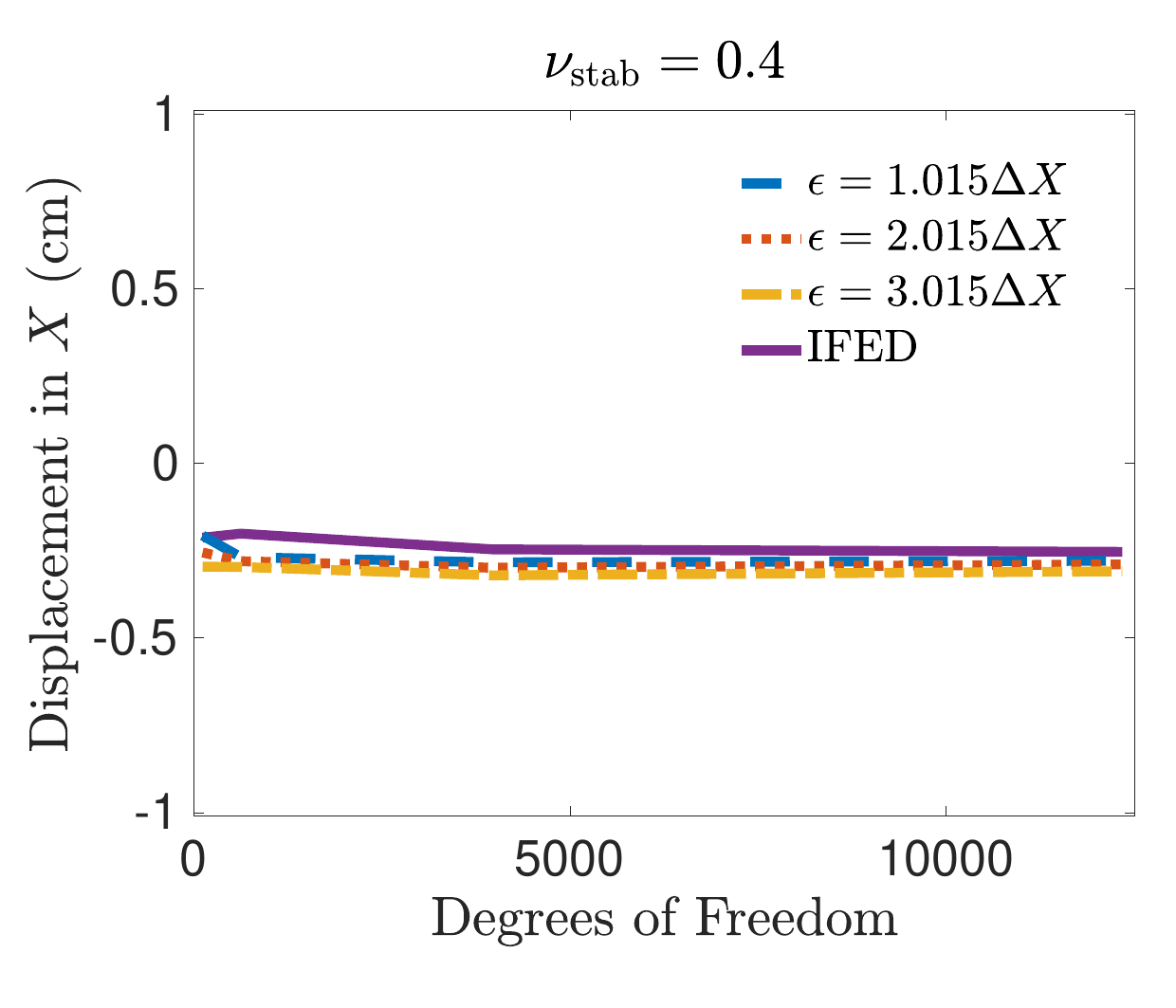}
          		 \label{f:aniso_torsion_disp_1}
        \end{subfigure}
         \begin{subfigure}{.33\textwidth}
          		\includegraphics[width=\textwidth]{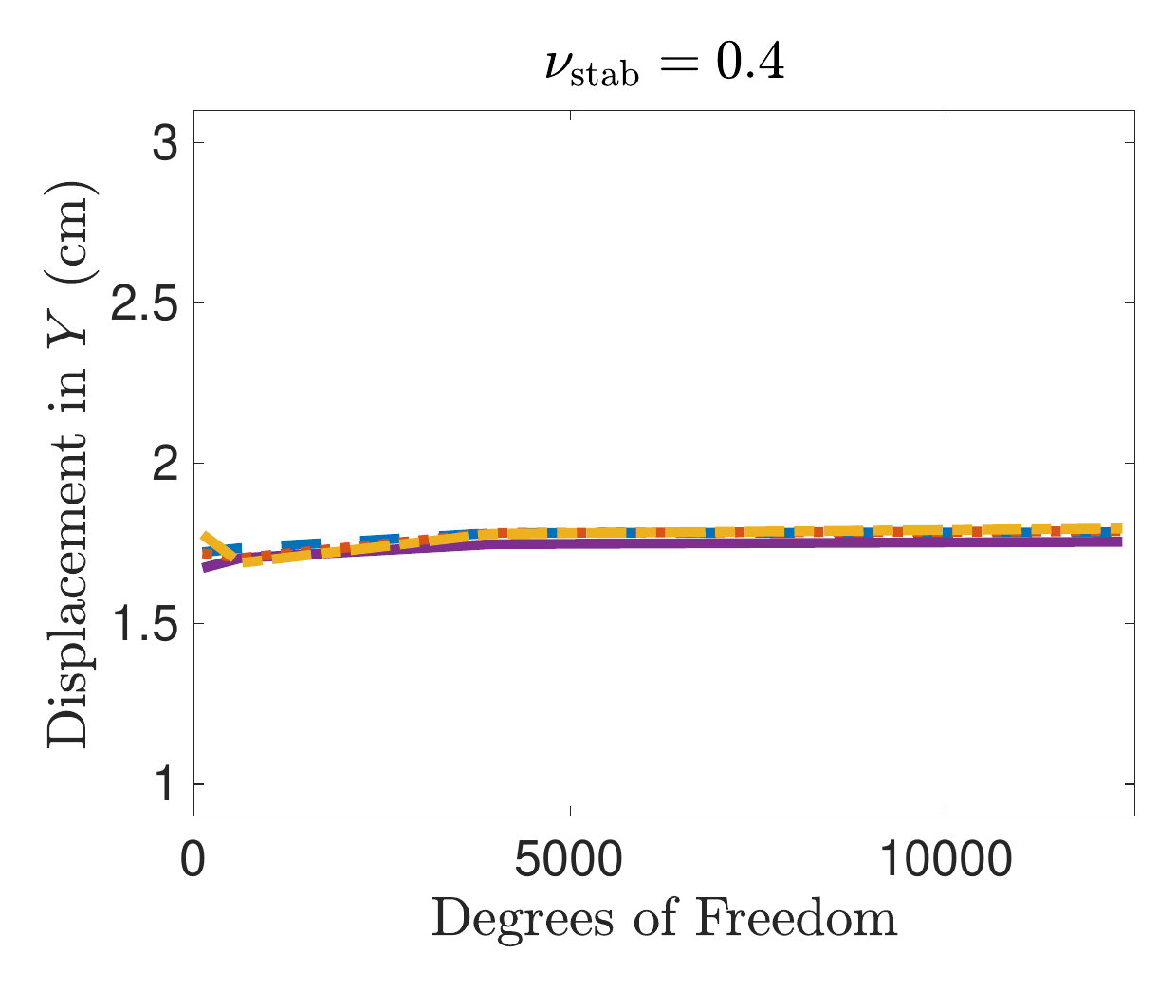}
          		 \label{f:aniso_torsion_disp_2}
        \end{subfigure} 
        \begin{subfigure}{.33\textwidth}
          		\includegraphics[width=\textwidth]{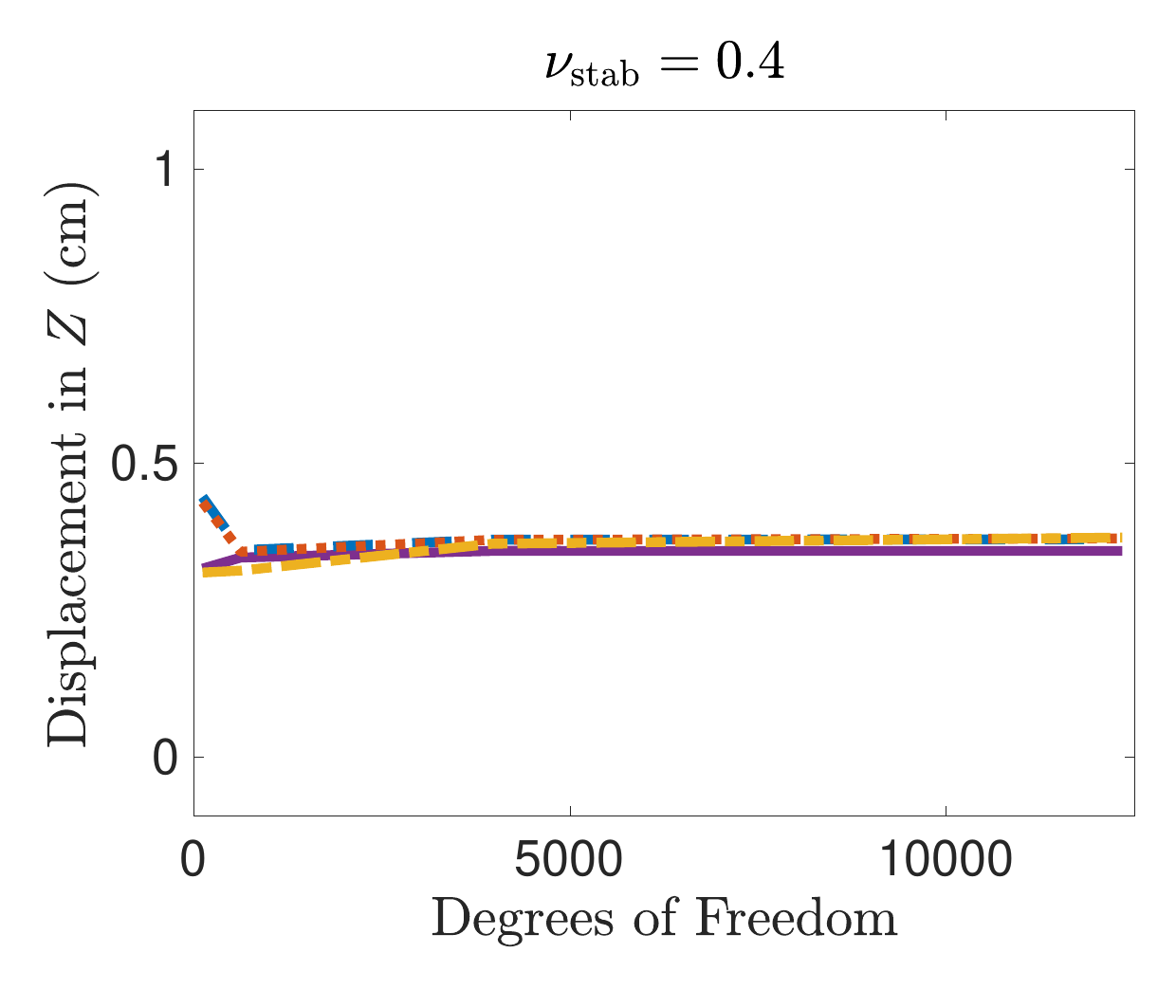}
          		 \label{f:aniso_torsion_disp_3}
        \end{subfigure}
    \end{tabular}
    \caption{Displacements in the $x$, $y$, and $z$ directions of the point of interest, highlighted in Fig.~\ref{f:aniso_torsion_schematics}, for different choices of horizon size $\horizonsize$. The number of solid DoF ranges from $117$ to $12337$.}
    \label{f:aniso_Torsion_disp}
\end{figure}

\begin{figure}[t!]
\centering
   \begin{tabular}{c}
         \begin{subfigure}{.4\textwidth}
         		\includegraphics[width=\textwidth]{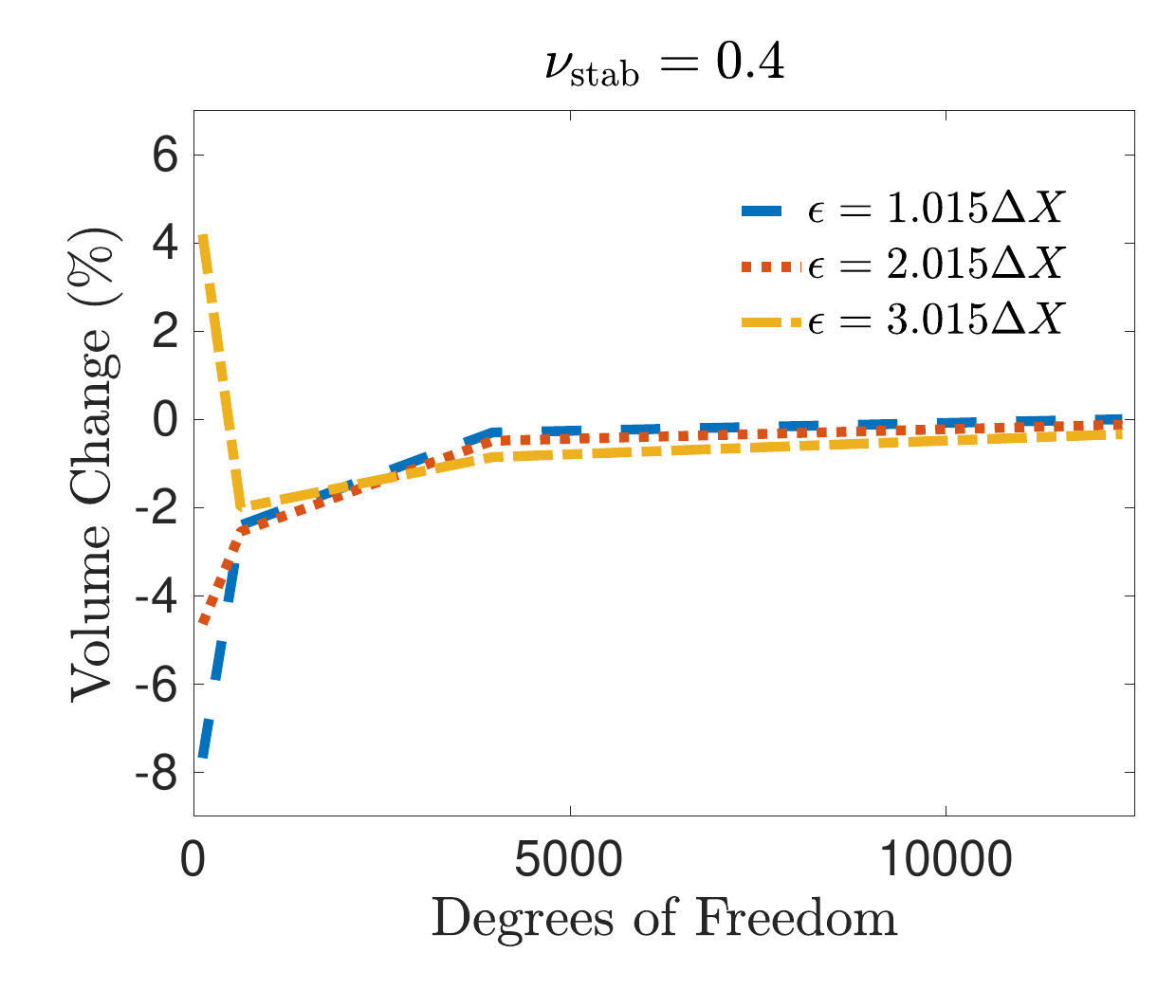}
        \end{subfigure} 
    \end{tabular}
    \caption{Volume change of the beam for different choices of horizon size $\horizonsize$. The number of solid DoF ranges from $117$ to $12337$.}
    \label{f:aniso_Torsion_vol}
\end{figure}

Fig.~\ref{f:aniso_torsion_deformation} illustrates the deformations of the beam under torsion along with the Jacobian determinant of the non-local deformation gradient tensor at each material point. 
Fig.~\ref{f:aniso_Torsion_disp} shows the maximum displacement of the point of interest located at $(3.5,0,0)$ in Fig.~\ref{f:aniso_torsion_schematics} at the steady states for various sizes of $\horizonsize$ under grid refinement. 
The displacements using IPD are comparable to the result obtained by the IFED method and converge under grid refinement to a value of approximately $-0.28 \, \text{cm}$ in the $x$ direction, $1.79 \, \text{cm}$ in the $y$ direction, and $0.37 \, \text{cm}$ in the $z$ direction.
The IFED method slowly converges compared to the IPD method.
With $28033$ solid DoF, the IFED method yields a result comparable to that obtained using the IPD method with $12337$ solid DoF.
Fig.~\ref{f:aniso_Torsion_vol} shows the volume change of the beam for different numbers of the solid DoF, ranging from $0.0035 \%$ to $7.67\%$. 
The smallest volume change obtained using the IFED method is $0.11 \%$, which is comparable to the IPD simulations.
Volume conservation can be improved by increasing the value of the numerical Poisson's ratio closer to $0.5$ in both IPD and IFED simulations.

\subsection{Adventitial strip}
\label{s:Aortic_tissue}
\begin{figure}[t!]
\centering
    \includegraphics[width=.9\textwidth]{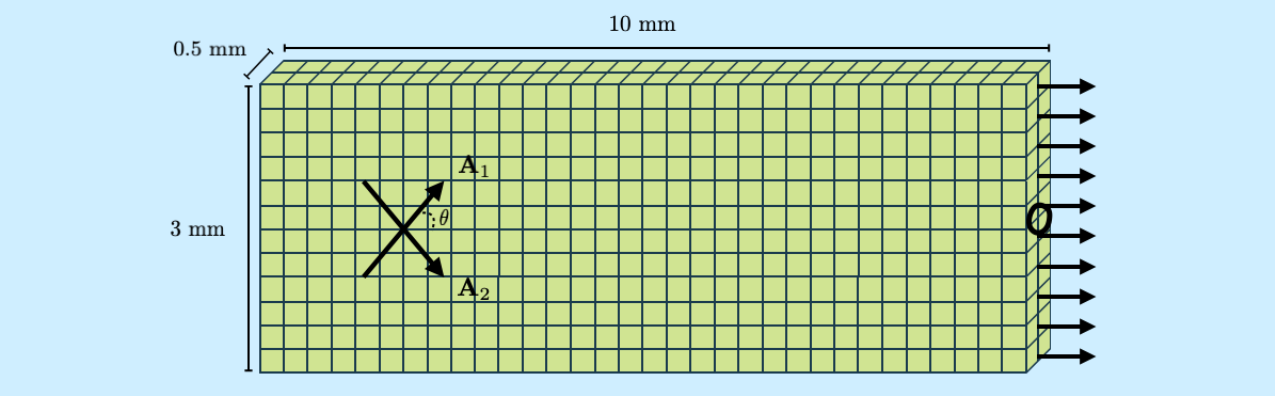}
    \caption{Schematic diagram for the three-dimensional adventitial strip. Two families of fibers are embedded in a neo-Hookean ground material. The computational domain is $\Omega = [0,L]^3$, with $L = 20 \, \text{mm}$, and the three-dimensional strip is placed at the center of the domain. Zero fluid velocity is enforced on the outer boundaries of the computational domain. The fiber directions are set to $\a_i = \left(\cos\theta, (-1)^i \sin\theta, 0 \right), \ i =1,2$ with $\theta = 49.98^\circ$.}
    \label{f:Aortic_tissue_schematics}
\end{figure}

We use a three-dimensional adventitial strip to perform uniaxial tensile tests to investigate realistic anisotropic responses of arterial walls.
This benchmark is adopted from the classical FE-based study of hyperelastic arterial models by Gasser, Ogden, and Holzapfel~\cite{Gasser_2005}.
The computational domain is the cube $\Omega = \left[0, L\right]^3$, with $L = 20 \,  \text{mm}$. 
The density and viscosity of the surrounding fluid are set to $\rho = 1.0 \,  \frac{\text{g}}{\text{cm}^3}$ and  $\mu = 0.04 \, \frac{\text{dyn$\cdot$s}}{\text{cm}^2}$, respectively. 
The left surface of the strip is fixed in place, a uniaxial tension of $1.0 \ \text{N}$ is applied to the right surface, and both tensile directional surfaces are not allowed to deform during the simulation.
We use the HGO model introduced in Sec.~\ref{s:constitutive_laws}.
The material parameters are set to $G = 7.64 \times 10^4 \  \frac{\text{dyn}}{\text{cm}^2}$, $k_1 = 966.6\times 10^4 \  \frac{\text{dyn}}{\text{cm}^2}$, and $k_2 = 524.6$, and the fiber directions are set to $\a_i = \left(\cos\theta, \left(-1\right)^i \sin\theta, 0 \right), \ i =1,2$ with $\theta = 49.98^\circ$.
We do not consider the fiber dispersion in this benchmark, i.e., $\kappa = 0$, and bond breakage is not allowed during the simulation.
For nearly incompressible elasticity of the tissue strip, we only test $\nu_{\text{stab}} = 0.4$.
Fig.~\ref{f:Aortic_tissue_schematics} provides a schematic of this test case. 

 \begin{figure}[t!]
\centering
	\includegraphics[width=0.27\textwidth]{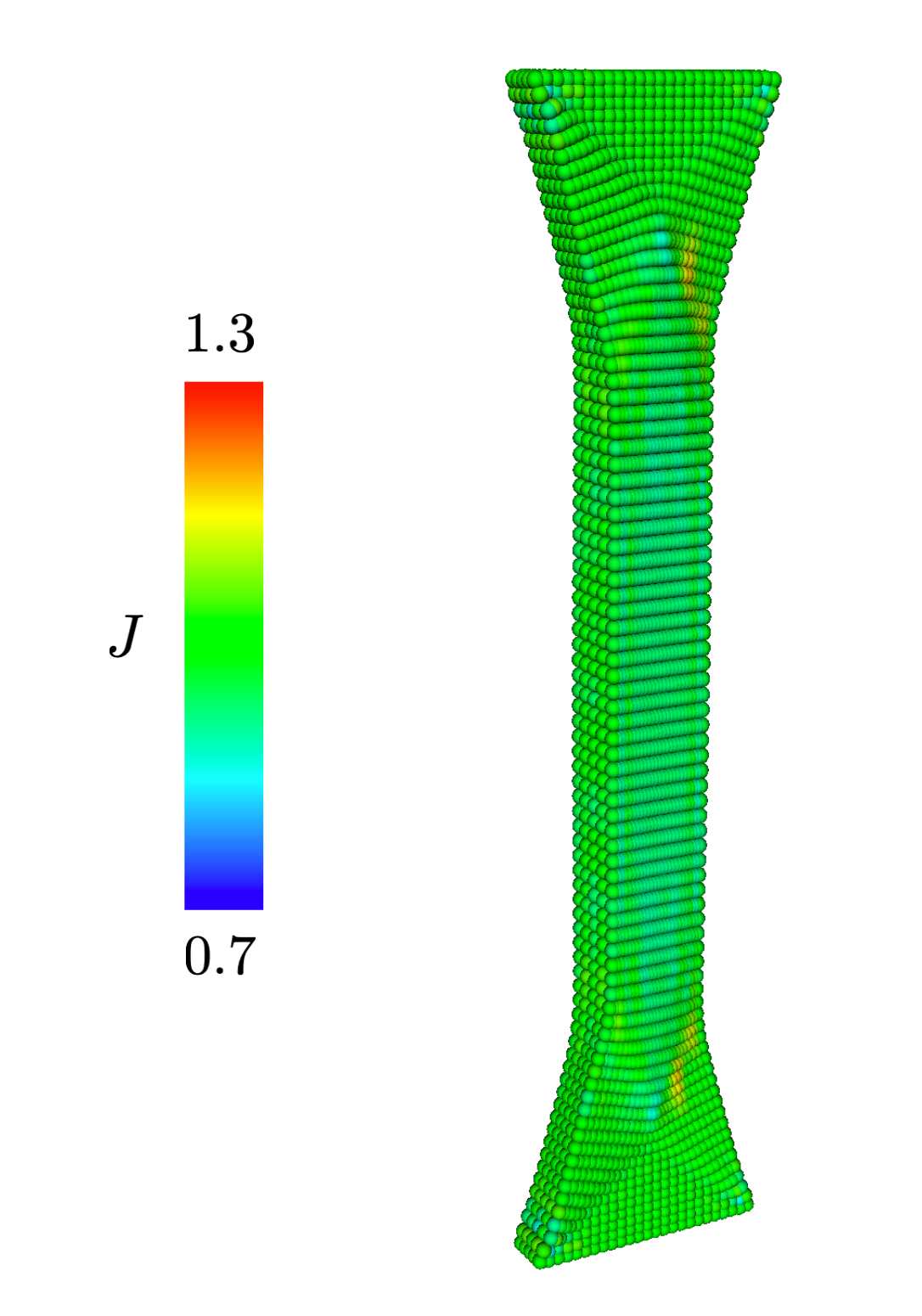} 
    \caption{A deformation of a three-dimensional tissue strip under a uniaxial tensile load of $1.0 \, \text{N}$ with the values of $J$ at material points using the HGO material model with $G = 7.64 \times 10^4 \  \frac{\text{dyn}}{\text{cm}^2}$, $k_1 = 966.6\times 10^4 \  \frac{\text{dyn}}{\text{cm}^2}$, and $k_2 = 524.6$, and the fiber directions are set to $\a_i = \left(\cos\theta, \left(-1\right)^i \sin\theta, 0 \right), \ i =1,2$ with $\theta = 49.98^\circ$. The simulation uses $10125$ solid DoF, $\horizonsize = 2.015 \Delta X$, and $\nu_\text{stab} = 0.4$.  }
    \label{f:Aortic_tissue_deformation}
\end{figure}

\begin{figure}[t!]
\centering
    \includegraphics[width=.45\textwidth]{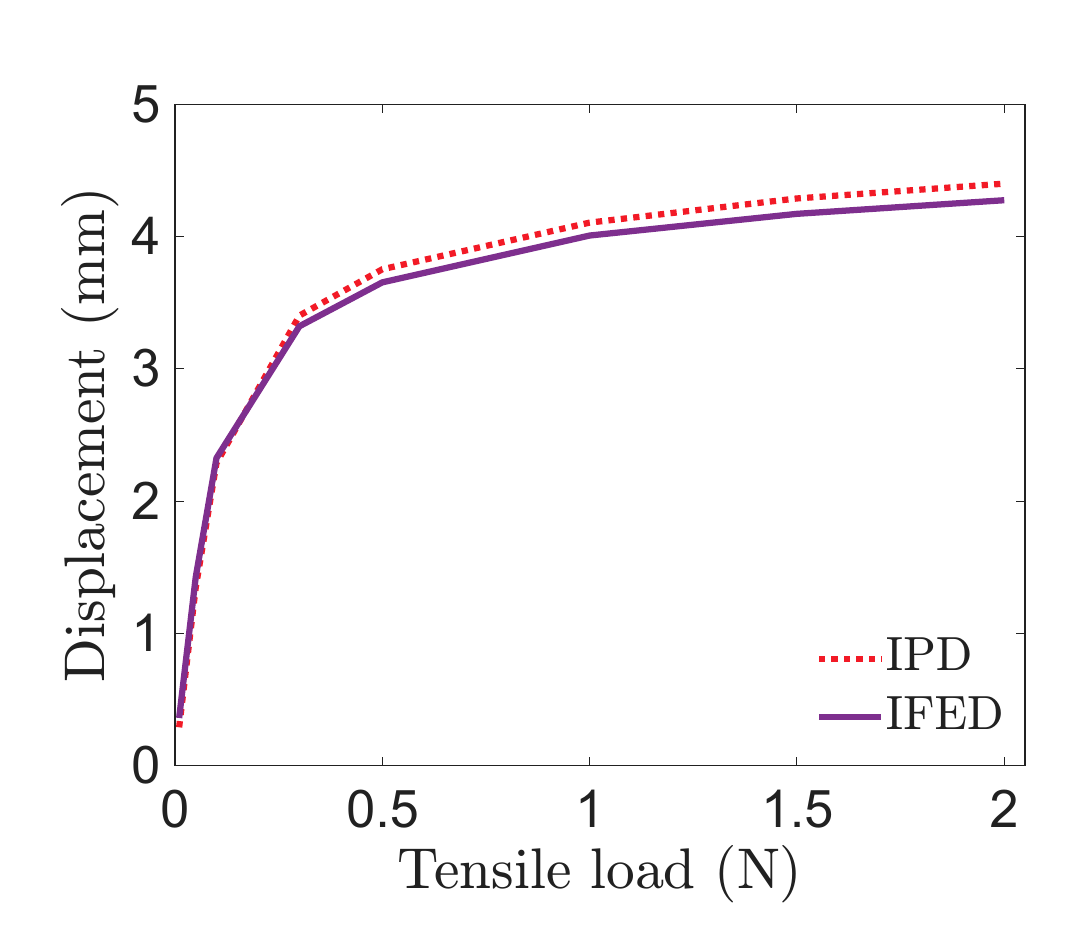}
    \caption{Displacements of the center point of the right surface of adventitial tissue, highlighted in Fig.~\ref{f:Aortic_tissue_schematics}, obtained using the IPD and IFED methods. The peridynamic horizon size is set to $\horizonsize = 2.015 \Delta X$ and the numerical Poisson's ratio is set to $\nu_\text{stab} = 0.4$ for both IPD and IFED simulations. The solid DoF are $10125$. }
    \label{f:Aortic_tissue_stretch}
\end{figure}

Fig.~\ref{f:Aortic_tissue_deformation} illustrates the deformations of the strip under the uniaxial tensile load of $1.0 \, \text{N}$ along with the Jacobian determinant of the non-local deformation gradient tensor at each material point. 
We further quantify the response over a range of applied loads, and Fig.~\ref{f:Aortic_tissue_stretch} shows the maximum displacement of the center point of the right surface in Fig.~\ref{f:Aortic_tissue_schematics} at the steady states for different external loads. 
The maximum difference in the axial displacement between the IPD and IFED results is $0.13 \, \text{mm}$, which is expected to decrease with further grid refinement.

\subsection{Anisotropic failure}
\label{s:aniso_failure}
\begin{figure}[t!]
\centering
    \includegraphics[width=.45\textwidth]{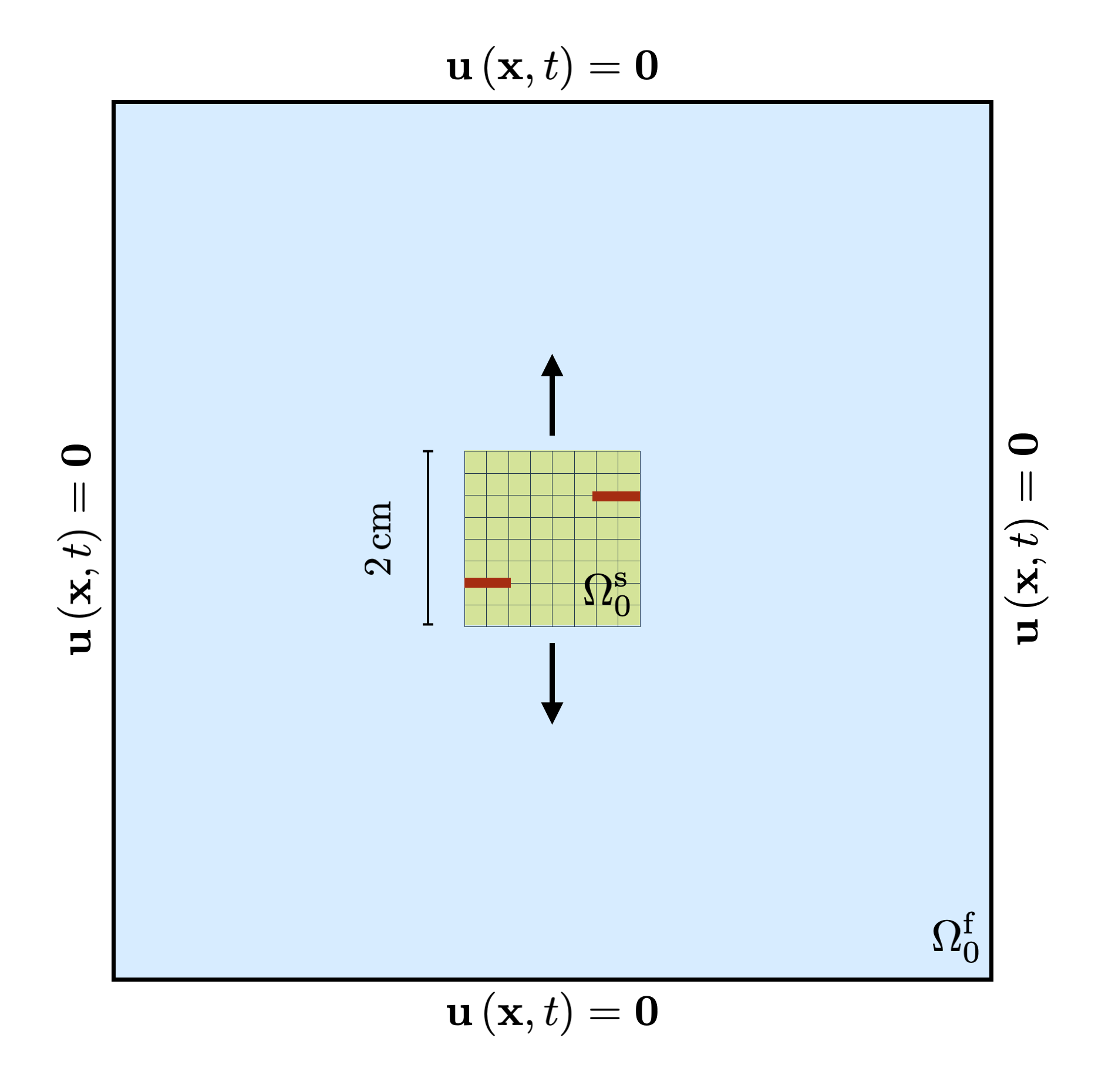}
    \caption{Schematic diagram for the two-dimensional anisotropic hyperelastic sheet with pre-existing notches. One family of fibers is embedded in a neo-Hookean ground material. The computational domain is $\Omega = [0,L]^2$, with $L = 10 \, \text{cm}$, and the two-dimensional sheet is placed at the center of the domain. Zero fluid velocity is enforced on the outer boundaries of the computational domain. We test three different fiber directions, $(\cos\theta_i,\sin\theta_i)$, for $i=1,2,3$ with $\theta_1 = 0^\circ$, $\theta_2 = 22.5^\circ$, and $\theta_3 = 45^\circ$.}
    \label{f:aniso_failure_schematics}
\end{figure}
We use a two-dimensional anisotropic hyperelastic sheet with two pre-existing notches to investigate anisotropic material damage and failure under uniaxial tensile loads.
The computational domain is $\Omega = \left[0, L\right]^2$, with $L = 10 \,  \text{cm}$, and the notches with the length of $0.5 \, \text{cm}$ are located on the top-right and bottom-left, respectively.
We use the six-point IB kernel function and set $M_{\text{FAC}} = 1$.
Displacement boundary conditions with a fixed velocity of $1 \, \text{cm}/\text{s}$ are applied to top and bottom surfaces in upward and downward directions, respectively, until each side moves $2 \, \text{cm}$.
Both tensile directional surfaces are not allowed to deform during the simulation. 
We use the HGO model with $G = 7.64 \times 10^2 \  \frac{\text{dyn}}{\text{cm}^2}$,  $k_1 = 966.6\times 10^2 \  \frac{\text{dyn}}{\text{cm}^2}$, and $k_2 = 524.6$, which is softer than the adventitial tissue used in Sec.~\ref{s:Aortic_tissue} to facilitate large deformations for the failure tests.
Note that the fiber dispersion is not considered in this benchmark.
The final time is $T_{\text{f}} = 5 \, \text{s}$.
Fig.~\ref{f:aniso_failure_schematics} provides a schematic of this failure test.

To investigate damage and crack growth depending on the fiber directions, we perform three failure tests using different fiber orientations, $\theta_1 = 0^\circ$, $\theta_2 = 22.5^\circ$  and $\theta_3 = 45^\circ$. 
We use the numerical Poisson's ratio of $\nu_\text{stab} = 0.4$ and the horizon size of $\horizonsize = 3.015 \Delta X$.
For the ductile failure model, the critical stretch values are set to $\sca = 1.5$ and $\scb = 1.8$.
Additionally, to localize the failure process in the central region of the tissue block, bond damage and failure are not allowed for the bonds that lie entirely above or below the notches.

\begin{figure}[t!]
\centering
	\begin{tabular}{ccc}
        \begin{subfigure}{.15\textwidth}
          		\includegraphics[width=\textwidth]{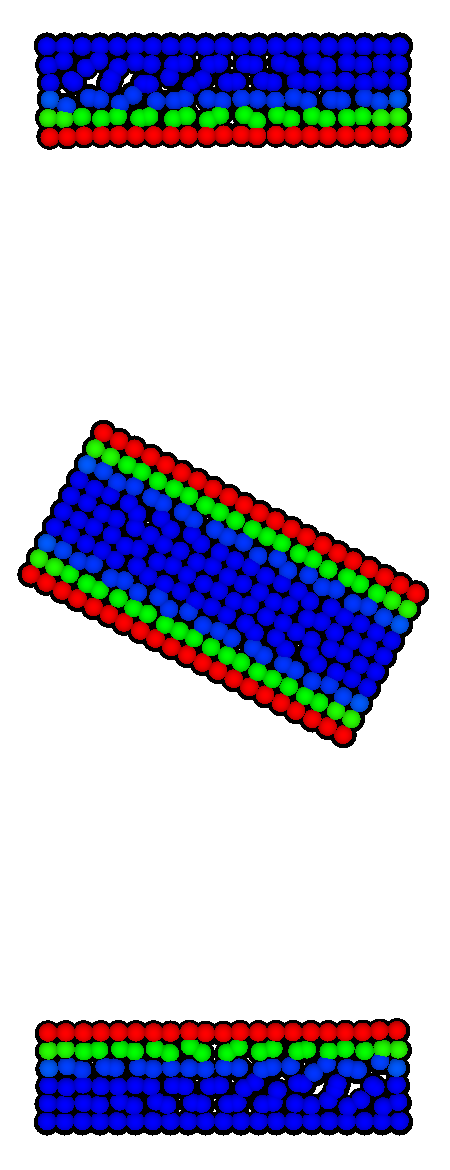}
          		 \caption{}
          		 \label{f:aniso_failure_0}
        \end{subfigure}
        \hspace{.07\textwidth}
        \begin{subfigure}{.15\textwidth}
          		\includegraphics[width=\textwidth]{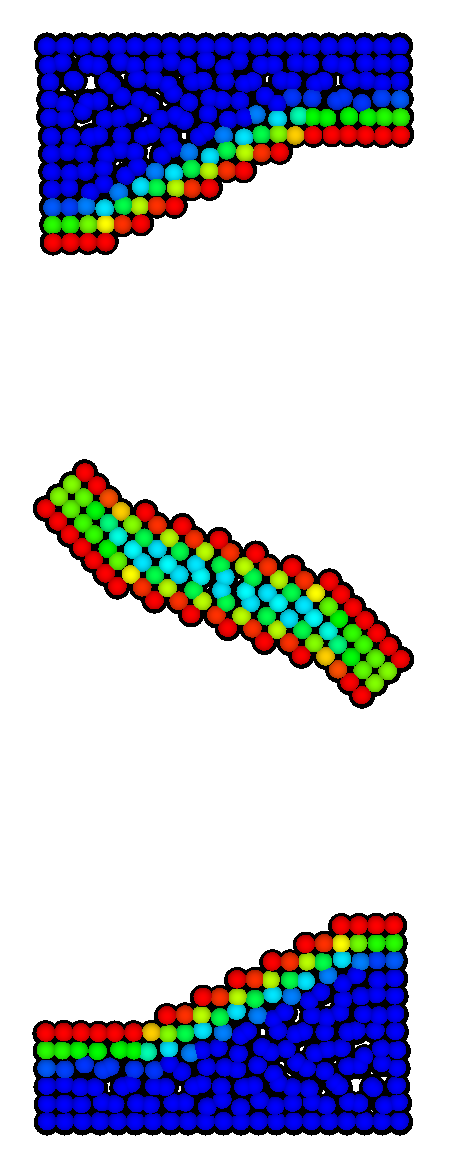}
          		 \caption{}
          		 \label{f:aniso_failure_22}
        \end{subfigure}
        \hspace{.07\textwidth}
         \begin{subfigure}{.15\textwidth}
          		\includegraphics[width=\textwidth]{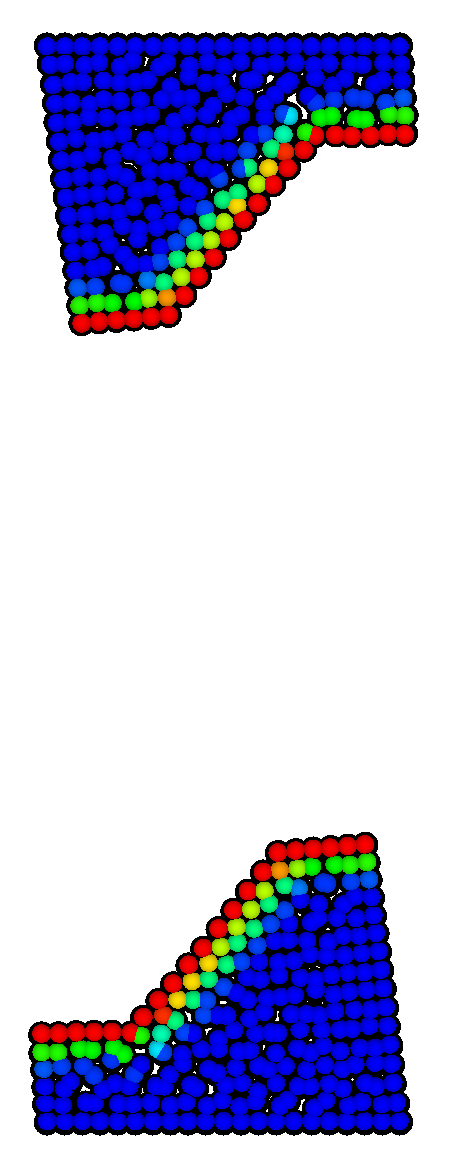}
          		 \caption{}
          		 \label{f:aniso_failure_45}
        \end{subfigure}
        \hspace{.02\textwidth}
        \begin{subfigure}{.09\textwidth}
        	\centering
          		\includegraphics[width=\textwidth]{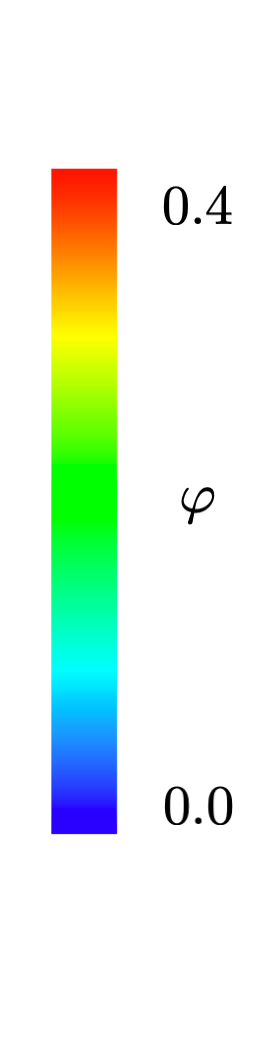}
          		\vspace{.15in}
        \end{subfigure} 
    \end{tabular}
\caption{Three different failures of the two-dimensional sheet with the local damage in NOSB-PD. Note that $\varphi = 1$ indicates all initially connected bonds are broken. The deformations are computed using $462$ solid DoF, $\horizonsize = 3.015 \Delta X$, and $\nu_{\text{stab}} = 0.4$. The fiber direction is set to (a) $\theta_1 = 0^\circ$, (b) $\theta_2 = 22.5^\circ$, and (c) $\theta_3 = 45^\circ$.}
\label{f:aniso_failure}
\end{figure}

\begin{figure}[t!]
\centering
\includegraphics[width=.8\textwidth]{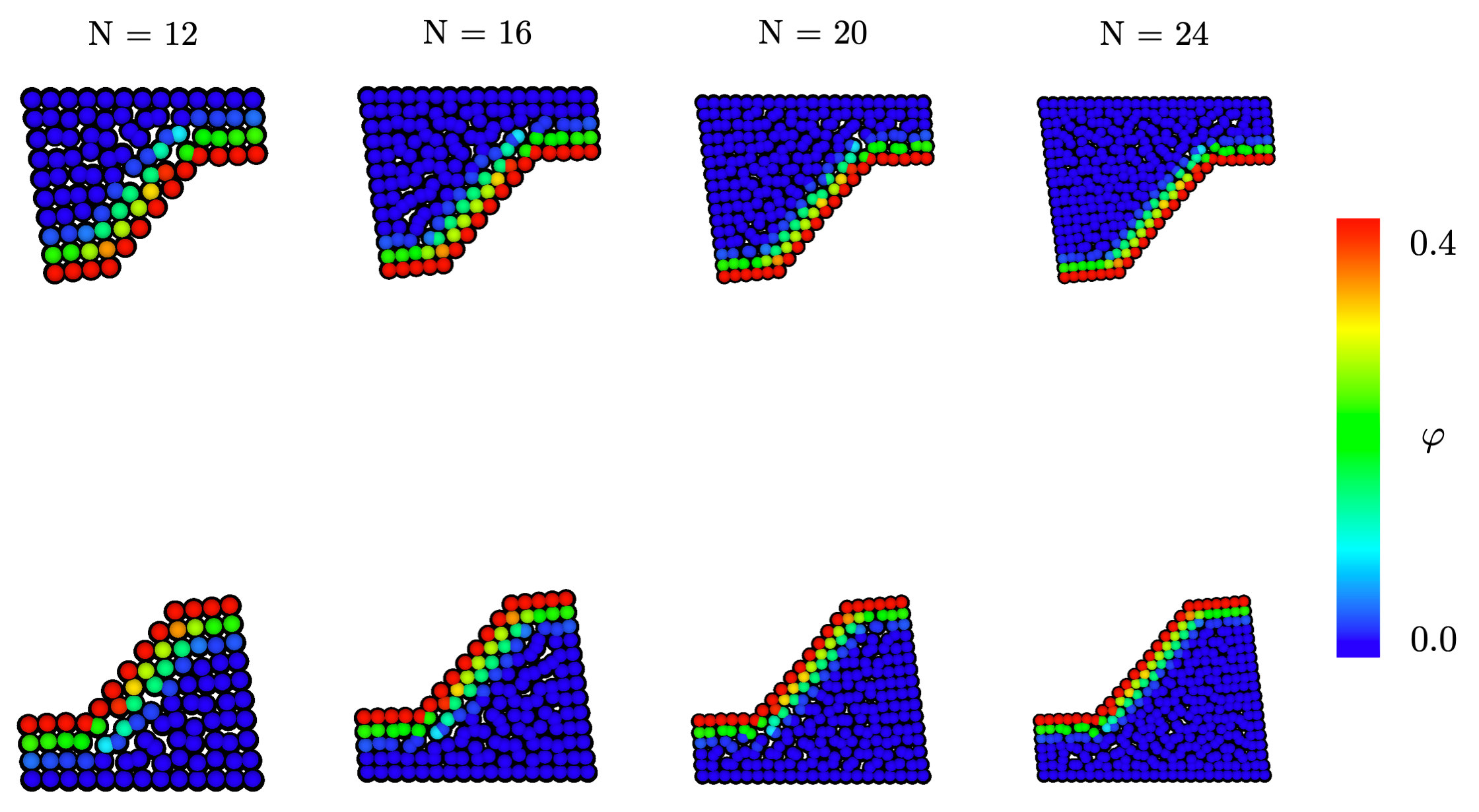}
\caption{Grid-refinement comparison of the crack trajectory in the two-dimensional sheet shown in the deformed configuration with the local damage. The fiber direction is $\theta_3 = 45^\circ$. $N = 12$ corresponds to 182 solid DoF, $N = 16$ corresponds to 306 solid DoF, $N = 20$ corresponds to 462 solid DoF, and $N = 24$ corresponds to 650 solid DoF.}
\label{f:aniso_crack_convergence}
\end{figure}

\begin{figure}[t!]
\centering
	\begin{tabular}{ccc}
        \begin{subfigure}{.33\textwidth}
          		\includegraphics[width=\textwidth]{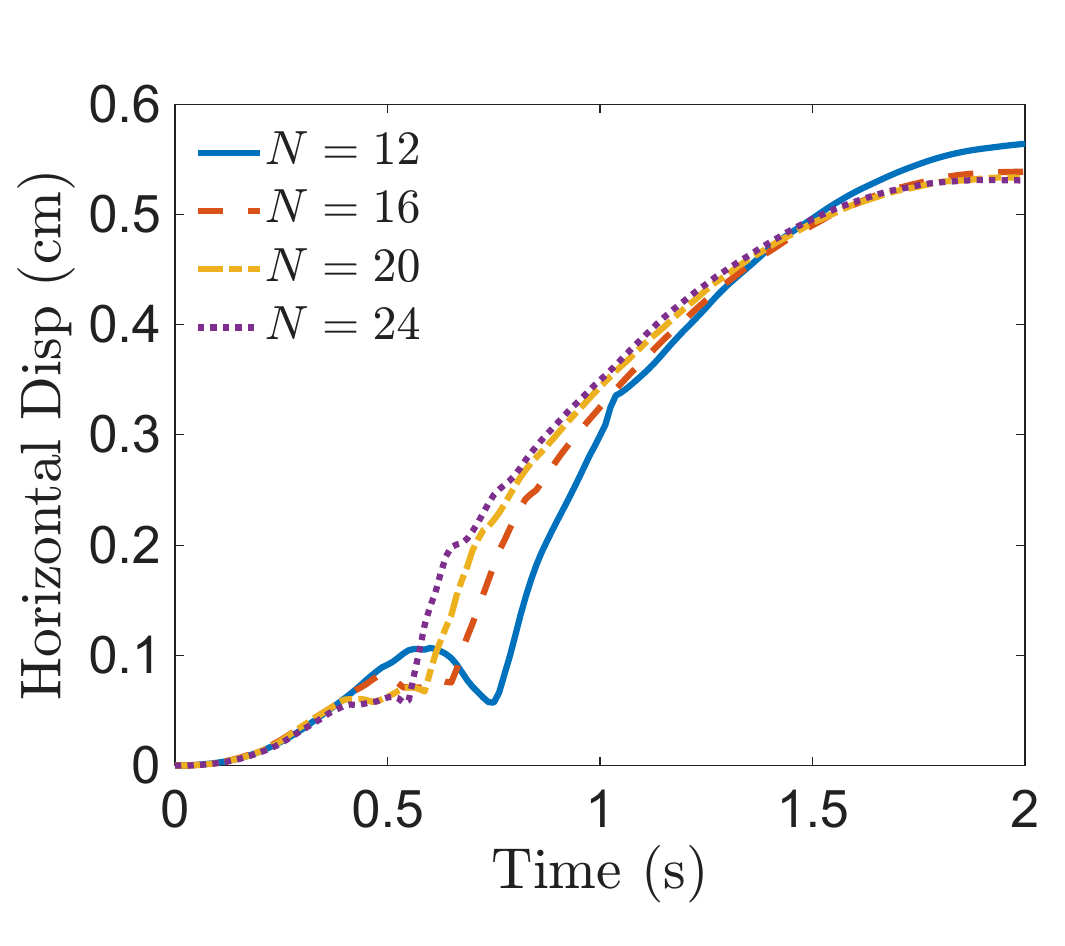}
          		 \caption{}
          		 \label{f:aniso_failure_disp1}
        \end{subfigure}
        \begin{subfigure}{.33\textwidth}
          		\includegraphics[width=\textwidth]{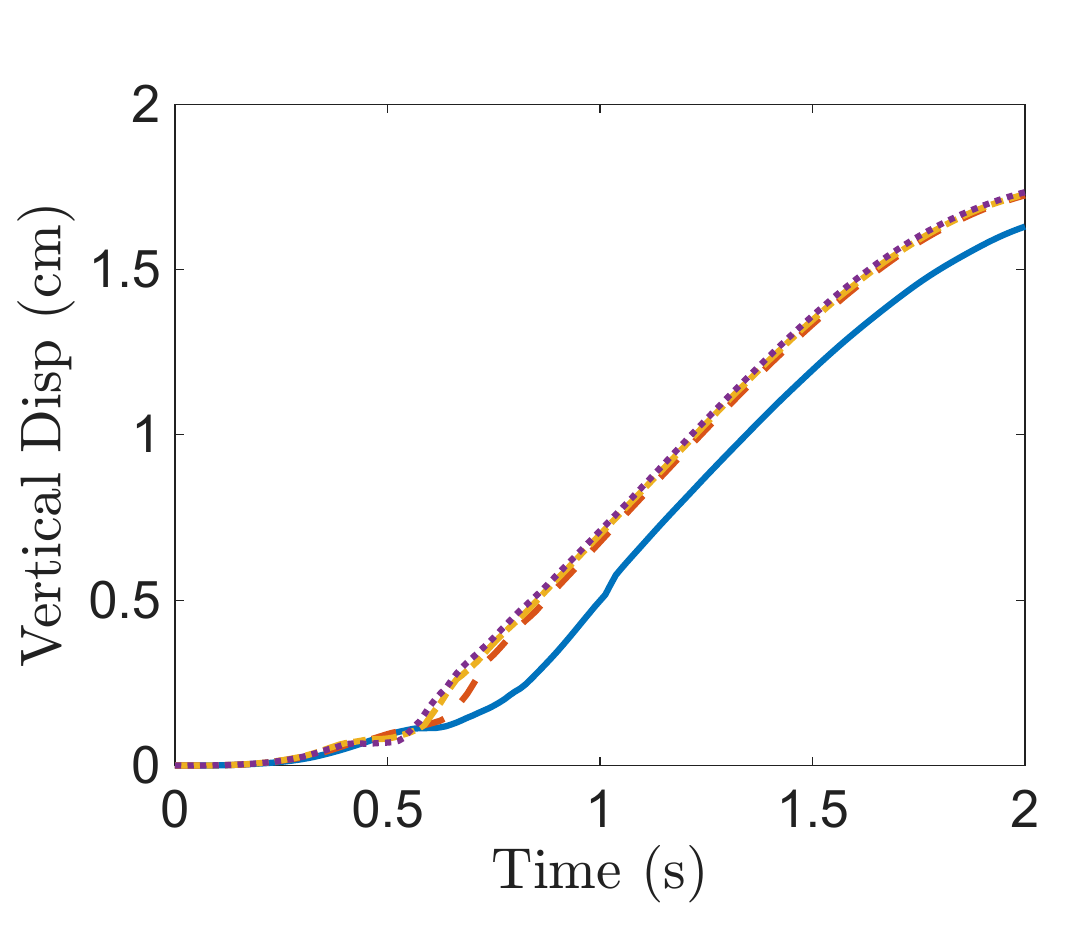}
          		 \caption{}
          		 \label{f:aniso_failure_disp2}
        \end{subfigure} 
         \begin{subfigure}{.33\textwidth}
          		\includegraphics[width=\textwidth]{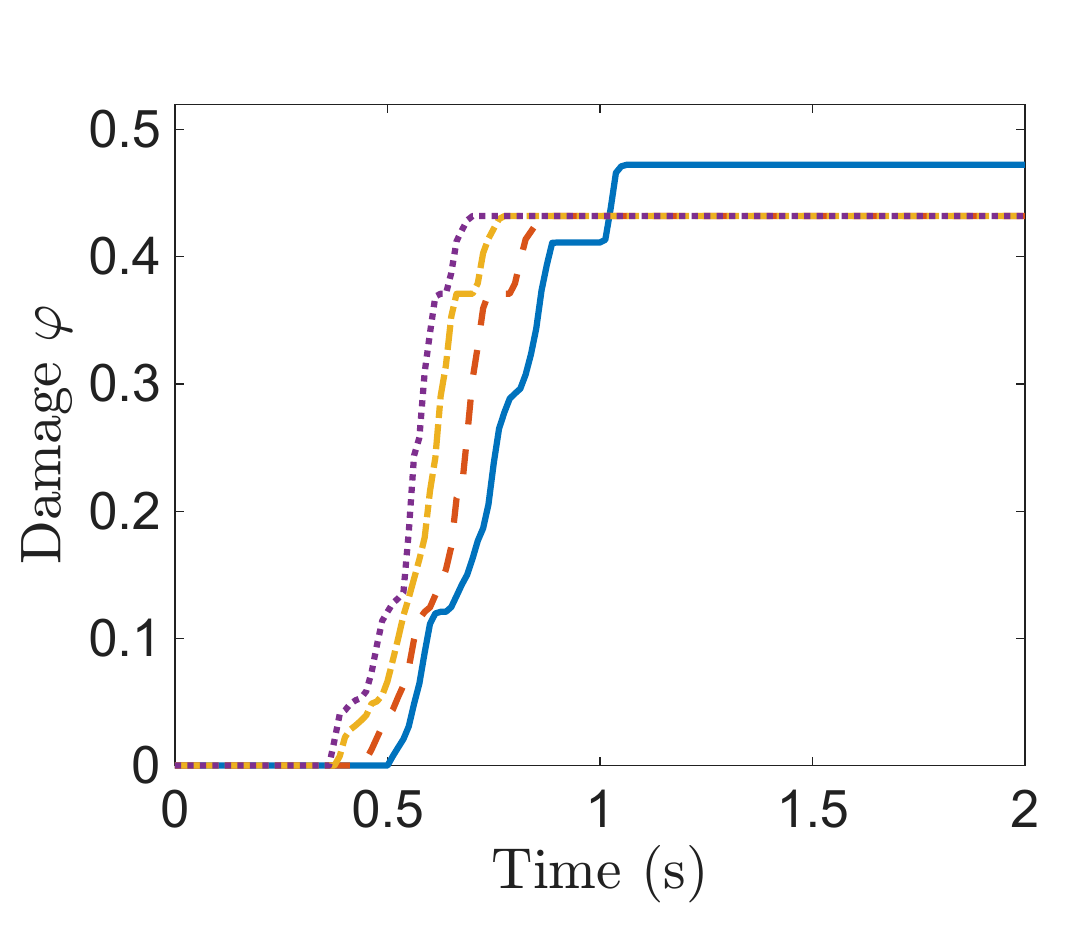}
          		 \caption{}
          		 \label{f:aniso_failure_damage}
        \end{subfigure} 
    \end{tabular}
\caption{(a) Horizontal displacements, (b) vertical displacements, and (c) local damage growth at the center of the sheet during the failure process under grid refinement with $\theta_3 = 45^\circ$. $N = 12$ corresponds to 182 solid DoF, $N = 16$ corresponds to 306 solid DoF, $N = 20$ corresponds to 462 solid DoF, and $N = 24$ corresponds to 650 solid DoF.}
\label{f:aniso_failure_convergence}
\end{figure}

Fig.~\ref{f:aniso_failure} demonstrates the effect of an embedded fiber direction on the failure process along with the local damage $\varphi$. 
The crack propagation path is well aligned with the fiber direction in the immersed structure for all three cases.
Fig.~\ref{f:aniso_crack_convergence} shows that the predicted crack trajectory remains consistent under grid refinement and aligns with the prescribed fiber direction at all resolutions.
We further demonstrate grid convergence of the failure process by tracking the horizontal and vertical displacements and damage growth at the center particle of the two-dimensional sheet with the fiber direction $\theta = 45^\circ$; see Fig.~\ref{f:aniso_failure_convergence}. 
Note that the observed grid consistency is specific to the selected $M_{\mathrm{FAC}}$, IB kernel, and PD horizon, and the sensitivity of the results to other choices is not examined in this study.

\subsection{Anisotropic elastic band}
\label{s:ansio_els_bnd}
\begin{figure}[t!]
\centering
    \includegraphics[width=.8\textwidth]{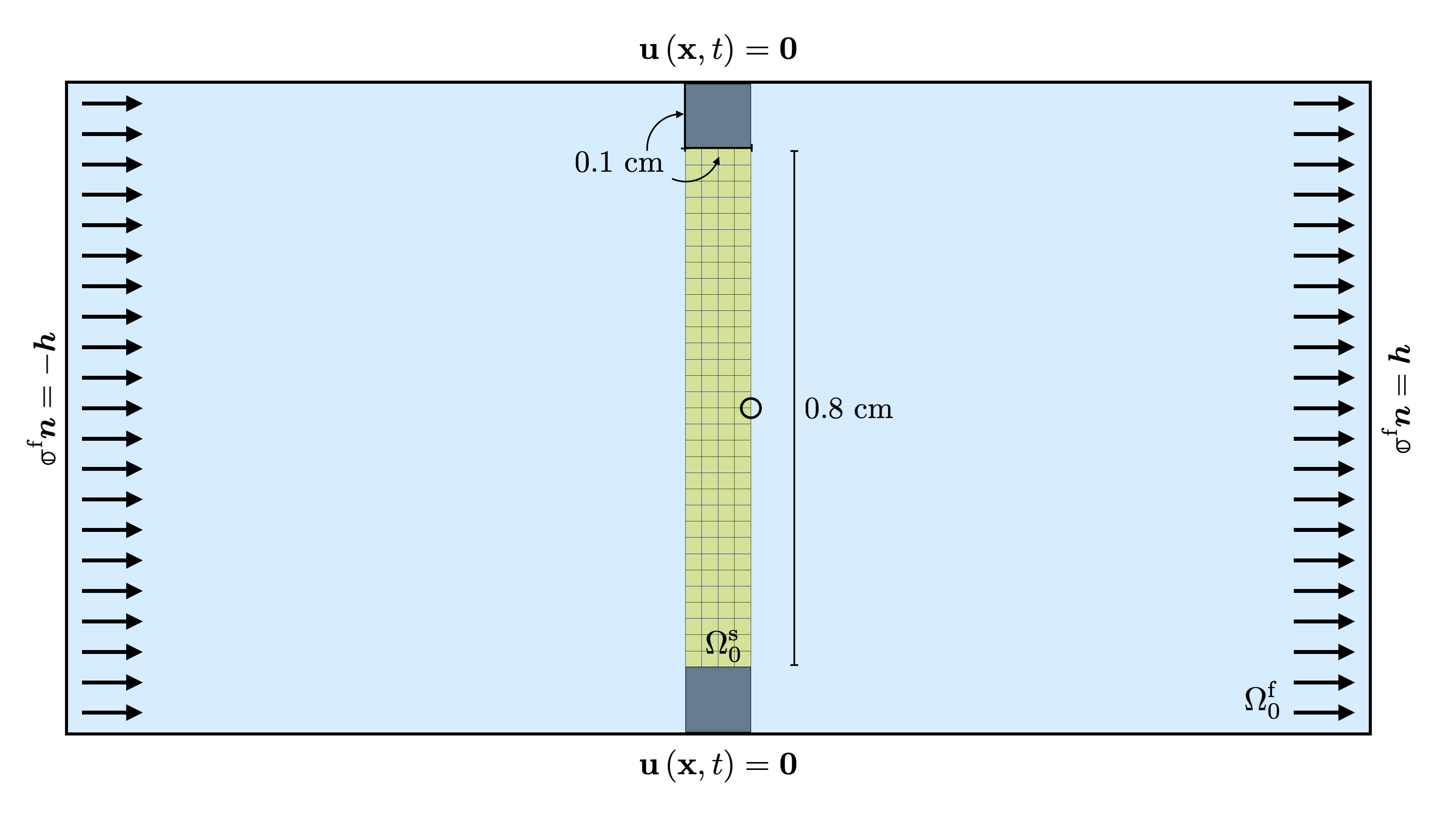}
    \caption{Schematic diagram for the anisotropic elastic band benchmark (Sec.~\ref{s:ansio_els_bnd}).  The initial configurations of the immersed structure and a fluid are denoted by $\Omega_0^{\text{s}}$ and $\Omega_0^{\text{f}}$, respectively. The entire computational domain is $\Omega = \Omega_0^{\text{s}} \cup \Omega_0^{\text{f}}$. Zero fluid velocity is enforced on the top and bottom boundaries of the computational domain, and fluid traction boundary conditions are applied to the left and right boundaries. Fluid traction is set to $\bm{h} \left(t\right) = \left(30 , 0\right)\, \frac{\mathrm{dyn}}{\mathrm{cm}^2}$.}
    \label{f:ELS_BND_schematics}
\end{figure}

To investigate purely fluid-driven deformations and subsequent material damage and failure processes in anisotropic materials, we use the two-dimensional elastic band test~\cite{KIM2023112466}. 
The computational domain is $\Omega = \left[0, 2L\right] \times \left[0, L\right]$, with $L = 1 \, \mathrm{cm}$. 
Fluid traction boundary conditions are imposed on the boundaries of the computational domain as $\vec{\bbsigma}^{\text{f}} \left(\x,t\right) \bm{n} \left(\x\right) = \bm{h}\left(t\right)$ and $\vec{\bbsigma}^{\text{f}} \left(\x,t\right) \bm{n} \left(\x\right) = - \bm{h}\left(t\right)$ on the left and right, respectively, in which $\vec{\bbsigma}^{\text{f}} $ is the fluid stress tensor and $\bm{h} \left(t\right) = \left(30 , 0\right) \, \frac{\mathrm{dyn}}{\mathrm{cm}^2}$.
The simulation time is $T_{\text{f}} = 0.3 \, \mathrm{s}$. 
Zero fluid velocity conditions are applied to the top and bottom boundaries of the computational domain. 
Both top and bottom surfaces of the elastic band are attached to fixed blocks, which prevent the fluid flow between the wall and the flexible band. 
We use a composite B-spline delta function~\cite{GRUNINGER2026114472}.
Fig.~\ref{f:ELS_BND_schematics} provides a schematic of this benchmark. 

To investigate the effect of a fiber against the fluid flow, we use two different fiber orientations, $\theta_1 = 90^\circ$ (perpendicular) and $\theta_2 = 0^\circ$ (parallel) with the modified standard reinforced model as in Sec.~\ref{s:aniso_Benchmark_Compression}.
We use the shear modulus of $G = G_{\text{f}} = 200\, \frac{\mathrm{dyn}}{\mathrm{cm}^2}$ and the horizon size of $\horizonsize = 3.015 \Delta X$.
For ductile failure, the critical stretch values  are set to $\sca = 4.5$ and $\scb = 5.4$.

\begin{figure}[t!]
\centering
   \begin{tabular}{c}
        \begin{subfigure}{.7\textwidth}
          		\includegraphics[width=\textwidth]{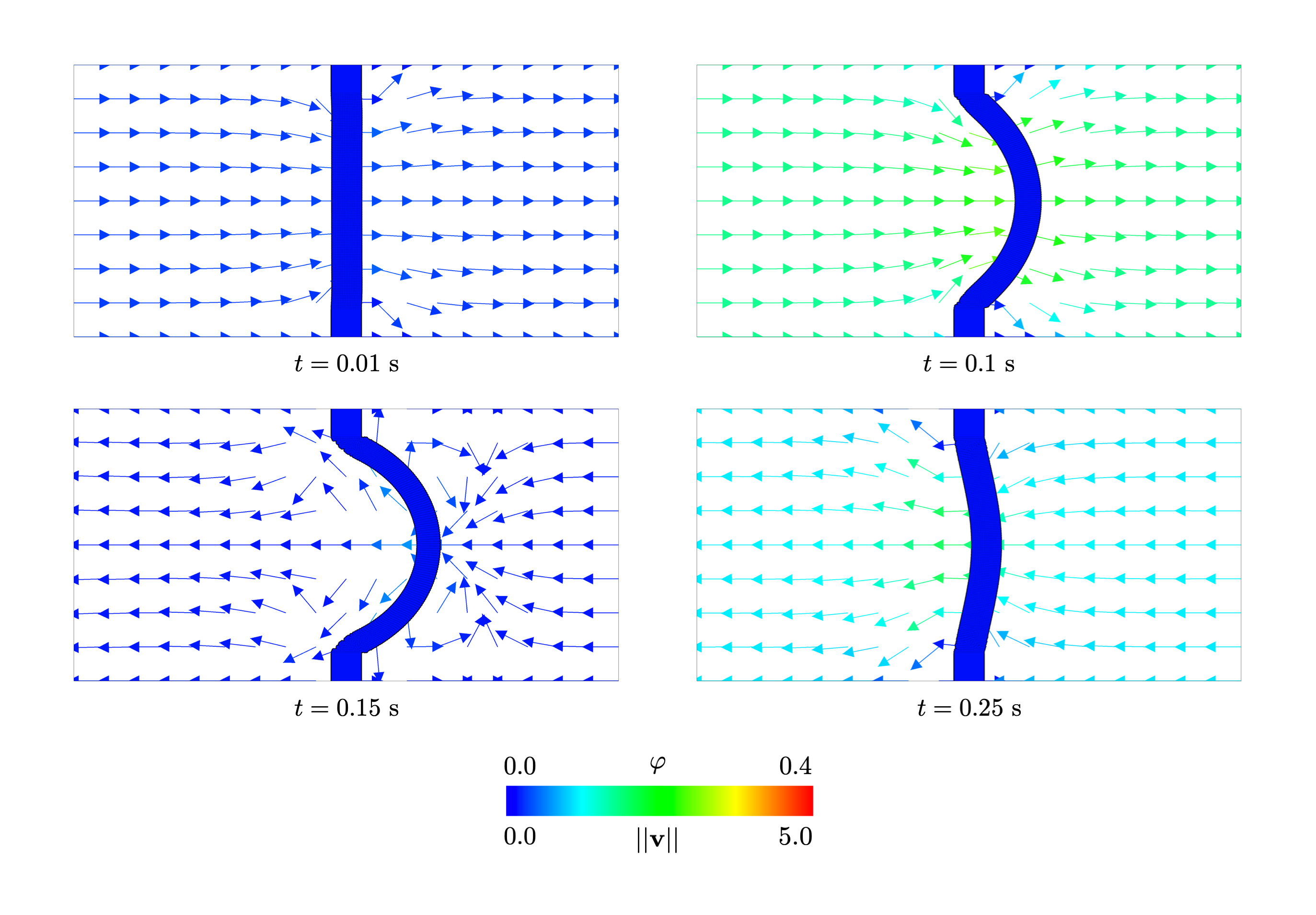}
          		\caption{}
          		 \label{f:ELS_BND_90}
        \end{subfigure}\\
         \begin{subfigure}{.7\textwidth}
        			\includegraphics[width=\textwidth]{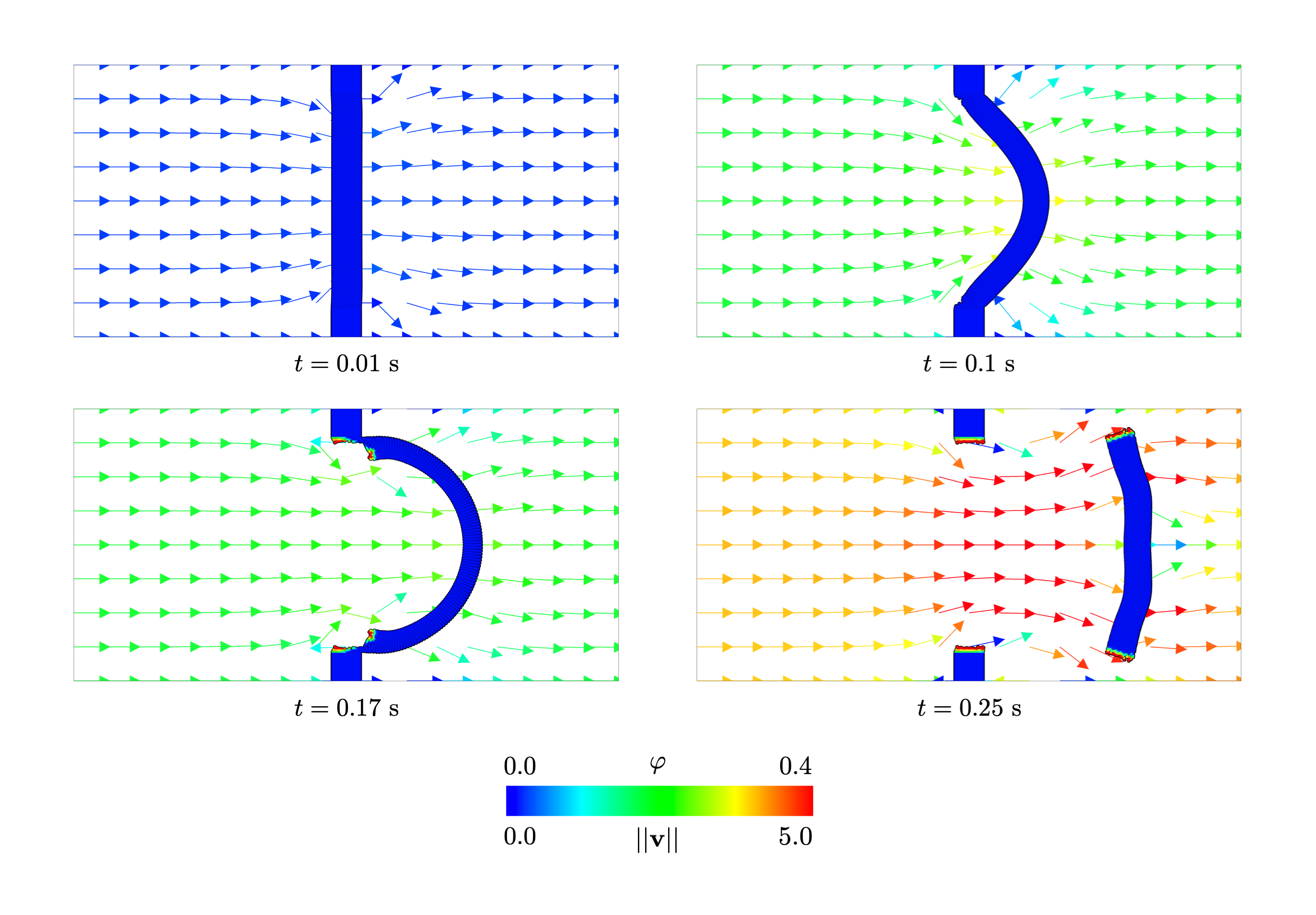}
          		\caption{}
          		 \label{f:ELS_BND_0}
        \end{subfigure} 
    \end{tabular}
    \caption{Deformations of the elastic band with two different fiber directions (a) $\theta_1 = 90^\circ$ and (b) $\theta_2 = 0^\circ$. Note that $\varphi = 0$ implies that all initial bonds are connected, and $\varphi = 1$ implies that all initial bonds are disconnected. The deformations are computed using $2193$ solid DoF, $\horizonsize = 3.015 \Delta X$, and $\nu_\text{stab} = 0.4$.}
    \label{f:ELS_BND_deformations}
\end{figure}

\begin{figure}[t!]
\centering
\includegraphics[width=.4\textwidth]{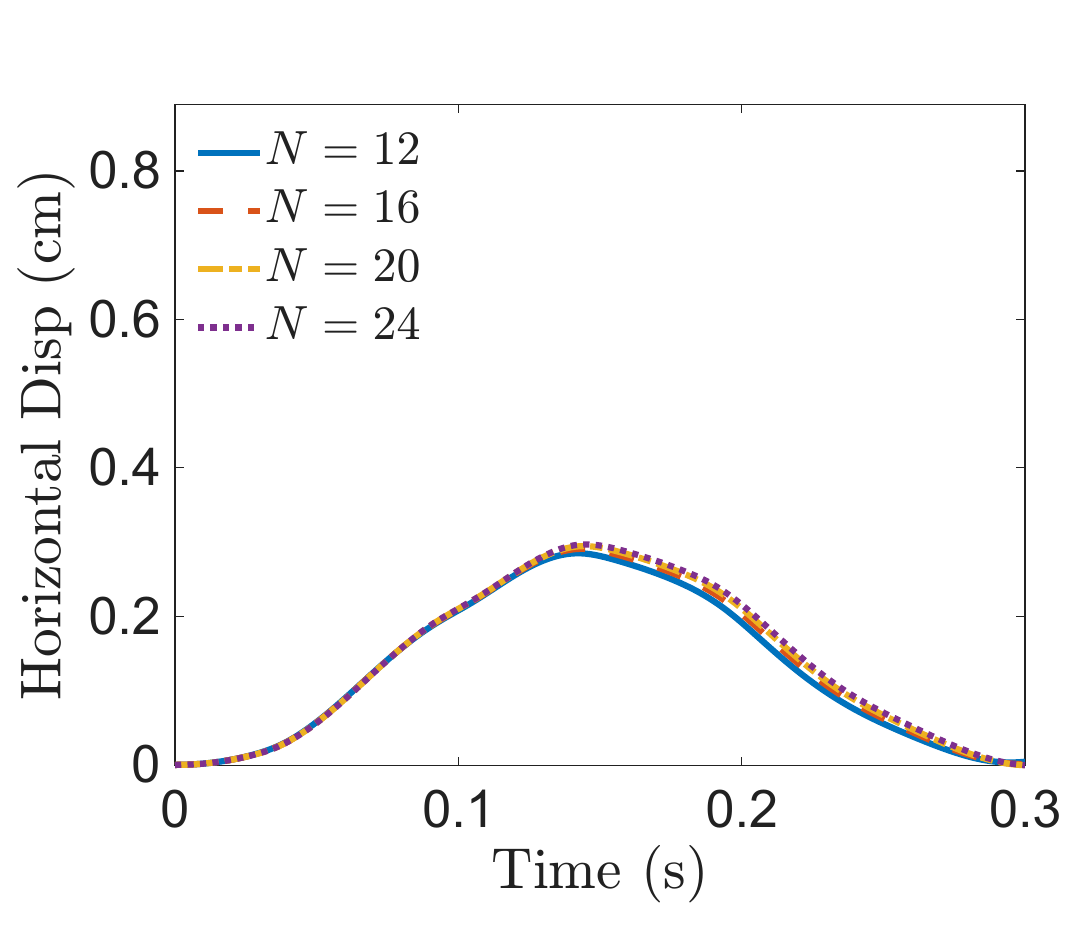}
\caption{Horizontal displacements of the point of interest, highlighted in Fig.~\ref{f:ELS_BND_schematics}, under grid refinement. The fiber direction is set to $\theta_1 = 90^\circ$. $N=12$ corresponds to $1261$ solid DoF, $N=16$ corresponds to $2193$ solid DoF, $N=20$ corresponds to $3381$ solid DoF, and $N=24$ corresponds to $4825$ solid DoF.}          		 \label{f:ELS_BND_disp}
\end{figure}

\begin{figure}[t!]
\centering
\includegraphics[width=\textwidth]{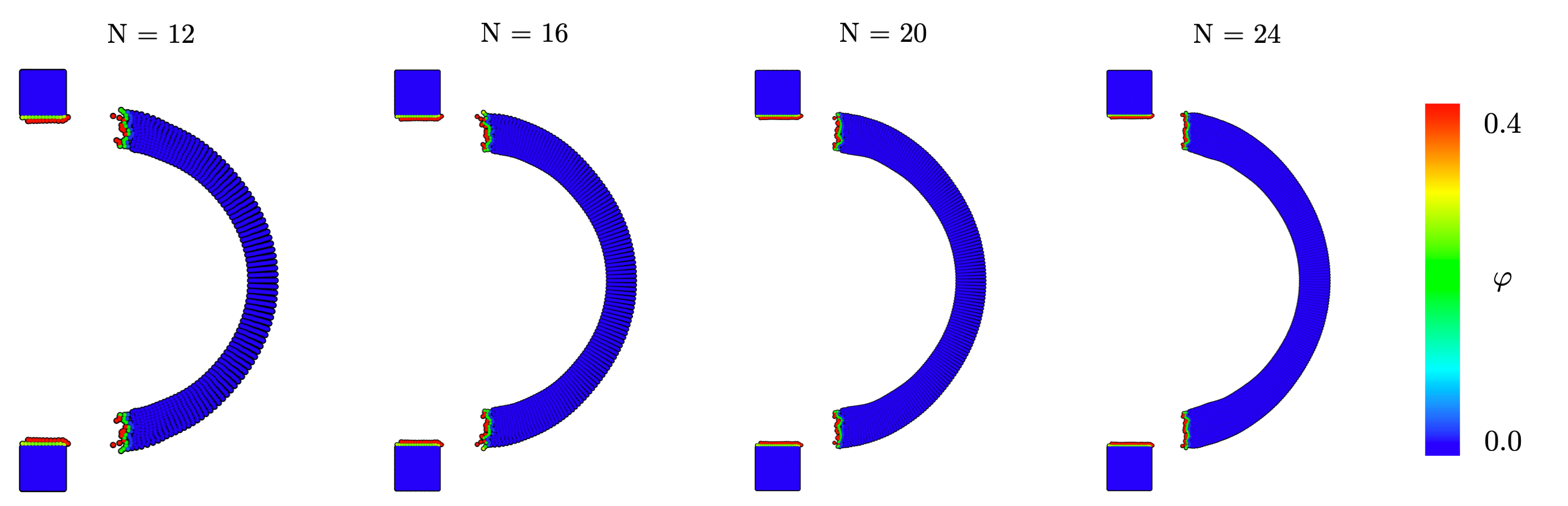}
\caption{Grid-refinement comparison of the crack trajectory in the two-dimensional elastic band shown in the deformed configuration with the local damage. For each resolution, the snapshot is taken at the instant when the band becomes completely detached from the fixed blocks. The fiber direction is $\theta_2 = 0^\circ$. $N=12$ corresponds to $1261$ solid DoF, $N=16$ corresponds to $2193$ solid DoF, $N=20$ corresponds to $3381$ solid DoF, and $N=24$ corresponds to $4825$ solid DoF.}          		 
\label{f:ELS_BND_crack_convergence}
\end{figure}

\begin{figure}[t!]
\centering
   \begin{tabular}{cc}
        \begin{subfigure}{.4\textwidth}
          		\includegraphics[width=\textwidth]{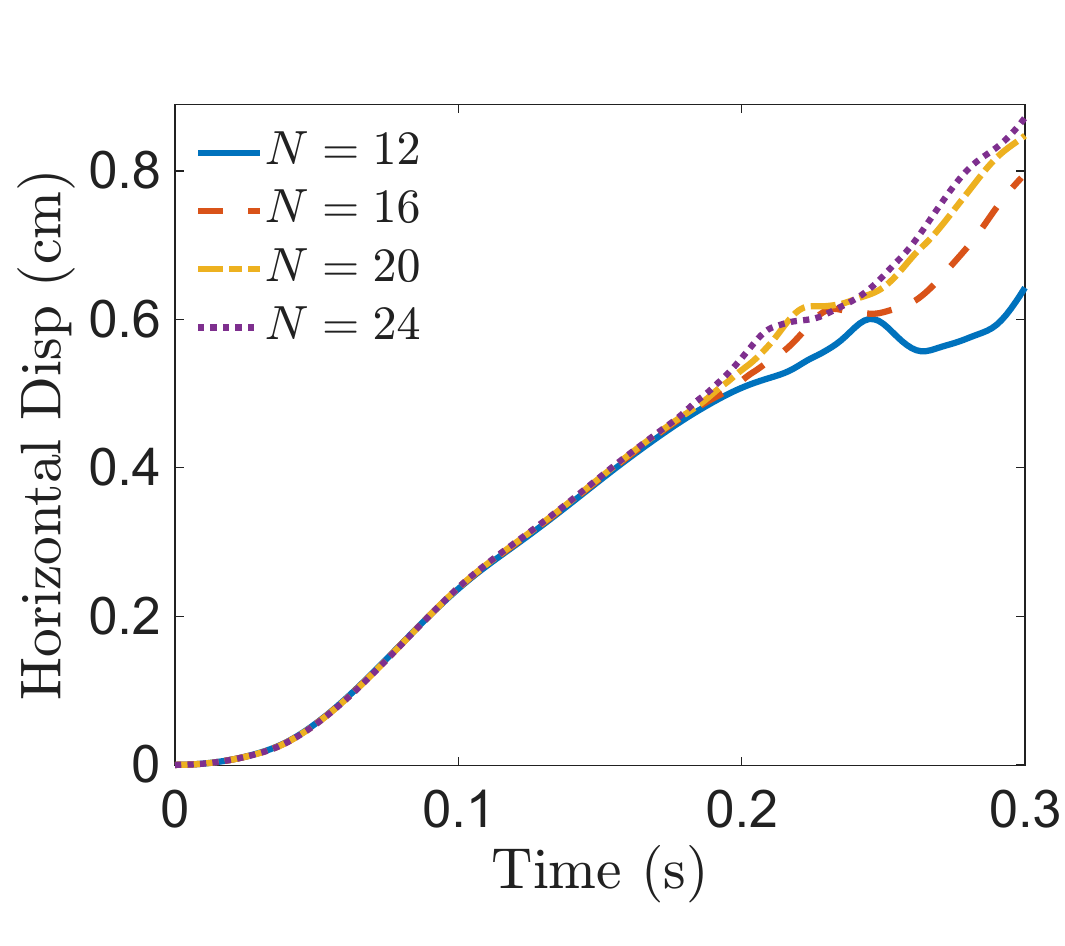}
          		\caption{}
          		 \label{f:ELS_BND_rupture_disp1}
        \end{subfigure} \hspace{.05\textwidth}
         \begin{subfigure}{.4\textwidth}
        			\includegraphics[width=\textwidth]{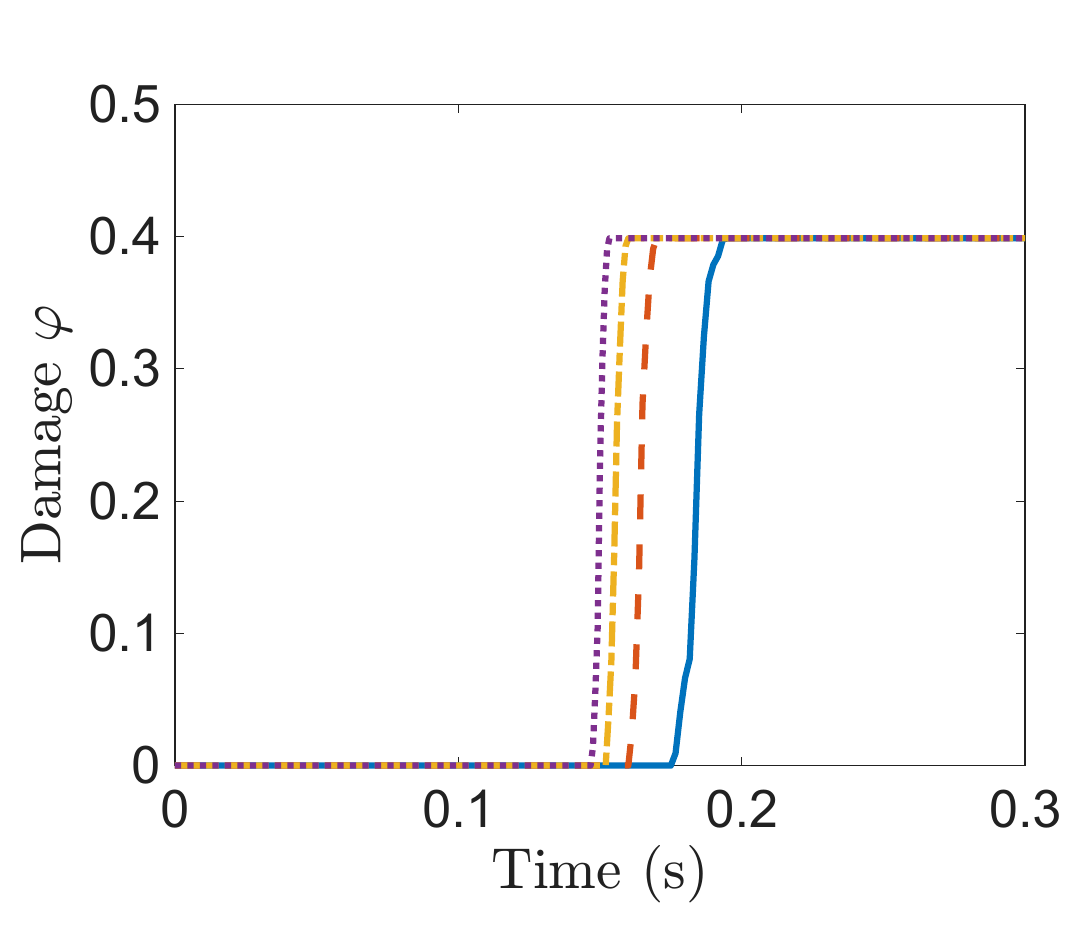}
          		\caption{}
          		 \label{f:ELS_BND_rupture_damage}
        \end{subfigure} 
    \end{tabular}
    \caption{(a) Horizontal displacements of the point of interest, highlighted in Fig.~\ref{f:ELS_BND_schematics}, under grid refinement. (b) Local damage growth at the top left corner of the detached band during the failure process under grid refinement. The fiber direction is set to $\theta_2 = 0^\circ$. $N=12$ corresponds to $1261$ solid DoF, $N=16$ corresponds to $2193$ solid DoF, $N=20$ corresponds to $3381$ solid DoF, and $N=24$ corresponds to $4825$ solid DoF.}
\end{figure}

Fig.~\ref{f:ELS_BND_deformations} shows the deformation of the elastic band with different choices of the fiber direction.
If the fiber is parallel to the fluid flow, the crack formation is initiated near the junctions between the fixed blocks and the band, and the band becomes entirely detached from the blocks when the bonds exceed the critical bond stretch.
On the other hand, the structural deformation is more resistant to the fluid flow if the fiber is perpendicular to the flow direction, and the immersed structure does not experience any structural damage or failure.
Fig.~\ref{f:ELS_BND_disp} and Fig.~\ref{f:ELS_BND_rupture_disp1} show the horizontal displacements of the point of interest, highlighted in Fig.~\ref{f:ELS_BND_schematics}, for different grid spacings. 
With the fiber oriented parallel to the fluid flow, we further demonstrate grid consistency by comparing the crack trajectories and local damage distributions at the instant when the band first becomes completely detached from the fixed blocks; see Fig.~\ref{f:ELS_BND_crack_convergence}.
Fig.~\ref{f:ELS_BND_rupture_damage} shows the local damage growth at the top left corner of the detached band during the failure process under grid refinement.

\section{Conclusions}
This paper develops an improved immersed peridynamics method to simulate fluid-driven material deformation, damage, and failure in anisotropic hyperelastic materials.
The initial version of our IPD formulation is limited to a uniform distribution of PD points, which cannot accurately represent a complex structural geometry and results in a stair-step geometry, as shown in Sec.~\ref{s:2d_cooks_discretizations}.
To accurately model anisotropic materials under FSI, this extended formulation integrates non-uniform structural discretizations, anisotropic constitutive models, and a ductile failure mechanism into the IPD framework.
We modify the volume terms in the discrete IPD formulation, as in Eqs.~\eqref{discrete_nonlocal_deformation_gradient}--\eqref{discrete_pd_force}, to achieve accurate PD volumetric force densities. 
Specifically, we adopt FE mesh generation techniques to accurately mesh the structural body and evaluate nodal volumes occupied by PD material points in the non-uniform discretizations of the immersed structure.
Such modifications will ultimately enable the IPD method to simulate realistic material behaviors in more complex FSI problems and further facilitate seamless coupling with standard FE-based methods for large-scale problems.
Our extended IPD results in Sec.~\ref{s:2d_cooks_discretizations} demonstrate that the non-uniformly distributed volumes of PD points yield improved deformational accuracy and volume conservation for classical solid mechanics benchmarks compared with the uniform lattice under grid refinement.

Further, we incorporate fiber-reinforced material models to capture deformations of anisotropic materials in Secs.~\ref{s:aniso_Benchmark_Compression}--\ref{s:aniso_Benchmark_Torsion}, and our IPD simulation results demonstrate $\horizonsize$-convergence to the FE-based IB simulations and achieve comparable accuracy and volume conservation with similar numbers of structural DoF.
The effect of the size of the peridynamic horizon is also investigated for anisotropy, and the IPD results are relatively insensitive to the size of the peridynamic horizons for non-failure tests. 
In Sec.~\ref{s:Aortic_tissue}, we adopt a material model for the adventitial tissue of an artery, demonstrating that the IPD method is highly capable of modeling actual biomaterials with an accuracy comparable to the classical FE-based method.
We examine damage and failure in anisotropic materials with a ductile failure criterion in Secs.~\ref{s:aniso_failure}--\ref{s:ansio_els_bnd} and demonstrate that the crack propagation paths are strongly dictated by the internal material structure, in contrast to isotropic behaviors, specifically aligned with the direction of the reinforcing fibers.
Although the present failure model captures the influence of constitutive anisotropy on local damage growth and crack
propagation, it does not explicitly represent direction-dependent fracture toughness. 
Future work will incorporate orientation-dependent or energy-based failure criteria calibrated using tissue-specific fracture data.

Correspondence-based PD formulations are known to suffer from zero-energy modes associated with the non-local deformation gradient tensor~\cite{silling2017stability}. 
In the numerical tests considered in this study, we do not observe spurious oscillatory deformation patterns characteristic of the zero-energy instabilities for the selected volumetric stabilization, IB kernel, mesh factor ratio, and PD horizon. 
This observation is limited to the numerical settings and benchmark problems examined here and neither demonstrates that zero-energy modes are mathematically eliminated nor identifies a general stabilization mechanism.
Future work will assess the stability of the IPD formulation across different discretizations, numerical parameters, and IB kernel supports to characterize the conditions for mitigating zero-energy modes and establish corresponding stability bounds.

Although the extended IPD framework demonstrates significant capabilities for simulating fluid-driven material failure, several methodological and numerical advancements are required to further elevate the mathematical robustness and computational efficiency of the framework for future applications in realistic fluid-driven fracture mechanics, such as aortic dissection.
A crucial next step involves leveraging the FE-based structural representation introduced in this work to facilitate seamless kinematic and dynamic coupling with FE-based methods. 
Such a hybrid coupling approach would allow bulk deformations to be resolved using highly efficient FE-based formulations, while reserving the IPD method strictly for computationally intensive localized damage and fracture. 
Such a numerical advancement is essential to significantly reduce the overall computational cost and enable the simulation of large-scale physiological systems.

Beyond computational efficiency, the physical fidelity of fracture mechanics in the IPD framework remains a primary area for future refinement. 
Although the current critical stretch-based failure criterion provides straightforward kinematics for bond breakage, adopting a mathematically rigorous energy-based failure formulation~\cite{foster2011energy} will better capture the complex energy dissipation mechanisms and crack tip singularities.
Furthermore, addressing the hydrodynamic resolution at the fluid-structure interface remains a critical numerical challenge. 
The current IB approach relies on regularized delta functions, which inherently introduce numerical smearing of the velocity and pressure fields across physical boundaries. 
Integrating sharp-interface techniques, such as an immersed interface method (IIM)~\cite{LeVeque1994,LI2001822}, will enforce exact stress jump conditions across the boundaries, providing the high-order interfacial accuracy required to properly resolve the highly localized and complex fluid dynamics within narrow, propagating crack surfaces. 
Crucially, such integration must address the inherent mismatch between the local nature of the IIM and the non-local PD kinematics, which requires mathematically consistent mapping strategies to distribute sharp surface tractions across the peridynamic horizon to mitigate artificial surface effects.
Addressing these numerical and mathematical aspects will solidify the IPD model as a highly robust computational FSI framework for advanced FSI fracture modeling.

\section*{Acknowledgements}
We gratefully acknowledge research support through NIH Award HL 117063 and NSF Awards OAC 1450327, OAC 1652541, OAC 1931516, and DMS 1929298.
Keon Ho Kim is grateful for research support through the American Heart Association Postdoctoral Fellowship 25POST1366103.
We thank Amneet P. S. Bhalla for the initial IPD implementation in IBAMR.
Numerical simulations were performed using facilities provided by the University of North Carolina at Chapel Hill through the Research Computing division of UNC Information Technology Services.

\end{document}